\documentclass[reqno]{amsart}
\theoremstyle{plain}
\usepackage{latexsym}
\usepackage{amssymb, color}
\usepackage{enumerate}
\usepackage{array}
\usepackage{epsfig}
\usepackage[sc]{mathpazo}
\definecolor{orange}{rgb}{1,0.5,0}
\definecolor{orange2}{rgb}{.7,.2,.3}

\def\XXint#1#2#3{{\setbox0=\hbox{$#1{#2#3}{\int}$}
     \vcenter{\hbox{$#2#3$}}\kern-.5\wd0}}

\newcommand{\PI}{\mathcal{P}}

\newcommand{\Mb}{\mathbb{M}}

\newcommand{\Sbb}{\mathbb{S}}

\newcommand{\Xb}{\mathbb{X}}

\newcommand{\Ib}{\mathbb{I}}

\newcommand{\lt}{\left}
\newcommand{\rt}{\right}
\newcommand{\nl}{\newline}

\newcommand{\nn}{\nonumber}

\newcommand{\lm}{\lambda}

\newcommand{\qd}{\quad}
\newcommand{\spt}{\mathrm{Spt}}

\newcommand{\ep}{\epsilon}

\newcommand{\wt}{\widetilde}

\newcommand{\AI}{\mathcal{A}}

\newcommand{\FI}{\mathcal{F}}
\newcommand{\MI}{\mathcal{M}}

\newcommand{\UI}{\mathcal{U}}

\newcommand{\OI}{\mathcal{O}}

\newcommand{\LI}{\mathcal{L}}

\newcommand{\ti}{\tilde}

\newcommand{\BI}{\mathcal{B}}

\newcommand{\R}{\mathrm {I\!R}}
\newcommand{\K}{\mathcal{K}}

\newcommand{\ca}[1]{\mathrm{Card}\lt(#1\rt)}
\newcommand{\dia}{\diamondsuit}

\newcommand{\red}[1]{{\textcolor{black}{#1}}} 
\newcommand{\blue}[1]{{\textcolor{black}{#1}}} 
\newcommand{\rred}[1]{{\textcolor{black}{#1}}} 
\newcommand{\bblue}[1]{{\textcolor{black}{#1}}} 
 
\newcommand{\ggreen}[1]{{\textcolor{black}{#1}}}
\newcommand{\xblue}[1]{{\textcolor{black}{#1}}} 
\newcommand{\xred}[1]{{\textcolor{black}{#1}}} 
\newcommand{\xxred}[1]{{\textcolor{black}{#1}}} 
 
\newcommand{\xxblue}[1]{{\textcolor{black}{#1}}} 
\newcommand{\zblue}[1]{{\textcolor{black}{#1}}}
\newcommand{\zzblue}[1]{{\textcolor{black}{#1}}}
\newcommand{\zzred}[1]{{\textcolor{black}{#1}}}
\newcommand{\ared}[1]{{\textcolor{black}{#1}}}
\newcommand{\ablue}[1]{{\textcolor{black}{#1}}}
\newcommand{\cblue}[1]{{\textcolor{black}{#1}}}
\newcommand{\xgreen}[1]{{\textcolor{black}{#1}}}
\newcommand{\pblue}[1]{{\textcolor{black}{#1}}}
\newcommand{\pg}[1]{{\textcolor{black}{#1}}}
\newcommand{\org}[1]{{\textcolor{black}{#1}}}
\newcommand{\orgg}[1]{{\textcolor{black}{#1}}}
\newcommand{\pgr}[1]{{\textcolor{black}{#1}}}
\newcommand{\fred}[1]{{\textcolor{black}{#1}}} 
\newcommand{\pb}[1]{{\textcolor{black}{#1}}}
\newcommand{\Nb}{\mathbb{N}}
\newcommand{\ap}{\alpha}
\newcommand{\ol}{\overline}

\newtheorem{a1}{Lemma}
\newtheorem{a2}[a1]{Theorem}

\newtheorem{a5}[a1]{Proposition}

\newtheorem{deff}[a1]{Definition}

\theoremstyle{remark}
\newtheorem{remark}{Remark}

\begin{document}
\title[Null Lagrangian Measures in \pg{subspaces}, compensated compactness and conservation laws]{Null Lagrangian Measures in \pg{subspaces}, compensated compactness and conservation laws}

\author[A. Lorent, G. Peng]{Andrew Lorent, Guanying Peng}
\address{A. L.\\Mathematics Department\\University of Cincinnati\\2600 Clifton Ave.\\ Cincinnati, OH 45221; 
G. P.\\Department of Mathematics\\University of Arizona\\617 N. Santa Rita Ave.\\ Tucson, AZ 85721.
}
\email{lorentaw@uc.edu, gypeng@math.arizona.edu.}
\subjclass[2000]{49J99, 28A25, 35L65}
\keywords{Null Lagrangian measures, compensated compactness, conservation laws, entropy solutions}
\maketitle

\begin{abstract}

Compensated compactness is \bblue{an important method} used to solve nonlinear PDEs, \zzred{in particular in the} study of \blue{hyperbolic} conservation laws. One of the simplest formulations of a compensated compactness problem is to ask for conditions on a \zzred{compact} set $\mathcal{K}\subset M^{m\times n}$ such that 
\begin{equation}
\label{abstreq1}
\zzblue{\lim_{j\rightarrow \infty} \|\mathrm{dist}(Du_j,\mathcal{K})\|_{L^p}= 0 \text{ and } \sup_{j}\|u_j\|_{W^{1,p}}<\infty\; \Rightarrow \{Du_{j}\}_{j}\text{ is precompact in } L^p}.
\end{equation}
Let $M_1,M_2,\dots, M_q$ denote the set of all minors of $M^{m\times n}$. A \xgreen{sufficient} condition for (\ref{abstreq1}) is that any \zzblue{probability} measure $\mu$ supported on $\mathcal{K}$ satisfying 
\begin{equation}
\label{abstreq2}
\int M_k(X) d\mu (X)=M_k\lt(\int X d\mu (X)\rt)\text{ for all } k
\end{equation}
is a Dirac measure. We call measures that satisfy (\ref{abstreq2}) \em Null Lagrangian Measures \rm and following \cite{mul}, we denote \blue{the set of} Null Lagrangian Measures supported on $\mathcal{K}$ by $\mathcal{M}^{pc}(\K)$. For general $m,n$, a necessary and sufficient condition for triviality of $\mathcal{M}^{pc}(\K)$ \org{was an open question} even in the case where $\K$ is a linear subspace \blue{of $M^{m\times n}$}. \org{We answer this question and} provide a necessary and sufficient condition for any linear subspace $\K\subset M^{m\times n}$. \pg{The ideas also allow us to show} that \red{for any 
$\pg{d}\in \lt\{1,2,3\rt\}$, $\pg{d}$-dimensional subspaces  $\K\subset M^{m\times n}$ support non-trivial Null Lagrangian Measures if and only if $\K$ has Rank-$1$ connections. This is known to be false for $\pgr{d\ge 4}$ \pg{from} \cite{bhat}.}

\red{Further }using the ideas developed we are able to answer  a question of Kirch\bblue{h}eim, M\"{u}ller and \v{S}ver\'{a}k \cite{kms}. Let 
$P_1(u,v):=\lt(\begin{matrix}
u & v\\
a(v) & u\\
ua(v) & \frac{1}{2}u^2+ F(v)\\
\end{matrix}\rt)$ \pgr{and} $\K_1 := \lt\{ P_1(u,v): u,v\in\R\rt\}$ for some function $a$ and its primitive $F$. The set $\K_1$ arises in the study of entropy solutions to the $2\times 2$ system of conservation laws 
\begin{equation*}
u_t=a(v)_x \qd\text{ and }\qd v_t=u_x.
\end{equation*}
In \cite{kms}, the authors asked what are the conditions on the function $a$ such that $\mathcal{M}^{pc}(\K_1\cap U)$ consists of Dirac measures, where 
$U$ is an open neighborhood of an arbitrary matrix in $\K_1$. Given $\alpha=(\alpha_1,\alpha_2)\in \R^2$, if $a'(\alpha_2)>0$ then we construct non-trivial measures in 
$\mathcal{M}^{pc}(\K_1\cap B_{\delta}\lt(P_1(\alpha)\rt))$ for any $\delta>0$. On the other hand if $a'(\alpha_2)<0$ then for 
sufficiently small $\delta>0$,  we show that $\mathcal{M}^{pc}(\K_1\cap B_{\delta}\lt(P_1(\alpha)\rt))$  consists of Dirac measures.

\end{abstract}

\section{Introduction}

Compensated compactness \zzred{(coupled with a-priori $L^{p}$ bounds) is an important method} of solving nonlinear PDEs. Amongst its most celebrated successes are the proofs of the first existence theorems for solutions of systems of hyperbolic conservation laws \bblue{with large data} by Tartar \cite{ta1}, \cite{ta2} and DiPerna \cite{dp2}, \cite{dp}. One of the simplest and most natural formulations of compensated compactness is to ask for conditions on a \zzblue{compact} set of matrices $\K\subset M^{m\times n}$ such that for any sequence $\{u_j\}\subset W^{1,p}(\Omega;\blue{\R^m})$, \zzblue{$1\leq p<\infty$,} defined on a \zzblue{bounded} domain $\Omega\subset\R^n$, if 
\begin{equation}
\label{eccz1}
\zzblue{\lim_{j\rightarrow \infty} \|\mathrm{dist}(Du_j,\mathcal{K})\|_{L^p(\Omega)}= 0  \qd\text{and}\qd u_j\rightharpoonup u\text{ in } W^{1,p}(\Omega;\blue{\R^m}) \text{ as }j\rightarrow \infty,}
\end{equation}
then \zzblue{there exists a subsequence such that}
\begin{equation}
\label{eccz2}
\zzblue{Du_{j_k}\rightarrow{Du} \text{ in } L^{p}(\Omega;\R^m) \text{ as }k\rightarrow \infty.}
\end{equation}
It turns out that \zzblue{a necessary and sufficient condition} on \zzred{a compact set} $\K$ for hypothesis (\ref{eccz1}) to imply (\ref{eccz2}) \red{is the following}: for any probability measure $\mu$ with $\spt\mu\subset \K$, if
\begin{equation*}
\int f(X) d\mu(X)\zzblue{\geq}f\lt(\int X d\mu(X) \rt) \text{for all Quasiconvex functions } f,
\end{equation*}
then $\mu$ is a Dirac measure. \zzred{Firstly note \ablue{that} by Corollary 3 \ablue{in} \cite{mulzhang}, since $\ablue{\K}$ is compact, we can without loss of 
generality assume $\lt\{u_j\rt\}\subset W^{1,\infty}(\Omega;\R^m)$. Then }this follows from Theorem \ablue{4.7} in \cite{mul} (see \zzblue{\cite{kind,kind2}} for the original source) and the fundamental theorem of 
Young measures (Theorem 3.1 and Corollary 3.2 in \cite{mul}).  However Quasiconvex functions are very hard to understand 
\footnote{Indeed one of the most important problems in Calculus of Variations is the question of whether in $2\times 2$ matrices, Rank-$1$ convex functions are Quasiconvex \cite{ball2}, \cite{astala}.}, so more commonly a smaller class of functions known as \em Polyconvex \rm functions are considered. These functions were introduced by Ball \cite{ball} in \rred{his} fundamental work on existence of minimizers of elasticity 
\rred{functionals}.  Given $X\in M^{m\times n}$, let $\hat{X}$ denote the vector of all minors of $X$. A polyconvex function is a function 
$f:M^{m\times n}\rightarrow \R$ that can be written as $f(X)=g(\hat{X})$ where $g$ is convex. 

Following \cite{mul}, \blue{given $\K\subset M^{m\times n}$}, we denote 
\begin{equation*}
\MI^{pc}(\K):=\lt\{\begin{array}{cc}\nu\in \mathcal{P}(M^{m\times n}):& \mathrm{Spt}(\nu)\subset \K, \int f(X) d\nu(X)\pb{\geq}f\lt(\overline{X}\rt) \text{ for all}\\
& \text{polyconvex functions } f, \text{where } \overline{X}=\int X d\nu(X)\end{array}\rt\}.
\end{equation*}
A function $g:M^{m\times n}\rightarrow \R$ is a \em Null Lagrangian \rm if $g(X)$ is an affine combination of the minors of $X\in M^{m\times n}$. Clearly 
if $g$ is a \em Null Lagrangian \rm then both $g$ and $-g$ are polyconvex. Therefore $\mu\in \MI^{pc}(\K)$ if and only if 
\begin{equation*}
\int M(X) d\mu(X)=M\lt(\int X d\mu(X)\rt)\text{ for all minors } M.
\end{equation*}
For this reason we shall call measures $\mu\in \MI^{pc}(\K)$ \em Null Lagrangian Measures\rm. 

As we will briefly sketch, the heart of a number of well known compensated compactness results is a proof that for some 
submanifold $\K$ in the space of matrices, $\mathcal{M}^{pc}(\K)$ consists of Dirac measures. There is overall little understanding of 
what general conditions  a set $\K$ has to have in order for $\mathcal{M}^{pc}(\K)$ to \xxblue{consist} of Dirac measures only, i.e., to be trivial. Even in the case when $\K$ is a linear subspace in the space of matrices, it \pg{was} an open problem to determine necessary and sufficient conditions 
on $\K$ for $\mathcal{M}^{pc}(\K)$ to be trivial \footnote{This was asked to the first author by V.  \v{S}ver\'{a}k  during a brief sabbatical visit to Minnesota in 2016.}. \red{\pblue{Our} Theorem \ref{C1} answers this question}. First we require some \pblue{definition}. 
\begin{deff}
\label{def1}
A set $\pg{S}\subset \R^n$ is a \emph{cone} \pblue{if} $\lm x\in \pg{S}$ \pblue{whenever} $x\in \pg{S}$ and $\lm>0$. A subset $V\subset\R^n$ is called a \emph{(real) algebraic set}  if $V$ is the locus of common zeros of some collection of polynomial functions on $\R^n$. An \emph{algebraic cone} \pblue{in $\R^n$} is a cone that is also an algebraic set.
\end{deff}
\begin{remark}
	\pblue{We say $\pg{V}\subset M^{m\times n}$ is an algebraic cone if $\pg{V}$ identified as a subset of $\R^{mn}$ is an algebraic cone.}
\end{remark}

\red{\begin{a2}
\label{C1}
Let $K\subset M^{m\times n}$ be a linear subspace.  Let $M_1,\dots, M_{\pg{q_1}}:M^{m\times n}\rightarrow \R$ be the set of all minors in $M^{m\times n}$ \pblue{and $M_{\pg{q_1}+1},\dots,M_{\pg{q_1}+mn}:M^{m\times n}\rightarrow \R$ be the projections onto the entries in $M^{m\times n}$}. Then $\mathcal{M}^{pc}(K)$ consists of Dirac measures if 
and only if 
\begin{eqnarray}
\label{sseqc2}
&~&\text{for each non-trivial algebraic cone  }\pg{V}\subset K\text{ there exists }\beta\in \pblue{\R^{\pg{q_1}\pblue{+mn}}\setminus\{0\}}\nn\\
&~&\qd\qd\text{ such that }\sum_{k=1}^{\pg{q_1}\pblue{+mn}} \beta_k M_k(\pblue{\zeta})\geq 0\text{ for all }\pblue{\zeta}\in \pg{V}\text{ and }\sum_{k=1}^{\pg{q_1}\pblue{+mn}} \beta_k M_k\not\equiv 0\text{ on }\pg{V}.
\end{eqnarray}
\end{a2}}
\pg{Our} \red{Theorem \ref{C1} is actually a special case of a more general result for measures that commute with a class of 
	homogeneous polynomials. Since the statement \pg{of the more general theorem} requires more background notations we postpone it until Section \ref{SS3}} (\pg{see Theorem \ref{T1}}). 

We say that a set $\xblue{\Sigma}\subset M^{m\times n}$ has \em Rank-$1$ connections \rm if and only if there exist $A,B\in \xblue{\Sigma}$ such that 
	$A\not=B$ and $\mathrm{Rank}(A-B)=1$. Note that if $K\subset M^{m\times n}$ is a subspace that satisfies (\ref{sseqc2}), then $K$ has no Rank-$1$ connections. Indeed, were this not the case, there would be a \pg{Rank-$1$ line $V\subset K$ which forms a non-trivial algebraic cone in $K$ such that} $M_k(A)=0$ for all $k=1,2,\dots, \pg{q_1}$. \pg{However every linear combination of the projection mappings $M_{q_1+1},\dots,M_{q_1+mn}$ either is trivial or changes sign on $V$,} which contradicts condition (\ref{sseqc2}). \pg{Note that in \cite{sv2},  \v{S}ver\'{a}k proved the \zblue{beautiful} result that for connected sets 
		$\K\subset M^{2\times 2}$, $\mathcal{M}^{pc}(\K)$ is trivial if and only if $\K$ does not 
		contain Rank-$1$ connections.} \pg{So a natural question is whether condition (\ref{sseqc2}), and thus triviality of $\mathcal{M}^{pc}(K)$, is equivalent to $K$ having no Rank-$1$ connections for subspaces $K$.} We have the following:
\red{\begin{a2}
\label{CC1} 
Let $\pblue{d}\in \lt\{1,2,3\rt\}$ and $K\subset M^{m\times n}$ be \pblue{a $d$-dimensional subspace}, then 
$\MI^{pc}(K)$ consists of Dirac measures if and only if $K$ does not contain \pblue{Rank-$1$} connections. 
\end{a2}}
 \pg{Such equivalence relation} is false even for subspaces $K\subset M^{m\times n}$ with $\orgg{\dim(K)\geq 4}$ (see the Appendix \ref{S11.2} \pblue{for a counter example given in \cite{bhat}}). Thus Theorem \pblue{\ref{CC1}} \orgg{is} optimal. \pg{Note that} for the applications that we 
have developed in this paper (and an application in a previous preprint version \cite{lppv2}, Section 9), condition (\ref{sseqc2}) is actually more useful and informative.

\zzred{Note that the set of $k$-dimensional subspaces in $M^{m\times n}$ is essentially the \xgreen{Grassmannian} space $G(k,mn)$. As is 
well known, $G(k,mn)$ forms \cblue{a} $k(mn-k)$-\cblue{dimensional} smooth compact connected manifold, see \cite{mas} or Section \cblue{\ref{t4}}, Lemma \ref{graslem1}. 
As such we say \cblue{that} a property holds  ``generically'' for $k$-\cblue{dimensional} subspaces in $M^{m\times n}$ if the set of subspaces for which \cblue{the property} does not 
hold can be covered by the Lipschitz images of a finite \cblue{collection} of submanifolds in $\R^{k(mn-k)}$ of dimension less than $k(mn-k)$, see 
Definition \ref{grasdef1}. \cblue{Using this point of view, we have the following} \ared{result}:
\begin{a2}
	\label{C17}
	Suppose \cblue{$k,m,n$ are positive integers with} $m, n\geq 2$ and $k \leq \frac{1}{2}\min\lt\{m, n\rt\}$. Then for a ``generic'' $k$-dimensional 
subspace $K\subset M^{m\times n}$, $\mathcal{M}^{pc}(K)$ consists of Dirac measures and hence $K$ has no Rank-$1$ connections. 
\end{a2}
Contrast this with the interesting result of Bhattacharya, Firoozye, James and Kohn (Proposition 4.4 in \cite{bhat}), \cblue{in which it is shown that $l(m,n)\leq mn-n$, where $l(m,n)$ denotes the 
\cblue{maximum} possible dimension of a linear subspace in $M^{m\times n}$ that does not have Rank-$1$ connections}. So \cblue{Theorem} \ref{C17} is completely false for higher dimensional subspaces. The bound $k \leq \frac{1}{2}\min\lt\{m,n\rt\}$ is surely not sharp, \cblue{and} an interesting 
and possibly accessible question is to determine the sharp bound \cblue{on $k$ such that the conclusion of Theorem \ref{C17} holds true. Although Theorem \ref{C17} is not a consequence of} 
Theorem \ref{C1}, the ideas of its proof \cblue{are} very closely related to \cblue{those of} the proof of Theorem \ref{C1}.}\nl

One of the motivations for studying Null Lagrangian Measures 
supported on \pgr{subspaces} is that such results might rather directly yield insights into how to prove triviality \org{or non-triviality} of $\mathcal{M}^{pc}(\K\cap U)$ where 
$\K\pgr{\subset M^{m\times n}}$ is any smooth submanifold and $U$ is a small neighborhood around an arbitrary \pgr{point $\zeta\in\K$}. \org{Since $\K\cap U$ \pgr{can be} arbitrarily well approximated by its tangent 
plane at $\zeta$, we might expect condition (\ref{sseqc2}) to be relevant in \pgr{understanding} the structure of $\mathcal{M}^{pc}(\K\cap U)$ and indeed this turns out to be the case}. In the following subsection we apply these insights to study a well known $2\times 2$ system of conservation laws. \org{As we will outline, the study of the weak solutions of the system that arise via compensated compactness is intimately connected with the structure of $\mathcal{M}^{pc}(\K)$ for a smooth submanifold $\K$ in matrix space.  In the case where the system is adjoined by a single additional  entropy inequality, \pgr{it is a model problem for systems of conservation laws in higher dimensions} and the related set $\K_1\subset M^{3\times 2}$ has \pgr{numerous open questions} about its structure (\orgg{\cite{kms}, Section 7}). \pgr{U}sing the ideas developed in the proof of Theorem \ref{C1}, we answer the question \fred{of} the 
structure  of $\mathcal{M}^{pc}(\K_1)$.}

\subsection{Connections and applications to conservation laws}\label{s1.2}

As mentioned above one of the main \zzred{successes} of compensated compactness is the proof of existence theorems for \blue{hyperbolic} conservation laws. To sketch this briefly, the standard way to solve a scalar equation is to add a viscosity term and obtain a solution 
to 
\begin{equation}
\label{eccz20}
u_t^{\ep}+\bblue{G}(u^{\ep})_x=\ep u^{\ep}_{xx}\text{ in }\xxblue{(0,\infty)\times\R}.
\end{equation}
Assuming $\lt\{u^{\ep}\rt\}_{\ep}$ is bounded in $L^{\infty}(\xxblue{(0,\infty)\times\R})$ we can extract a 
subsequence $u^{\ep_k}\overset{*}{\rightharpoonup} u$ in $L^{\infty}(\xxblue{(0,\infty)\times\R})$. Letting $\nu_{t,x}$ be the Young measure \blue{associated with the weak* convergence}, i.e., $u(t,x)=\int_{\R} y \;d\nu_{t,x}$, we have $G(u^{\ep_k})\overset{*}{\rightharpoonup} \overline{G}$ in 
$L^{\infty}(\xxblue{(0,\infty)\times\R})$ where $\overline{G}(t,x)=\int G(y)  \;d\nu_{t,x}$. 

Now for any convex function $\Phi:\R\rightarrow \R$, define $\Psi(y):=\int_{\bblue{0}}^{y} \Phi'(s) G'(s) ds$. The pair 
$(\Phi,\Psi)$ is called an entropy/entropy flux pair. The key point is that by virtue of the Div-Curl lemma we know that 
\begin{equation}
\label{eccz21}
\int \lt(G(y)\Phi(y)-y \Psi(y)\rt) \;d\nu_{t,x}=\overline{G}(t,x)\overline{\Phi}(t,x)-u(t,x)\overline{\Psi}(t,x),
\end{equation}
where $\overline{\Phi}(t,x)=\int \Phi(y) d\nu_{t,x}$ and  $\overline{\Psi}(t,x)=\int \Psi(y) d\nu_{t,x}$. Define 
$P_{\Phi}:\R\rightarrow M^{2\times 2}$ by $P_{\Phi}(\bblue{z}):=\lt(\begin{array}{cc} G(z) & z \\ \Psi(z) & \Phi(z) \end{array}\rt)$ 
and the measure $\mu_{\Phi}$ on the set $\K_{\Phi}:=\lt\{P_{\Phi}(z):z\in \R\rt\}$ 
by \red{$\mu_{\Phi}:=(P_{\Phi})_{\sharp} \nu_{t,x}$, the push forward of $\nu_{t,x}$ by the \pg{mapping} $P_{\Phi}$}. By (\ref{eccz21}), $\mu_{\Phi}\in \mathcal{M}^{pc}(\K_{\Phi})$. So to prove triviality of $\nu_{t,x}$ \rred{it} suffices to prove $\mathcal{M}^{pc}(\K_{\Phi})$ is trivial 
for any choice of convex function $\Phi$. As this is such a wide class, \blue{for a lot of scalar conservation laws, one can find an appropriate convex function $\Phi$ for which $\mathcal{M}^{pc}(\K_{\Phi})$ is trivial, and hence the Young measures are trivial. \bblue{The fact that $u$ is a weak solution to (\ref{eccz20}) \xxblue{without viscosity term} follows from triviality of Young measures in a standard way.}}

For systems of conservation laws (\cblue{other than $2\times 2$ systems or scalar conservation laws}) there are only finitely many entropy/entropy flux pairs $(\Phi_1,\Psi_1)$,  
$(\Phi_2,\Psi_2)$, \dots, $(\Phi_m,\Psi_m)$. By analogous argument to the \xgreen{scalar} case, \xgreen{the Young measures can be pushed forward into $\mathcal{M}^{pc}(\K)$, where $\K$ is the subset of matrices whose rows consist of the conservation laws and the entropy/entropy flux pairs $(\Phi_j,\Psi_j)$}. Then triviality of the Young measures and hence \zzred{(given appropriate a-priori $L^p$ bounds)} proof of existence of solutions via compensated compactness comes down to proving \xgreen{triviality of} $\mathcal{M}^{pc}(\K)$. 

One of the best known 
results in this area is the work of DiPerna \cite{dp2} on strong convergence of solutions to a class of systems of two genuinely nonlinear conservation laws in one dimension, where the hypotheses are compactness in $W^{-1,2}$ of \em every \rm entropy/entropy flux pair acting on the approximating solutions. As a particular example, the result applies to the system with the form of the Lagrangian equations of elasticity \bblue{given by}
	\begin{equation}\label{eq201}
	\begin{cases}
	v_t - u_x=0,\\
	u_t - a(v)_x = 0
	\end{cases}
	\end{equation}
	for some given smooth function $a:\R\rightarrow \R$ that is strictly convex and increasing. Possibly motivated by the question of compactness for higher dimensional systems, in another well known work DiPerna \cite{dp} proves a local existence result for the system (\ref{eq201}) with just two entropy/entropy flux pairs. Following \cite{dp} we introduce the natural entropy/entropy flux pair $(\eta_1, q_1)$ associated to the system (\ref{eq201}). More precisely, we define
\begin{equation*}
\eta_1(u,v) := \frac{1}{2}u^2+F(v), \qd q_1(u,v) := -ua(v),
\end{equation*}
where $F(\xi)=\int_0^{\xi} a(s)ds$. As in \cite{kms}, we consider entropy solutions of (\ref{eq201}) defined as $L^{\infty}$ functions $(u,v)$ satisfying
\begin{equation}\label{eq203}
\begin{cases}
v_t - u_x=0,\\
u_t - a(v)_x = 0,\\
(\eta_1)_t + (q_1)_x \leq 0
\end{cases}
\end{equation}
in the sense of distributions. Adding a viscosity term to the first two equations of (\ref{eq203})  we obtain the pair $(u^{\ep},v^{\ep})$ that solves 
\begin{equation*}
v^{\ep}_t-u^{\ep}_x=\ep v^{\ep}_{xx}, \;\; u^{\ep}_t-a\lt(v^{\ep}\rt)_x=\ep u^{\ep}_{xx}.
\end{equation*}
Assuming appropriate bounds on \bblue{$u^{\ep}, v^{\ep}, u^{\ep}_x, v^{\ep}_x$} (see (5.38) of \cite{evans}), we obtain the system 
\blue{(\ref{eq203})} with right hand side precompact in $W^{-1,2}_{loc}$. Hence as we have sketched for scalar equations, we have $(u^{\ep},v^{\ep})\overset{*}{\rightharpoonup} (u,v)$ \bblue{in $L^{\infty}$} and the Young \pg{measures} can be pushed forward into the set $\K_1$ where 
\begin{equation}\label{eq205}
\K_1 := \lt\{\lt(\begin{matrix}
u & v\\
a(v) & u\\
ua(v) & \frac{1}{2}u^2+F(v)\\
\end{matrix}\rt): u,v\in\R\rt\}.
\end{equation}
By use of the Div-Curl lemma we have \pg{measures} in $\mathcal{M}^{pc}(\K_1)$. 

\pg{Thus understanding the structure} of $\K_1$ plays a fundamental role in understanding the system (\ref{eq203}). If there is so little rigidity of the structure of $\K_1$ \blue{ that certain subset $\K_1^{rc}$ of $\K_1^{pc}$ ($\K_1^{pc}$ and $\K_1^{rc}$ are called the Polyconvex hull and Rank-$1$ convex hull of $\K_1$, respectively, see Section 4.4 \xxblue{in} \cite{mul}) is sufficiently non-trivial}, then a very different kind of non-trivial solution to (\ref{eq203}) can be obtained as a differential inclusion into $\K_1$\footnote{ Note that a differential inclusion into set $\K_1$ gives a solution to (\ref{eq203}) with the inequality replaced by an equality.}. There have been enormous interests and spectacular progresses in reformulating PDEs as differential inclusions and obtaining solutions via \em convex integration \rm\cite{camles1}, \cite{camles2}, \cite{camles3}, \cite{phil}. Some of the initial impetus for these works come from 
the pioneering work on Calculus of Variations by \cite{mulsv1}, \cite{mulsv2}, \cite{sychev}, \cite{kirch9}, \cite{ki}. For this reason Kirchheim, M\"{u}ller and \v{S}ver\'{a}k  \cite{kms} asked the following question with respect to \xgreen{the} system (\ref{eq203}) and its associated differential inclusion \bblue{into the set $\K_1$} \xgreen{under more general conditions on the function $a$, namely,} what are the natural assumptions on the function $a$ such that the following statement is true:\nl 

\noindent\textbf{(S1)} For each point $\zeta\in \red{\K_1}$, there exists a neighborhood $U\subset M^{3\times 2}$ of $\zeta$ such that $\mathcal{M}^{pc}(\K_1\cap \overline{U})$ is trivial.\nl

\blue{For the system (\ref{eq201}) without implementing any entropy/entropy flux pairs, the statement \bf (S1) \rm \bblue{for the corresponding set $\K_0:=\lt\{\lt(\begin{smallmatrix}  u & v \\ a(v) & u\end{smallmatrix}\rt):u,v\in \R\rt\}$} is well understood using results in \bblue{\cite{sv2}}}. On the other hand, it is proved in \cite{dp} that a set \pblue{analogous to $\K_1$ obtained 
by inclusion of an additional dual entropy/entropy flux pair} satisfies statement \bf (S1) \rm \red{if the function $a$ has the properties $a'>0$ and $a''\ne 0$}. However, this question (as well as some other related properties) for the set $\K_1$ defined in (\ref{eq205}) (\ared{which is associated with }system (\ref{eq203}) with just one entropy/entropy flux pair) remained open. (For more details, see \cite{kms}, Section 7.) 

For the convenience of later discussions, we parametrize the set $\mathcal{K}_1$  by the mapping 
\begin{equation}
\label{auxxeq7.4}
P_1(u,v):=\lt(\begin{matrix}
u & v\\
a(v) & u\\
ua(v) & \frac{1}{2}u^2+ F(v)\\
\end{matrix}\rt).
\end{equation}
In this notation, $\mathcal{K}_1= \lt\{ P_1(u,v): u,v\in\R\rt\}$. In Section \ref{s7}, given a point $P_1\lt(\xgreen{(\ti \alpha_1,\ti\alpha_2)}\rt)\in\K_1$, we will show that statement \bf (S1) \rm is false if $a'(\ti\alpha_2)>0$ and true if $a'(\ti\alpha_2)<0$. Specifically, we have

\begin{a2}\label{t102}
	\blue{Suppose $a\in C^{2}(\R)$. Given $\bblue{\ti\alpha}\xgreen{=(\ti \alpha_1,\ti \alpha_2)}\in\R^2$, if $a'(\bblue{\ti\alpha}_2)>0$, then there exist non-trivial measures in $\mathcal{M}^{pc}\lt(\K_1\cap B_{\delta}(P_1(\bblue{\ti\alpha}))\rt)$ for all $\delta>0$. On the other hand, if $a'(\bblue{\ti\alpha}_2)<0$, then there exists $\delta_0>0$ depending on the function $a$ and $\bblue{\ti\alpha}_2$ such that $\mathcal{M}^{pc}\lt(\K_1\cap B_{\delta}(P_1(\ti\alpha))\rt)$ is trivial for all $0<\delta\leq\delta_0$. }
\end{a2}

\red{Indeed the second part of Theorem \ref{t102} can be made a bit stronger. More precisely, recall that 
$\K_0:=\lt\{\lt(\begin{smallmatrix}  u & v \\ a(v) & u\end{smallmatrix}\rt):u,v\in \R\rt\}$. Given $\bblue{\ti\alpha}\in\R^2$, if $a'(\bblue{\ti\alpha}_2)<0$ then $\mathcal{M}^{pc}\lt(\K_0\cap B_{\delta}\lt(\lt(\begin{smallmatrix} \ti{\alpha}_1 & \ti{\alpha}_2 \\ a(\ti{\alpha}_2) & \xxblue{\ti{\alpha}_1} \end{smallmatrix}\rt)\rt)\rt)$ is trivial for sufficiently small $\delta>0$ depending on $a$ and $\bblue{\ti\alpha}_2$. As $\mathcal{M}^{pc}(\K_1)$ can be naturally embedded into $\mathcal{M}^{pc}(\K_0)$, this implies the second part of the theorem} (see the proof of Theorem \ref{t102} in Section \ref{s7}).  Theorem \ref{t102} is closely related to Theorem \ref{C1}. Indeed, one can check directly that for the submanifold $\K_1$ given in (\ref{eq205}), there does not exist non-trivial linear combination of all three minors that remains non-negative.  Nevertheless, it should be noted that the set $\K_1$ given in (\ref{eq205}) is a nonlinear submanifold in the space of $3\times 2$ matrices whose nonlinear structure poses \blue{extremely delicate} issues. \blue{As a result, the arguments needed are \rred{significantly} beyond those used for \red{subspaces}. \pg{Our} proof of the first part in Theorem \ref{t102} \pg{is constructive and} allows to produce infinitely many \pg{non-trivial} elements in $\MI^{pc}(\K_1)$ (see Theorem \ref{t103}).}

\subsection{\xblue{Acknowledgments}} The first \xblue{author} is very grateful to V. \v{S}ver\'{a}k  for many very helpful conversations during a \xblue{two-week} visit to Minnesota in \xblue{November} of 2016. \xblue{T}hese conversations essentially introduced us to this topic and \xblue{led} to some \xblue{initial} ideas for Theorem \ref{C1}.  Both authors are very grateful for \ared{a great deal of} very helpful \xblue{correspondence} since then. \xred{We would also like to thank S. M\"{u}ller for providing us with a proof of Lemma \ref{LAUX33}. The proof provided is considerably simpler and more elegant than our original proof.} \pblue{The second author would like to thank Y. Shi for helpful discussions on elementary algebraic geometry.}

\section{A more general formulation of Theorem \ref{C1}} 
\label{SS3}

\pblue{In this section, we give a more general formulation of Theorem \ref{C1} in terms of homogeneous polynomials (see Theorem \ref{T1} below). To state the theorem, we need some preparation.} 

\pblue{Given a set $S\subset\R^n$}, let $\PI(S)$ denote the space of probability measures supported on $S$. Given a collection of homogeneous polynomials $\FI=\lt\{f_1,f_2, \dots, f_{M_0}\rt\}$ on $\R^n$, let $\pblue{\LI(\FI)}$ denote the 
set of linear combinations of the functions in $\FI$, i.e.,
\begin{equation}
\label{eqa1}
\LI(\FI):=\lt\{\sum_{i=1}^{M_0} \lm_i f_i+\lm_0:\pblue{f_i\in\FI}, \lm_0,\lm_1,\dots, \lm_{\pblue{M}_0}\in \R\rt\}.
\end{equation}

\begin{deff}
\label{def2}
We say \pblue{that a collection of homogeneous polynomials} $\FI$ \pblue{on $\R^n$} satisfies property $R$ if
\begin{equation*}
f(z-\pblue{z}_0)\in \LI(\FI)\text{ for any }f\in \FI \blue{\text{ and any } z_0\in\R^n}.
\end{equation*}
\end{deff}

\begin{deff}
\label{def3}
We define the set of Null Lagrangian Measures with respect to a set of \pblue{homogeneous} polynomials $\mathcal{F}$ \pblue{on $\R^n$}  by 
\begin{equation*}
\mathbb{M}^{pc}_{\mathcal{F}}:=\lt\{\mu\in \mathcal{P}(\R^n): f\lt( \int z \; d\mu(z)\rt)=\int f(z) \; d\mu(z)\;\text{ for all }f\in \mathcal{F} \rt\},
\end{equation*}
and further we define $\mathbb{M}^{pc}_{\mathcal{F}}(\varpi):=\lt\{\mu\in  \mathbb{M}^{pc}_{\mathcal{F}}: \int z\; d\mu(z)=\varpi \rt\}$.
\end{deff}

\pblue{Now we are ready to state the more general formulation of Theorem \ref{C1}. Recalling the definition of algebraic cone in Definition \ref{def1}, we have}

\begin{a2}
\label{T1}
Let $\FI=\lt\{f_1,f_2,\dots, f_{M_0}, \blue{f_{M_0+1}, \dots, f_{M_0+n}}\rt\}$ be a collection of homogeneous polynomials \pblue{on $\R^n$} satisfying $f_{M_0+j}(z)=z_j$ \pblue{for $j=1,\dots,n$} and the property $R$ \pblue{as in Definition \ref{def2}}. Then \pg{$\Mb^{pc}_{\FI}$ consists of Dirac measures} if and only if 
\begin{eqnarray}
\label{eqcc1}
&~&\text{for each \blue{non-trivial algebraic} cone  }\pg{V}\subset \R^{n}\text{ there exists }y\in \R^{M_0+n}\backslash\lt\{0\rt\}\nn\\
&~&\qd\qd\qd\qd\text{ such that }\sum_{k=1}^{M_0\blue{+n}} y_k f_k\geq 0\text{ and }\sum_{k=1}^{M_0\blue{+n}}  y_k f_k\not\equiv 0\text{ on }\blue{\pg{V}}. 
\end{eqnarray}
\end{a2}

\pblue{The reason why we are interested in homogeneous polynomials is clear: minors in $M^{m\times n}$ are simply homogeneous polynomials. In Section \ref{sec4}, our efforts will be devoted to proving Theorem \ref{T1}. As can be seen in Section \ref{sec5}, Theorem \ref{C1} is \pgr{a} \orgg{fairly straightforward} consequence of the above theorem.}

\section{Proof Sketch}

In this section we will sketch briefly the main ideas of the proofs of \pg{our main theorems}.

\subsection{Sketch of proofs of Theorems \ref{T1} and \ref{C1}} \label{s2.1} To illustrate the key ideas, we sketch the proof in the special case \pg{where $\varpi=0$}. Let $\mu\in \mathbb{M}^{pc}_{\FI}(0)$. \bblue{By definition}, we have that 
\begin{equation}
\label{zzzdde1}
\int f_k(z) d\mu(z) = f_k\lt(\int z\; d\mu(z)\rt)=0 \text{ for }k=1,2,\dots, M_0+n.
\end{equation}
If the condition (\ref{eqcc1}) is satisfied, then we can find some $\pg{y\in\R^{M_0+n}\setminus\{0\}}$ such that $g(z):=\sum_{k=1}^{M_0+n} y_k f_k\geq 0$ on $\R^n$. \pg{It is not hard to show that the highest degree terms in $g$, denoted by $g_1$ which is homogeneous, is also non-negative and non-trivial on $\R^n$}. \bblue{By (\ref{zzzdde1}), we have} $\int g_1(z) d\mu(z)=0$, and therefore $\spt\mu\subset \pg{V}:=\lt\{z:g_1(z)=0\rt\}$ \pb{and $V$} is an algebraic cone. \bblue{By assumption, we can find} another linear combination 
that is non-trivial and non-negative on $\pg{V}$. This way we can iteratively reduce the support of $\mu$ onto cones of smaller and smaller dimensions until the support is reduced to the origin.

The necessity part of the proof is a bit more intricate. Suppose there exists a \pg{cone} $\pg{V}$ such that 
\begin{equation}
\label{rreqdd1}
\sum_{k=1}^{M_0+n} y_k f_k\text{ \pg{changes sign} on }\pg{V}\text{ for every }y\in \R^{M_0+n}\setminus\{0\}. 
\end{equation}
To construct a non-trivial measure in $\mathbb{M}^{pc}_{\FI}(0)$ it suffices to find points $\zeta_1,\zeta_2,\dots, \zeta_{m_0}\in \R^n$ and \bblue{weights $\gamma_1,\gamma_2,\dots, \gamma_{m_0}\geq 0$} satisfying
\begin{equation}
\label{zzzdde3}
\sum_{l=1}^{m_0}\gamma_l f_k(\zeta_l) = 0\text{ for all }k=1,2,\dots,M_0+n\text{ and }\sum_{l=1}^{m_0}\gamma_l=1.
\end{equation}
Then defining $\mu:=\sum_{l=1}^{m_0} \gamma_l \delta_{\zeta_l}$ we have $\mu\in\mathbb{M}_{\mathcal{F}}^{pc}(0)$. Indeed, if we find solutions to (\ref{zzzdde3}), then simply because 
the set of functions $\FI$ contains the projections $f_{M_0+j}(z)=z_j$ we automatically have $\ol\mu=\sum_{l=1}^{m_0} \gamma_l \zeta_l=0$. \pg{Further the equations for $f_k$ for $k=1,\dots,M_0$ imply that $\mu$ commutes with these functions.}

Now define $a(\zeta):=\lt(f_1(\zeta),
\dots, f_{M_0+n}(\zeta)\rt)$, $\AI:=\lt\{a(\zeta):\zeta\in \pg{V}\backslash \lt\{0\rt\}\rt\}$ and 
$b=0\in\R^{M_0+n}$. Finding $\zeta_1,\zeta_2,\dots, \zeta_{m_0}$ and $\gamma_1,\gamma_2,\dots, \gamma_{m_0}\geq 0$ that 
satisfy (\ref{zzzdde3}) is equivalent to showing $b\in \mathrm{Conv}(\AI)$. Suppose this was false, then by the 
Hyperplane Separation Theorem we must be able to find some $c\in \R$ and $y\in \R^{M_0+n}$ such that 
$y\cdot w\geq c$ for all $w\in \mathrm{Conv}(\AI)$ and $y\cdot b\pg{\leq}c$. However for any such $y$, by (\ref{rreqdd1}) there must 
exist some $\zeta_y\in \pg{V}\backslash \lt\{0\rt\}$ such that $\sum_{k=1}^{M_0+n} y_kf_k(\zeta_y)=a(\zeta_y)\cdot y<0$, which implies that $c\leq a(\zeta_y)\cdot y <0 = b\cdot y$. Thus $b$ cannot be separated from $\mathrm{Conv}(\AI)$ by \pg{any} hyperplane and so $b\in \mathrm{Conv}(\AI)$. There are nevertheless technicalities to ensure that the linear combination given by (\ref{rreqdd1})  is not trivial. These are overcome by restricting to a basis of $\FI$ on $V$. 

\subsubsection{Sketch of proof of Theorem \ref{C1}}
Let $\sigma:\R^M\rightarrow K$ be a linear isomorphism \pb{where $M=\mathrm{dim}(K)$}. We define $f_k(z):=M_k(\sigma(z))$ for $k=1,2,\dots, q_1$, where $M_k$ are minors in $M^{m\times n}$ and thus $f_k$ are homogeneous polynomials. Further define $f_{q_1+j}(z):=z_j$ for $j=1,\dots,M$. By properties of 
determinants (see Lemmas \ref{LA2} and \ref{LA3}) it is not hard to see this set of functions satisfy property $R$. It is also straightforward 
to show that measures \pg{in $\mathbb{M}^{pc}_{\FI}$} can be pushed forward via $\sigma$ to form measures \pg{in} $\MI^{pc}(K)$. As such Theorem 
\ref{C1} is essentially a corollary to Theorem \ref{T1}.

\subsection{Sketch of proof of Theorem \ref{CC1}} From Theorem \ref{C1} we have learned that given $\mu\in \MI^{pc}(K)$, to show 
that the support of $\mu$ can be \pg{reduced} to a lower dimensional cone we need only to find a linear combination of minors 
that is non-negative and non-trivial on $K$. Further if the minors we use are $2\times 2$ minors then our cone is actually a subspace. 

Let $\sigma:z\in\R^{d}\mapsto \lt(\begin{array}{ccc} a_{11}\cdot z &   \dots & a_{1n}\cdot z \\
\dots & \dots \\
a_{m1}\cdot z & \dots & a_{mn}\cdot z \end{array}\rt)\in K$ be a linear isomorphism. A major simplification comes from the following observation: by performing row and column operations on the matrix $\sigma(z)$ we arrive at a matrix $\ti{\sigma}(z):=\lt(\begin{array}{ccc} \ti{a}_{11}\cdot z &   \dots & \ti{a}_{1n}\cdot z \\
\dots & \dots \\
\ti{a}_{m1}\cdot z & \dots & \ti{a}_{mn}\cdot z \end{array}\rt)$, and defining $\wt{K}:=\lt\{ \ti{\sigma}(z):z\in \R^{d}\rt\}$, we arrive at a 
\em different subspace\rm. \pg{If $K$ has no Rank-$1$ connections, then $\wt K$ also has no Rank-$1$ connections.} Further 
\begin{equation}
\label{rreqqbb20}
\mathrm{Span}\lt\{M_1(\sigma(z)), \dots, M_{q_0}(\sigma(z)) \rt\}=\mathrm{Span}\lt\{M_1(\ti{\sigma}(z)), \dots, M_{q_0}(\ti{\sigma}(z)) \rt\},
\end{equation}
\pb{where $M_1,\dots,M_{q_0}$ are all $2\times 2$ minors in $M^{m\times n}$.} This is the content of Lemma \ref{r1lemma2}. Thus if we can find a sequence of row and column operations to reduce $\sigma(z)$ to a matrix $\ti{\sigma}(z)$ which has a simpler 
structure \pg{that allows to} find a linear combination of minors that is non-trivial and non-negative, then by (\ref{rreqqbb20}) \orgg{there must exist a} \pg{linear combination \pgr{which} also works for the subspace $K$}. 

 It turns out that the restrictions on $\lt\{a_{ij}\rt\}$ \pg{imposed} by $K$ having no Rank-$1$ connections are such that \pg{one can transform $\sigma(z)$ to some $\ti\sigma(z)$ such that (\ref{sseqc2}) can be checked relatively easily. The most delicate step in the proof is when $K$ is isomorphic to a three-dimensional subspace in $M^{3\times 3}$, in which case we need to invoke an argument of \v{S}ver\'{a}k to show that all \pb{three-dimensional} subspaces in $M^{3\times 3}_{sym}$ must contain Rank-$1$ connections. If $K\not\subset M^{3\times 3}_{sym}$, then carefully checking all $2\times 2$ minors in $K$ gives a linear combination satisfying (\ref{sseqc2}).}

\subsection{Sketch of proof of Theorem \ref{C17}}
\label{sktt1}
 The space of $k$-dimensional subspaces in $M^{m\times n}$ is trivially 
isomorphic to $G(k,p)$ for $p=mn$.  It is well known that $G(k,p)$ is a real analytic compact connected manifold of 
dimension $k(p-k)$. The charts for $G(k,p)$ can be found by fixing a \pg{pair of} transversal subspaces $W_0, W_1$ of $\R^{\pg{p}}$ where 
$\mathrm{dim}(W_0)=k$ and $\mathrm{dim}(W_1)=p-k$, then viewing the \pg{elements} of $G(k,p)$ as graphs of linear maps from 
$W_0$ to $W_1$. So fixing $W_0$, $W_1$ and choosing a basis for each space, each $A\in \R^{(p-k)k}\simeq M^{(p-k)\times k}$ \pg{defines} a linear mapping $T_A:W_0\rightarrow W_1$. \pg{Further define} $\phi_{W_0,W_1}(A):=\lt\{v+T_A(v):v\in W_0\rt\}\in G(k,p)$ \pg{and $\phi_{W_0,W_1}$ forms a chart for $G(k,p)$}. As we vary $W_0, W_1$ (smoothly varying our choice of basis) we obtain a complete set of charts. 

Now letting $M_1, M_2, \dots, M_{q_0}$ denote the set of all $2\times 2$ minors of $M^{m\times n}$, we define quadratics on 
$\R^k$ by $Q_j^A(y):=M_j\lt(\sum_{l=1}^k y_l (a_l+T_A(a_l)) \rt)$ for $j=1,2,\dots, q_0$ where $\lt\{a_1,a_2,\dots, a_k\rt\}$ is a 
basis of $W_0$. Each $Q_j^A$ can be represented by some $X^A_j\in M^{k\times k}_{sym}$.
The key point of the proof is the following: we are able to define a non-trivial real analytic function $\Lambda:\R^{k(p-k)}\rightarrow \R$ such that 
\begin{equation}
\label{gpseq2}
\mathrm{Span}\lt\{X^A_1, X^A_2, \dots, X^A_{q_0}\rt\}=M^{k\times k}_{sym}\text{ for all }\pg{A}\in \R^{\pg{(p-k)k}}\backslash \lt\{\pg{A}:\Lambda(A)=0\rt\}.
\end{equation}
\pg{Thus for all $A$ but the zero set of the analytic function $\Lambda$ \orgg{(which is ``small'')}, one can find a linear combination of $\lt\{X^A_1, X^A_2, \dots, X^A_{q_0}\rt\}$ that is positive definite, and this gives a linear combination of the $2\times 2$ minors that satisfies (\ref{sseqc2}).} Then the conclusions follow by very similar 
arguments to the proof of the sufficiency part of Theorem \ref{C1}. The existence of $\Lambda$ follows by identifying each \pg{symmetric} $X^A_j$ as a vector in $\R^{\frac{k(k+1)}{2}}$ and forming a matrix \pg{$\Pi(A)$} in $M^{\frac{(k+1)k}{2}\times q_0}$ with these vectors as columns. Then $\Lambda(A):=\det\lt(\Pi(A)\Pi(A)^T \rt)$ satisfies (\ref{gpseq2}). To show that $\Lambda$ is non-trivial, we notice that $\Lambda(A_0)\ne 0$ where $A_0$ defines the subspace \fred{\pb{$V_0$ given by} (\ref{graseq5}) of Lemma \ref{graslem2}}.

\subsection{Sketch of proof of Theorem \ref{t102}} \ared{As sketched briefly in the introduction,} the case where $\alpha'(\ti{\alpha}_2)<0$ follows easily from a well known result of \ared{\v{S}ver\'{a}k \cite{sv2}}. The case where 
$\alpha'(\ti{\alpha}_2)>0$ is the one that requires real work. \bblue{As in Subsection \ref{s2.1}, to streamline the sketch, we} consider the special case where $\ti{\alpha}=0$ \bblue{and $a(\ti{\alpha}_2)=0$}. \bblue{Given $s_0,t_0$ sufficiently small}, let 
\begin{equation*}
\zeta_0:= \lt(\begin{matrix} 0 & 0 \\ 0 & 0 \\ 0 & 0 \end{matrix}\rt), \; \zeta_1:= \lt(\begin{matrix} s_0 & 0 \\ 0 & s_0 \\ 0 & \frac{1}{2}s_0^2 \end{matrix}\rt), \;\; \zeta_2:= \lt(\begin{matrix} -s_0 & 0 \\ 0 & -s_0 \\ 0 & \frac{1}{2}s_0^2 \end{matrix}\rt)
\end{equation*}
and 
\begin{equation*}
 \zeta_3:= \lt(\begin{matrix} 0 & t_0 \\ a(t_0) & 0 \\ 0 & F(t_0) \end{matrix}\rt),\;\; \zeta_4:= \lt(\begin{matrix} 0 & -t_0 \\ a(-t_0) & 0 \\ 0 & F(-t_0) \end{matrix}\rt).
\end{equation*}
So $\zeta_0,\zeta_1,\dots, \zeta_4\in \mathcal{K}_1$.  \bblue{For $0<\ep<1$ sufficiently small, we construct non-trivial measures supported at the above five points, with weight $1-\ep$ at $\zeta_0$, and total weight $\ep$ at the other four points.} Let $D_1, D_2, D_3$ denote the $(1,2), (2,3), (1,3)$ minors of a $3\times 2$ matrix, respectively. We set the matrix
\begin{equation*} 
A:= \lt(\begin{matrix} D_1(\zeta_1) & D_1(\zeta_2) & D_1(\zeta_3) & D_1(\zeta_4) \\ D_2(\zeta_1) & D_2(\zeta_2) & D_2(\zeta_3) & D_2(\zeta_4) \\ D_3(\zeta_1) & D_3(\zeta_2) & D_3(\zeta_3) & D_3(\zeta_4) \\ 1 & 1 & 1 & 1 \end{matrix}\rt)
=\lt(\begin{matrix} s_0^2 & s_0^2 & -t_0a(t_0) & t_0a(-t_0) \\ 0 & 0 & a(t_0)F(t_0) & a(-t_0)F(-t_0) \\ \frac{1}{2}s_0^3 & -\frac{1}{2}s_0^3 & 0 & 0 \\ 1 & 1 & 1 & 1 \end{matrix}\rt).
\end{equation*}
\orgg{As a first step to obtain a non-trivial measure in $\MI^{pc}(\pgr{\K}_1)$, we construct a measure $\mu$ with \pgr{$\spt\mu\subset \lt\{\zeta_0,\zeta_1,\dots, \zeta_4\rt\}$}, $\int D_k(\pgr{\zeta}) d\mu=0$ for $k=1,2,3$ and $\mu(\pgr{\K_1}\backslash \lt\{\zeta_0\rt\})\pgr{=}\ep$. This is equivalent to finding some 
$\gamma_0\in \R^4_{+}$ such that $A\gamma_0=\pgr{(0,0,0,\ep)^T}$. To do this we use \pgr{the} Farkas-Minkowski Lemma (\org{see Corollary 7.1d, \cite{schr}}):
\begin{a1}[Farkas-Minkowski]
\label{fm}
Let $A\in M^{m\times n}$ be a matrix with columns $\lt\{a_1,a_2,\dots, a_n\rt\}$ and $b\in \R^{\pgr{m}}$. There exists $x\in \R^n_{+}$ such that $Ax=b$ if and only if 
$y\cdot b\geq 0$ for every vector $y\in \R^{\pgr{m}}$ with $y\cdot a_i\geq 0$ for $i=1,2,\dots, \pgr{n}$.
\end{a1}
By a careful analysis using the special structure of the points $\zeta_j$, we have that  $\sum_{i=1}^3 y_i \pgr{D_i(\zeta)}$ changes sign on $\lt\{\zeta_1,\zeta_2,\zeta_3,\zeta_4\rt\}$ for all non-trivial $y\in \R^3$. By arguments analogous to the last paragraph of Section \pgr{\ref{s2.1}} this allows us to apply \pgr{the} Farkas-Minkowski Lemma (indeed Farkas-Minkow\pgr{s}ki and the Hyperplane Separation Theorem are closely related results). So if we define $L^{\ep}(\gamma) := A\gamma - \lt(0,0,0,\ep\rt)^{T}$, then we have the existence of $\gamma_0\in\R^4_{+}$ such that $L^{\ep}(\gamma_0)=0$}. However what we need to solve for a measure in $\mathcal{M}^{pc}(\K_1)$ is $G^{\ep}(\gamma):=L^{\ep}(\gamma)-Q(\gamma)=0$, where
\begin{equation*}
Q(\gamma) := \lt(\begin{matrix} D_1\lt(\sum_{j=1}^{4}\gamma_j\zeta_j\rt) \\ D_2\lt(\sum_{j=1}^{4}\gamma_j\zeta_j\rt) \\ D_3\lt(\sum_{j=1}^{4}\gamma_j\zeta_j\rt) \\ 0 \end{matrix}\rt).
\end{equation*}
Since $G^{\ep}$ is a quadratic perturbation of an invertible function, it should seem reasonable that for small enough $\ep$, $G^{\ep}(\gamma)=0$ will have a solution. But to actually establish that the solution is non-negative we carry out an iterative argument inspired by the proof of the inverse function theorem. To this end, we start from the non-negative solution $\gamma_0$ of the linear part $L^{\ep}(\gamma)=0$, and use an iterative argument to solve for $\gamma_k$ in each step $k>0$ such that $\gamma_k$ converges to the actual solution of $G^{\ep}(\gamma)=0$. The convergence of this scheme is guaranteed by choosing $\ep$ sufficiently small. These are the contents of Lemmas \ref{l105} and \ref{l106}. 

What is slightly surprising is that to prove the general case we need to work instead with the set $\mathcal{K}^{\alpha}_1$ defined by (\ref{bvv22}) of Section \ref{prelims}. This set is essentially a stripping away of the quadratic part of $\mathcal{K}_1$ around a point $\alpha$ \pg{and similar ideas have been used by} DiPerna \cite{dp}. In some sense, the set $\K_1^{\alpha}$ plays the role of simplifying the problem by allowing the assumptions $\ti{\alpha}=0$ and $a(\ti{\alpha}_2)=0$.

\section{Proof of Theorem \ref{T1}}
\label{sec4}

\pg{The structure of real algebraic sets in $\R^n$ has been well studied.} In this section, we will make use of the following \emph{descending chain condition} for real algebraic sets, \pg{whose proof is a simple application of the classical Hilbert's Basis Theorem (see \cite{milnor}, page 9)}.
	\begin{a5}
		\label{p10}
		Any sequence $V_1\supsetneqq V_2\supsetneqq V_3\supsetneqq\dots$ of real algebraic sets must terminate after a finite number of steps.
	\end{a5}

Given a set of points $S\subset \R^n$, we \pblue{denote by $\mathrm{Conv}(S)$}  the \em convex hull \rm of $S$. \pblue{It is well known that, since $S$ is a subset of a finite dimensional space, its convex hull can be represented as} 
\begin{equation*}
\mathrm{Conv}(S)=\lt\{\blue{\sum_{i=1}^m \lm_i a_i: m\in \mathbb{N}, a_i\in S, \lm_i \in\R_{+}, \sum_{i=1}^m \lm_i=1}\rt\}.
\end{equation*}
%
\red{Given a vector $v\in \R^n$ we let $\lt[v\rt]_i$ be the $i$-th \pblue{component} of $v$.} Let $\FI$ be a finite collection of homogeneous polynomials on $\R^n$ and $\pg{V}$ be a \pblue{non-empty subset of} $\R^n$. We \pblue{denote by}
$\FI_{\pg{V}}$ a subset of $\FI$ such that $\lt\{f_{\lfloor \pg{V}}:f\in \FI_{\pg{V}}\rt\}$ forms a basis of the \fred{space} 
$\org{\mathrm{Span}\lt\{f_{\lfloor V}:f\in \FI\rt\}}$.

\begin{a1}
	\label{L1}
	Let $\FI=\lt\{f_1,f_2,\dots, f_{M_0}, \blue{f_{M_0+1}, \dots, f_{M_0+n}}\rt\}$ be a set of homogeneous polynomials \pblue{on $\R^n$} \blue{such that \pblue{$f_{M_0+j}(z)=z_j$ for $j=1,\dots,n$} (not necessarily satisfying property $R$)}.
	Suppose there exists \org{a set $V\pgr{\subset \R^n}$} \pgr{such that} \pblue{$\{0\}\subsetneqq \pg{V}$} \pgr{and}
	\begin{equation*}
	\text{for all }y\in \R^{M_0\blue{+n}}\backslash \lt\{0\rt\},\text{ the linear combination }\sum_{k=1}^{M_0\blue{+n}} y_k f_{k \lfloor \pg{V}}\text{ is either trivial or changes sign,}
	\end{equation*}
	then there exists non-trivial $\mu\in \Mb^{pc}_{\FI}(0)$.
\end{a1}

\begin{proof}
Let $\FI_\pg{V}=\lt\{f_{k_1}, f_{k_2}, \dots, f_{k_\pblue{N_1}}\rt\}$. \blue{Note that $\FI_\pg{V}$ must be non-empty, as otherwise $f_{M_0+j}(z)=z_j=0$ on $\pg{V}$ for all $j$ and hence \pblue{$\pg{V}=\{0\}$ which is a contradiction}.} Define
\begin{equation}
\label{eqa12}
a(\zeta):=\pblue{\lt(\begin{array}{c} f_{k_1}(\zeta)\\ f_{k_2}(\zeta)\\ \dots \\  f_{k_\pblue{N_1}}(\zeta)\end{array}\rt)} \pblue{\text{ for }\zeta\in\R^n}\text{ and }b:=\pblue{\lt(\begin{array}{c} 0  \\ 0 \\ \dots \\ 0 \end{array}\rt)}
\end{equation}
to be vectors in \pblue{$\R^{N_1}$} and let 
\begin{equation*}
\AI:=\lt\{a(\zeta):\zeta\in \blue{\pg{V}\setminus\{0\}} \rt\}.
\end{equation*}
Note that $b=a(0)$. We claim that $b\notin\AI$. Suppose not,  then there exists some $\zeta\in \pg{V}\setminus\{0\}$ such that $f_{k_j}(\zeta)=0$ for all $j=1,\dots,\pblue{N_1}$. As $\FI_\pg{V}$ forms a basis of $\org{\mathrm{Span}\lt\{f_{\lfloor V}:f\in \FI\rt\}}$, it follows that $f_{k}(\zeta)=0$ for all $k=1,\dots, M_0+n$. However, this implies that $\zeta_j=f_{M_0+j}(\zeta)=0$ for all $j\pblue{=1,\dots,n}$, and hence $\zeta=0$, which is a contradiction. 

We will show that $b\in \blue{\mathrm{Conv}(\AI)}$ \pblue{by using the Hyperplane Separation Theorem}. First note that
\begin{equation}
\label{eqa11.5}
\text{for all }y\in \pblue{\R^{N_1}\setminus\{0\}},\text{ the linear combination }\sum_{i=1}^{\pblue{N_1}} y_i f_{k_i \lfloor \pg{V}}\text{ changes sign }
\end{equation}
since $\lt\{f_{k_1}, f_{k_2}, \dots, f_{k_\pblue{N_1}}\rt\}$ is linearly independent \pblue{on $\pg{V}$}. Let $\pblue{\BI}=\lt\{b\rt\}$. Note that $\blue{\mathrm{Conv}(\AI)}$ and $\BI$ are both convex sets. \org{Let} \blue{$y\in \R^{\pblue{N_1}}\setminus\{0\}$ and $c\in\R$} \org{be} such that 
\begin{equation}
\label{eqa17}
w\cdot y\geq \blue{c}\text{ for all }w\in \blue{\mathrm{Conv}(\AI)}.
\end{equation}
\org{By} (\ref{eqa11.5}) there exists $\zeta_y\in \blue{\pg{V}\setminus\{0\}}$ such that 
\begin{equation}
\label{eqa19.4}
\sum_{i=1}^{\pblue{N_1}} y_i f_{k_i}(\zeta_y)<0. 
\end{equation}
Now as $a(\zeta_y)\in \AI$, by (\ref{eqa17}) and (\ref{eqa19.4}) we have that 
\begin{equation}
\label{eqa18}
0>\sum_{i=1}^{\pblue{N_1}}  y_i f_{k_i}(\zeta_y)=a(\zeta_y)\cdot y\geq \blue{c}.
\end{equation}
Thus $y\cdot b=0\overset{\org{(\ref{eqa18})}}{>}c$. By the Hyperplane Separation Theorem (\pblue{see, e.g.,} \cite{boyd} Exercise 2.22) this \blue{implies that $b\in\mathrm{Conv}(\AI)$}. 

As $b\in \blue{\mathrm{Conv}(\AI)}$, there exists $\lm_1,\lm_2, \dots, \lm_{p_0}\in \R_{+}$ and 
$\zeta_1,\zeta_2,\dots,\zeta_{p_0}\blue{\in \pg{V}\setminus\{0\}}$ such that 
$\sum_{i=1}^{p_0} \lm_i=1$ and 
\begin{equation}
\label{eqa21}
b=\sum_{i=1}^{p_0} \lm_i a(\zeta_i).
\end{equation}
Let $\mu:=\sum_{i=1}^{p_0} \lm_i \delta_{\zeta_i}$. \blue{Note that $\mu$ is non-trivial since $b\notin\AI$.} \blue{We claim that $\overline{\mu}=0$. To see this, it suffices to show that
	\begin{equation}\label{eqb1}
	\sum_{i=1}^{p_0} \lm_i \red{\lt[\zeta_{i}\rt]_j} = 0 \text{ for all } j=1,\dots,n.
	\end{equation} 
	As $\FI_\pg{V}$ is a basis of $\org{\mathrm{Span}\lt\{f_{\lfloor V}:f\in \FI\rt\}}$, we have
	\begin{equation*}
	f_{M_0+j} = \sum_{r=1}^{\pblue{N_1}}\alpha^j_r f_{k_r}\pblue{\text{ on }\pg{V}}\red{\text{ for some }\alpha^j\in \R^{\pblue{N_1}}}
	\end{equation*}
	and hence
	\begin{equation*}
	\sum_{i=1}^{p_0} \lm_i \red{\lt[\zeta_{i}\rt]_j} = \sum_{i=1}^{p_0} \lm_i f_{M_0+j}(\zeta_{i}) = \sum_{i=1}^{p_0} \sum_{r=1}^{\pblue{N_1}}\lm_i\alpha_r^{\pblue{j}} f_{k_r}(\zeta_i) =  \sum_{r=1}^{\pblue{N_1}}\alpha_r^{\pblue{j}} \lt(\sum_{i=1}^{p_0} \lm_i f_{k_r}(\zeta_i)\rt)\overset{(\ref{eqa21}), \org{(\ref{eqa12})}}{=} 0.
	\end{equation*}
	This shows (\ref{eqb1}) for all $j\in \lt\{1,2,\dots, n\rt\}$ and therefore $\overline{\mu}=0$.} Now 
\begin{equation*}
\red{\int f_{k_\pblue{r}}(\pblue{z}) d\mu(\pblue{z})=\sum_{i=1}^{p_0} \lm_i f_{k_r}(\zeta_i)\overset{(\ref{eqa12})}{=}\sum_{i=1}^{p_0} \pblue{\lm_i}\lt[a(\zeta_i)\rt]_r\overset{(\ref{eqa21})}{=} \lt[b\rt]_r\overset{ (\ref{eqa12})}{=}0\text{ for }r=1,2,\dots,N_1}
\end{equation*}
and thus $\mu\in \pblue{\Mb}_{\FI_\pg{V}}^{pc}(0)$. 

\pblue{Finally we show that $\mu \in \Mb^{pc}_{\FI}(0)$. For any $f\in \FI$ there exists $\beta\in\R^{N_1}$ such that $f_{\lfloor \pg{V}}=\sum_{i=1}^{N_1} \beta_i f_{k_i \lfloor \pg{V}}$. It follows that
	\begin{equation*}
		\begin{split}
			&\int_{\R^n} f(z) d\mu(z)=\int_{\pg{V}} f(z) d\mu(z)\\
			&\qd\qd=\sum_{i=1}^{N_1} \beta_i\int_{\pg{V}} f_{k_i}(z) d\mu(z)=\sum_{i=1}^{N_1} \beta_i f_{k_i}(0)=f(0),
		\end{split}
	\end{equation*}
	and hence $\mu\in \mathbb{M}^{pc}_{\FI}(0)$. }
\end{proof}

%
%
%

\begin{a2}
	\label{P1}
	Let $\FI=\lt\{f_1,f_2,\dots, f_{M_0}, \blue{f_{M_0+1}, \dots, f_{M_0+n}}\rt\}$ be a set of homogeneous polynomials \pblue{on $\R^n$} \blue{such that \pblue{$f_{M_0+j}(z)=z_j$ for $j=1,\dots,n$} (not necessarily satisfying property $R$)}. Then $\Mb^{pc}_{\FI}(0)=\lt\{\delta_{0}\rt\}$ if and only if (\ref{eqcc1})  of Theorem \ref{T1} holds true. 
\end{a2}

\begin{proof}
\pblue{Suppose (\ref{eqcc1}) of Theorem \ref{T1} is false, then clearly Lemma \ref{L1} gives a non-trivial $\mu\in \Mb_{\FI}^{pc}(0)$. } So in the following we assume (\ref{eqcc1}), and thus for 
every \blue{non-trivial algebraic} cone $\pg{V}$ there exists $y\in \R^{M_0\pblue{+n}}\backslash \lt\{0\rt\}$ such that 
\begin{equation*}
\sum_{k=1}^{M_0\pblue{+n}} y_k f_k\geq 0\text{ and }\sum_{k=1}^{M_0\pblue{+n}}  y_k f_k\not\equiv 0\text{ on }\blue{\pg{V}}. 
\end{equation*}
Given any $\mu\in \Mb^{pc}_{\FI}(0)$, we will show that $\mu=\delta_{0}$. 

Now take $\pg{V}_1=\R^{n}$, so there exists $y^1\in \R^{M_0\pblue{+n}}\backslash \lt\{0\rt\}$ such that 
$\sum_{k=1}^{M_0\pblue{+n}} y_k^1 f_k$ is non-negative and non-trivial on $\pg{V}_1$. Let 
\begin{equation}
\label{eqa7}
m_1:=\pblue{\mathrm{deg}\lt(\sum_{k=1}^{M_0\pblue{+n}} y_k^1 f_k\rt)}
\end{equation}
and
\begin{equation*}
\MI_1:=\lt\{\pblue{k\in \lt\{1,2,\dots, M_0\pblue{+n}\rt\}: \mathrm{deg}(f_k)=m_1}\rt\}.
\end{equation*}
We claim that
\begin{equation}
\label{eqa10.5}
\sum_{k\in \MI_1} y^1_k f_k\text{ is non-trivial and non-negative on }\pg{V}_1.  
\end{equation}
\pblue{First because of the definition of $m_1$, it is clear that $\sum_{k\in \MI_1} y^1_k f_k$ is non-trivial.} Now suppose it changes sign, then 
there exists $\zeta_1\in \pg{V}_1$ such that 
\begin{equation}
\label{eqa11}
\sum_{k\in \MI_1} y^1_k f_k(\zeta_1)<0. 
\end{equation}
\pblue{Note that all $f_k$'s with degree higher than $m_1$, if any, cancel out in $\sum_{k=1}^{M_0\pblue{+n}} y_k^1 f_k$.} Now, \pblue{letting $d_k:=\mathrm{deg}(f_k)$, we have}
\begin{equation}
\label{eq30}
\begin{split}
\sum_{k=1}^{M_0\pblue{+n}} y^1_k f_k(\lm \zeta_1)&=\sum_{k\in \lt\{1,2,\dots, M_0\pblue{+n}\rt\}\backslash \MI_1} y^1_k f_k(\lm \zeta_1)+\sum_{k\in \MI_1} y^1_k f_k(\lm \zeta_1)\\
&=\sum_{k\in \lt\{1,2,\dots, M_0\pblue{+n}\rt\}\backslash \MI_1} y^1_k \lm^{d_k} f_k(\zeta_1)+\sum_{k\in \MI_1} y^1_k \lm^{m_1} f_k(\zeta_1)\\
&=\lm^{m_1}\lt(\sum_{k\in \lt\{1,2,\dots, M_0\pblue{+n}\rt\}\backslash \MI_1} y^1_k \lm^{d_k-m_1} f_k(\zeta_1)+\sum_{k\in \MI_1} y^1_k  f_k(\zeta_1)  \rt)\\
&\overset{(\ref{eqa7})}{\leq} \frac{\lm^{m_1}}{2}\lt(\sum_{k\in\MI_1} y^1_k  f_k(\zeta_1)  \rt)\text{ for all large enough }\lm\pblue{>0}. 
\end{split}
\end{equation}
Together with (\ref{eqa11}) this contradicts the fact that $\sum_{k=1}^{M_0\pblue{+n}} y^1_k f_k$ is \pg{non-negative} on $\pg{V}_1$, and thus (\ref{eqa10.5}) is established. 

Since $\sum_{k\in \MI_1} y^1_k f_k$ is a homogeneous polynomial of degree $m_1$ the set 
\begin{equation*}
\pg{V}_2:=\lt\{z\in \R^n : \sum_{k\in \MI_1} y^1_k f_k(z)=0 \rt\} 
\end{equation*}
forms an \blue{algebraic} cone. \blue{Further, since $\sum_{k\in \MI_1} y^1_k f_k(z)$ is non-trivial on $\pg{V}_1=\R^n$, we have $\pg{V}_1\supsetneqq \pg{V}_2$.} Note that since $\mu\in \Mb_{\FI}^{pc}(0)$ we have that 
\begin{equation}
\label{eqa41}
\int \sum_{k\in \MI_1} y^1_k f_k(z) d\mu(z)=\sum_{k\in \MI_1} y^1_k f_k(\ol\mu)=\sum_{k\in \MI_1} y^1_k f_k(0)=0.
\end{equation}
So we must have that 
\begin{equation}
\label{eqa44}
\mu\lt(\R^n\backslash  \pg{V}_2\rt)=0. 
\end{equation}

If $\pg{V}_2=\{0\}$, then we are done because of (\ref{eqa44}). So suppose $\pg{V}_2$ is non-trivial. By hypothesis there exists $y^2\in \R^{M_0\pblue{+n}}\backslash \lt\{0\rt\}$ such that
\begin{equation}
\label{eqaa40}
\sum_{k=1}^{M_0\pblue{+n}} y^2_k f_k \text{ is non-trivial and non-negative on } \pg{V}_2.
\end{equation}
Now we repeat the arguments as above. Let 
\begin{equation*}
m_2:=\pblue{\mathrm{deg}\lt(\sum_{k=1}^{M_0\pblue{+n}} y_k^2 f_{k \lfloor \pg{V}_2}\rt)} \text{ and }\MI_2:=\lt\{\pblue{k\in \lt\{1,2,\dots, M_0\pblue{+n}\rt\}: \mathrm{deg}(f_k)=m_2}\rt\}.
\end{equation*}
Now we claim 
\begin{equation}
\label{eqa10}
\sum_{k\in \MI_2} y^2_k f_k\text{ is non-trivial and non-negative on }\pg{V}_2.  
\end{equation}
\pblue{Again by definition of $m_2$ we know that $\sum_{k\in \MI_2} y^2_k f_k$ is non-trivial on $\pg{V}_2$.} If it changes sign on $\pg{V}_2$, then there exists $\zeta_2\in \pg{V}_2$ such that $\sum_{k\in \MI_2} y^2_k f_k(\zeta_2)<0$. Since $\lm\zeta_2\in \pg{V}_2$ for any $\lm>0$,  we can argue in an identical manner to 
(\ref{eq30}) and conclude that for large enough $\lm$, 
$$
\sum_{k=1}^{M_0\pblue{+n}} y^2_k f_k(\lm \zeta_2)\leq \frac{\lm^{m_2}}{2}\lt( \sum_{k\in \MI_2} y^2_k f_k(\zeta_2)  \rt)<0.
$$
This contradicts (\ref{eqaa40}) and thus (\ref{eqa10}) is established. 

\blue{Now in exactly the same way as 
(\ref{eqa41}) we have that $\int \sum_{k\in \MI_2} y^2_k f_k(z)  d\mu (z)=0$. Hence letting 
\begin{equation*}
\pg{V}_3:=\lt\{z\in \pg{V}_2: \sum_{k\in \MI_2} y^2_k f_k(z)=0\rt\}
\end{equation*}
we have that $\mu(\pg{V}_2 \backslash \pg{V}_3)=0$. Further $\pg{V}_3$ is an algebraic cone satisfying $\pg{V}_2\supsetneqq \pg{V}_3$. If $\pg{V}_3=\{0\}$, then we are done. Otherwise, we can repeat the above process to obtain a descending chain of algebraic cones $\pg{V}_1\supsetneqq \pg{V}_2\supsetneqq \pg{V}_3\supsetneqq \dots$. By Proposition \ref{p10}, after a finitely many steps, the chain must stop. Let $\pg{V}_p$ be the last algebraic cone in the chain. We claim that $\pg{V}_p=\{0\}$. Assume not, then $\pg{V}_p$ is a non-trivial algebraic cone. By hypothesis and the arguments as above, there exists an algebraic cone $\pg{V}_{p+1}\subsetneqq \pg{V}_{p}$, which is a contradiction. As $\mathrm{Spt} \mu \subset \pg{V}_p$, we conclude that $\mu=\delta_{0}$.} This completes the proof of the theorem. 
\end{proof}

\begin{proof}[Proof of Theorem \ref{T1}]
	\pblue{For any $\varpi\in \R^n$, we define the translation $P^{\varpi}:\R^n\rightarrow \R^n$ by}
	\begin{equation}
	\label{eqbb35}
	P^{\varpi}(z):=z-\varpi.
	\end{equation}
	By Lemma \ref{L6}, \pblue{for a collection of polynomials $\FI$ satisfying property $R$}, we have $\mu\in \Mb_{\FI}^{pc}(\varpi)$ if and only if $\lt(P^{\varpi}\rt)_{\sharp} \mu\in \Mb_{\FI}^{pc}(0)$, where \pblue{$\lt(P^{\varpi}\rt)_{\sharp} \mu$ is the push forward of $\mu$ under the mapping $P^{\varpi}$}. So it suffices to show 
	that $\Mb_{\FI}^{pc}(0)$ consist of Diracs if and only if (\ref{eqcc1}) of Theorem \ref{T1} holds true. This is exactly the 
	content of Theorem \ref{P1}. \pg{Thus condition (\ref{eqcc1}) holds true if and only if $\Mb_{\FI}^{pc}(\varpi)$ is trivial for all $\varpi\in \R^n$, and hence if and only if $\Mb_{\FI}^{pc}$ consists of Dirac measures.}
\end{proof}

%
%

\section{Proof of Theorem \ref{C1}}
\label{sec5}
In this section, we give the proof of Theorem \ref{C1}. \pblue{The following notation will be used at multiple places throughout this paper}. Given a matrix $A\in M^{m\times n}$, let 
\begin{equation}
\label{eqp2}
\pblue{R_i(A)}\in M^{1\times n}\text{ denote the }i\text{-th row of }A
\end{equation}
and
\begin{equation}
\label{eqp3}
\pblue{\lt[A\rt]_{ij} \text{ denote the }(i,j)\text{-th entry of }A.}
\end{equation}
\pblue{In the following we will deal with submatrices whose sizes vary. So we introduce the following notation. }For positive integers $m,n$ let 
\begin{equation}
\label{eqn1}
\pblue{M}^{m,n}_1(A), \pblue{M}^{m,n}_2(A), \dots, \pblue{M}^{m,n}_{q(m,n)}(A)\text{ denote all the minors of a matrix }A\in M^{m\times n},
\end{equation}
\pblue{where $q(m,n)$ denotes the number of minors in $M^{m\times n}$}.

\begin{proof}[Proof of Theorem \ref{C1}]
Let $M=\mathrm{dim}(K)$. There exists a linear isomorphism $\sigma:\R^M\rightarrow K$ such that
\begin{equation}
\label{zzeqb35}
\sigma(z)=\lt(\begin{array}{cccc} a_{11}\cdot z & a_{12}\cdot z & \dots & a_{1n}\cdot z \\
a_{21}\cdot z & a_{22}\cdot z & \dots & a_{2n}\cdot z \\
\dots & \dots \\
a_{m1}\cdot z & a_{m2}\cdot z & \dots & a_{mn}\cdot z \end{array}\rt)
\end{equation}
for $a_{ij}\in\R^M$. We claim that 
\begin{equation}
\label{rreqb31}
\dim\lt(\mathrm{Span}\lt\{a_{ij}\rt\}\rt)=M.
\end{equation}
\red{Suppose this is false, then} there exists $z_0\in \pblue{\R^M\setminus\{0\}}$ such that $a_{ij}\cdot z_0=0$ for all $i\in\lt\{1,2,\dots, m\rt\}$, 
$j\in\lt\{1,2,\dots, n\rt\}$. Thus $\sigma(\pblue{z_0})=0$ \pblue{which contradicts the fact that $\sigma$ is an isomorphism.}  Hence by 
(\ref{rreqb31}) there exist $\{\lm_{k}^{ij}\}$ such that
\begin{equation}
\label{eqb2}
\blue{z_k=\sum_{i,j}\lm_{k}^{ij}\lt(a_{ij}\cdot z\rt)\text{ for all } k=1,2,\dots,M.}
\end{equation}
Define 
\begin{equation}
\label{eqbba12}
f_k(z):=M_k(\sigma(z))\text{ for }k=1,2,\dots, \pblue{\pg{q_1}}\text{ and }\blue{f_{\pg{q_1}+j}(z):= z_j}\text{ for }j=1,2,\dots,M,
\end{equation}
and let $\mathcal{F}:=\lt\{f_1,f_2,\dots, f_{\pg{q_1}+M}\rt\}$. \nl

\em Step 1. \rm Let $\nu \in \mathcal{P}(\R^M)$ and $\mu:=\sigma_{\sharp} \nu$, i.e., $\mu$ is the push forward of $\nu$ under the mapping $\sigma$, then $\nu\in \mathbb{M}^{pc}_{\mathcal{F}}$ if and only if $\mu\in \mathcal{M}^{pc}(K)$.

\em Proof of Step 1. \rm By change of variable formula for push forward measures (see \cite{afp}, P. 32) for $k\in \lt\{1,2,\dots, \pg{q_1}\rt\}$ we have
\begin{equation*}
	\int f_k(z) d\nu(z)=\int M_k(\sigma(z)) d\nu (z)=\int M_k(X) d\mu(X)
\end{equation*}
and 
\begin{equation*}
	M_k\lt(\int X d\mu(X)\rt)=M_k\lt(\int \sigma(z) d\nu (z)\rt) =M_k\lt(\sigma\lt(\int z d\nu (z) \rt)\rt)=f_k\lt(\ol{\nu}\rt).
\end{equation*}
\pblue{This establishes Step 1}. \nl

\em Step 2. \rm We will show that the set of functions $\mathcal{F}$ has property $R$. 

\em Proof of Step 2. \rm For each $k\in \lt\{1,2,\dots,  \pblue{\pg{q_1}}\rt\}$, $M_k$ is a minor and as such is the 
determinant of a $p_k\times p_k$ submatrix for some $p_k\in \lt\{2,\dots, \min\lt\{m,n\rt\}\rt\}$.  So for each $k$ there exists a linear mapping $\red{P_k}:M^{m\times n}\rightarrow M^{p_k\times p_k}$ defined by pairwise \pblue{distinct} sets 
$I:=\lt\{i_1,i_2,\dots, i_{p_k}\rt\}$ and $J:=\lt\{j_1,j_2,\dots, j_{p_k}\rt\}$ such that \pblue{$P_k(A)=\lt\{\lt([A]_{ij}\rt): i\in I, j\in J\rt\}$ for all $A\in M^{m\times n}$}.
Now using  Lemma \ref{LA3} for the third equality \red{(recalling definition (\ref{eqn1}))} we have \pblue{for any $z_0\in\R^M$}
\begin{equation}
\label{eqbba3.2}
\begin{split}
&f_k(z+z_0)=M_k\lt(\sigma(z)+\sigma(z_0)\rt)\\
&\qd=\det\lt(P_k(\sigma(z))+ P_k(\sigma(z_0))\rt)\\
&\qd\overset{(\ref{eqap12})}{=}\det\lt(P_k(\sigma(z))\rt)+ \det\lt(P_k(\sigma(z_0))\rt)+\sum_{\pblue{l}=1}^{q(p_k,p_k)} \PI_l\lt(P_k(\sigma(z_0))\rt)M_l^{p_k,p_k}\lt(P_k(\sigma(z))\rt)\\
&\qd= \pg{f_k(z) + f_k(z_0) + \sum_{\pblue{l}=1}^{q(p_k,p_k)} \PI_l\lt(P_k(\sigma(z_0))\rt)M_l^{p_k,p_k}\lt(P_k(\sigma(z))\rt)},
\end{split}
\end{equation}
\pblue{where $\PI_l\lt(P_k(\sigma(z_0))\rt)$ are polynomial functions of the entries of $P_k(\sigma(z_0))$}. For each \red{$\pblue{l}\in \lt\{\pblue{1},\dots, q(p_k,p_k)\rt\}$} we have
\begin{equation*}
M_l^{p_k,p_k}\circ P_k=M_{k_l}\text{ for some }\red{k_l\in \lt\{1,2,\dots, \pg{q_1+mn}\rt\}}.
\end{equation*}
\pg{If $1\leq k_l\leq q_1$, then $M_{k_l}(\sigma(z))=f_{k_l}(z)$. If $q_1+1\leq k_l\leq q_1+mn$, then $M_{k_l}(\sigma(z))$ is a projection mapping of the form $a_{ij}\cdot z$ by (\ref{zzeqb35}), and thus by (\ref{eqbba12}) is a linear combination of $\{f_{q_1+j}\}$ for $j=1,\dots,M$. Hence, we see from (\ref{eqbba3.2}) that $f_k(z+z_0)\in \LI(\FI)$ (defined in (\ref{eqa1})).}
As this is true for each \red{$k\in \lt\{1,2,\dots, \pg{q_1}\rt\}$ and is trivially true for $f_{\pg{q_1}+j}$ \pblue{for $j=1,\dots,M$}}, we have shown that $\FI$ has property $R$. This completes the 
proof of Step 2. \nl

\em Step 3. \rm We will show that for $\FI$ \pblue{consisting of the polynomials} defined by (\ref{eqbba12}), condition (\ref{sseqc2}) of Theorem \ref{C1} is equivalent 
to condition (\ref{eqcc1}) of Theorem \ref{T1}. 

\em Proof of Step 3. \rm First note that $\pg{V}\subset \R^n$ is a non-trivial \blue{algebraic} cone \blue{if and only if} $\wt{\pg{V}}:=\sigma(\pg{V})$ is a 
non-trivial \blue{algebraic} cone in $K$. \blue{Indeed, first assume that $\pg{V}\subset \R^n$ is a non-trivial \blue{algebraic} cone.} To see that $\wt{\pg{V}}$ is a cone, take $v\in \wt{\pg{V}}$, $\lm \pblue{>}0$, then 
$\sigma^{-1}(\lm v)=\lm \sigma^{-1}(v)\in \pg{V}$, so $\lm v\in \wt{\pg{V}}$. \blue{To see that $\wt{\pg{V}}$ is an algebraic set, let $g$ be any polynomial function that vanishes on $\pg{V}$ and define $\wt{g}(\zeta):=g(\sigma^{-1}(\zeta))$. As $g$ is a polynomial function \pg{of $\sigma^{-1}(\zeta)$ and each coordinate of $\sigma^{-1}(\zeta)$ can be represented as a linear combination of the entries of $\zeta$ by (\ref{eqb2}), it follows that $\wt{g}$ is a polynomial function of the entries of $\zeta$}. It is also clear that $\wt{g}$ vanishes on $\wt{\pg{V}}$. This shows that $\wt{\pg{V}}$ is an algebraic set in $K$. Conversely, assume that $\wt{\pg{V}}$ is an algebraic cone in $K$. Almost identical arguments as above show that $\pg{V}$ is an algebraic cone in $\R^M$. }

\red{Now suppose we have condition (\ref{sseqc2}) of Theorem \ref{C1}. \pblue{Then for any non-trivial algebraic cone $\wt{\pg{V}}\subset K$}, there exists 
$\beta\in \R^{\pg{q_1}\pblue{+mn}}\backslash \lt\{0\rt\}$ such that
\begin{equation*}
\sum_{k=1}^{\pblue{\pg{q_1}}} \beta_k M_{k}\lt(\sigma(z)\rt)+\sum_{k=\pg{q_1}+1}^{\pg{q_1}\pblue{+mn}} \beta_k\lt(\sigma(z)\rt)_k\geq 0 \text{ and }\sum_{k=1}^{\pg{q_1}} \beta_k M_{k}\lt(\sigma(z)\rt)+\sum_{k=\pg{q_1}+1}^{\pg{q_1}+mn} \beta_k\lt(\sigma(z)\rt)_k\not\equiv 0\text{ on }\pg{V}
\end{equation*}
\pblue{for $\pg{V}=\sigma^{-1}(\wt{\pg{V}})$.} By (\ref{eqbba12}) we have that for $\FI$, condition (\ref{eqcc1})} of Theorem \ref{T1} holds true. \red{Next suppose condition (\ref{eqcc1}) of Theorem \ref{T1} holds true for $\FI$. Note 
that $\sum_{k=\pg{q_1}+1}^{\pg{q_1}+M}   y_k f_k(z)\overset{(\ref{eqbba12}),  (\ref{eqb2})}{=} \sum_{k=\pg{q_1}+1}^{\pg{q_1}+M}   y_k \sum_{i,j} \lm^{ij}_k (a_{ij}\cdot z)$, so this together with \pblue{(\ref{eqbba12})} gives that the non-trivial and non-negative linear combination we have in $\FI$ is actually one that we can express as \pg{a} linear combination \pblue{in $\{M_{k}\}$ \pg{for $k=1,\dots,q_1+mn$}.} Hence we have condition (\ref{sseqc2}) of Theorem \ref{C1}.} This completes the proof of Step 3. \nl 

\em Proof of Theorem \ref{C1} completed. \rm Let $\mu\in \MI^{pc}(K)$. By Step 1, we have $\nu:=\pblue{\lt(\sigma^{-1}\rt)_{\sharp}} \mu\in \Mb^{pc}_{\FI}$. \pblue{So $\sigma$ establishes a one-to-one correspondence between $\MI^{pc}(K)$ and $\Mb^{pc}_{\FI}$. Since $\sigma$ is an isomorphism, it is clear that $\MI^{pc}(K)$ is trivial if and only if $\Mb^{pc}_{\FI}$ is trivial. By Step 2, $\FI$ is a set of homogeneous polynomials with property $R$ and from (\ref{eqbba12}) it is clear that $\FI$ satisfies the assumptions of Theorem \ref{T1}. By Theorem \ref{T1}, $\Mb^{pc}_{\FI}$ is trivial if and only if condition (\ref{eqcc1}) is satisfied, which is equivalent to condition (\ref{sseqc2}) holding true by Step 3. The conclusion of Theorem \ref{C1} hence follows from the above equivalence relations.}
\end{proof}

%
%

\section{Proof of Theorem \ref{CC1}}

In this section we give the proof of Theorem \ref{CC1}. Our main tool is the following

\begin{a2}
	\label{l16}
	Let $d\in\{1,2,3\}$ and $K\subset M^{m\times n}$ be a $d$-dimensional subspace without Rank-$1$ connections. \pblue{Denote by $M_1,\dots, M_{q_0}:M^{m\times n}\rightarrow \R$ the set of all \pg{$2\times 2$} minors in $M^{m\times n}$}. Then there exists some $\beta\in \pg{\R^{q_0}\setminus\{0\}}$ such that 
	\begin{equation}
	\label{sseqc1}
	\sum_{k=1}^{q_0} \beta_k M_k(X)\geq 0\text{ for all }X\in K\text{ and }\sum_{k=1}^{q_0} \beta_k M_k\not\equiv 0\text{ on }K.
	\end{equation}
\end{a2}

The proof of the above theorem requires some preparation. We begin by introducing some notation. Given $A\in M^{m\times n}$, we denote 
\begin{equation}
\label{eqr1a70}
M^{n_1,n_2}_{m_1,m_2}(A):=\det\lt(\begin{matrix} \pg{[A]}_{m_1 n_1} & \pg{[A]}_{m_1 n_2}  \\
 [A]_{m_2 n_1} & [A]_{m_2 n_2}  \end{matrix}\rt) 
\end{equation}
for $m_1\ne m_2\in \lt\{1,2,\dots, m\rt\}$,  $n_1\ne n_2\in \lt\{1,2,\dots, n\rt\}$. Let $K\subset M^{m\times n}$ be a \zblue{$d$-dimensional} subspace, then there exist $a_{ij}\in \R^d$ for $i=1,\dots,m$ and $j=1,\dots,n$ such that 
\begin{equation}
\label{eqr1a1}
P_K(z):=\lt(\begin{array}{ccccc} a_{11}\cdot z & a_{12}\cdot z & \dots & a_{1n}\cdot z \\
a_{21}\cdot z & a_{22}\cdot z & \dots & a_{2n}\cdot z\\
\dots \\
\dots \\
a_{m1}\cdot z & a_{m2}\cdot z & \dots & a_{mn}\cdot z   \end{array}\rt)
\end{equation}
is a \ared{linear isomorphism of $\R^d$ \cblue{onto} $K$ and hence is a parametrization}. Thus we have
\begin{equation*}
	K = \lt\{P_K(z): z\in\R^d\rt\}.
\end{equation*}
Note that every linear isomorphism $P:\R^d\rightarrow M^{m\times n}$ corresponds uniquely to some $P(z)$ in the form (\ref{eqr1a1}), which can be identified as an $m\times n$ matrix with entries in the polynomial ring $\R[z_1,\dots,z_d]$. For the rest of this section, we do not distinguish between such linear isomorphisms and the associated matrices in the form (\ref{eqr1a1}) with entries in $\R[z_1,\dots,z_d]$. We define an equivalence relation between linear isomorphisms from $\R^d$ into $M^{m\times n}$ as follows.
\begin{deff}
	\label{def16}
	Let $P_1, P_2:\R^d\rightarrow M^{m\times n}$ be two linear isomorphisms. We say that $P_1$ is equivalent to $P_2$, written as
	\begin{equation} 
	\label{eqp103}
	P_1\sim P_2,
	\end{equation}
	if $P_2(z)$ can be obtained from $P_1(z)$, both viewed as $m\times n$ matrices with entries in the polynomial ring $\R[z_1,\dots,z_d]$, by finitely many elementary row and column operations.
\end{deff}
The following result will be used repeatedly in the proof of Theorem \ref{l16}.
\begin{a1}\label{r1lemma2}
	Let $d, m, n$ be positive integers such that $\min\{m,n\}\geq 2$ and $d\leq mn$. Denote by $M_1, M_2, \dots, M_{q_0}$ all $2\times 2$ minors of $M^{m\times n}$. Further let $P_1,P_2:\R^d\rightarrow M^{m\times n}$ be two linear isomorphisms such that $P_1\sim P_2$ in the sense of Definition \ref{def16}. Then, denoting $K_j:=P_j(\R^d)$ for $j=1,2$, 
	we have that $K_1$ has no Rank-$1$ connections if and only if $K_2$ has no Rank-$1$ connections. Further, we have
	\begin{equation}
	\label{eqr1a5}
	\mathrm{Span}\lt\{M_1(P_1(z)), \dots,  M_{q_0}(P_1(z))\rt\}=\mathrm{Span}\lt\{M_1(P_2(z)), \dots,  M_{q_0}(P_2(z))\rt\}
	\end{equation}
	as subsets of the polynomial ring $\R[z_1,\dots,z_d]$.
\end{a1}

\begin{proof} 
	By induction, it suffices to consider the case where $P_2(z)$ is obtained from $P_1(z)$, both viewed as $m\times n$ matrices with entries in  $\R[z_1,\dots,z_d]$, by an elementary row or column operation. We only show the case where $P_2(z)$ is obtained from $P_1(z)$ by an elementary row operation, as the proof for column operation is identical. 
	
	\xxred{Note that, for any fixed $z_0\in\R^d$, $P_1(z_0)$ and $P_2(z_0)$ are $m\times n$ matrices with entries in $\R$. As $P_2(z_0)$ is obtained from $P_1(z_0)$ by an elementary row operation, it is clear that
		\begin{equation}
		\label{eqr1a5.3}
		\mathrm{Rank} (P_1(z_0))= \mathrm{Rank} (P_2(z_0)).
		\end{equation}
		Since $z_0\in \R^d$ is arbitrary, it follows that}
	\begin{equation*}
		\begin{split}
			&K_1\text{ has no Rank-$1$ connections}\nn\\
			&\qd\iff \mathrm{Rank}(P_1(z))\geq 2\text{ for any }z\not =0\nn\\
			&\qd\overset{(\ref{eqr1a5.3})}{\iff} \mathrm{Rank}(P_2(z))\geq 2\text{ for any }z\not =0\nn\\
			&\qd\iff K_2\text{ has no Rank-$1$ connections}.
		\end{split}
	\end{equation*}
	
	Next, as $P_2(z)$ is obtained from $P_1(z)$, both viewed as $m\times n$ matrices with entries in $\R[z_1,\dots,z_d]$, by an elementary tow operation, we have that 
	\begin{equation}
	\label{eqr1a5.7}
	R_i(P_2(z))=\sum_{i'=1}^m c_{i,i'} R_{i'}(P_1(z)),
	\end{equation}
	\xxblue{where \pg{recall that} $R_i(A)$ denotes the $i$-th row of a matrix $A$}. It follows that
	\begin{eqnarray*}
		&~&M_{i_0,i_1}^{j_0,j_1}\lt(P_2(z)\rt)=\det\lt(\begin{array}{cc} \lt[P_2(z)\rt]_{i_0,j_0} & \lt[P_2(z)\rt]_{i_0,j_1} \\  \lt[P_2(z)\rt]_{i_1,j_0} &   \lt[P_2(z)\rt]_{i_1,j_1}\end{array}\rt)\nn\\
		&~&\qd\qd=\lt(\lt[P_2(z)\rt]_{i_0,j_0},\lt[P_2(z)\rt]_{i_0,j_1} \rt)\land \lt(\lt[P_2(z)\rt]_{i_1,j_0},\lt[P_2(z)\rt]_{i_1,j_1}\rt)\nn\\
		&~&\qd\qd= \xxred{\lt(\sum_{i'=1}^m c_{i_0,i'} \lt(\lt[P_1(z)\rt]_{i',j_0},\lt[P_1(z)\rt]_{i',j_1}\rt)\rt)\land 
			\lt(\sum_{i''=1}^m c_{i_1,i''} \lt(\lt[P_1(z)\rt]_{i'',j_0}, \lt[P_1(z)\rt]_{i'',j_1}\rt)\rt)}\nn\\
		&~&\qd\qd=\sum_{i',i''=1}^m c_{i_0,i'}c_{i_1,i''} M_{i',i''}^{j_0,j_1}(P_1(z)).
	\end{eqnarray*}
	This shows that all $2\times 2$ minors of $P_2(z)$ are inside the span of the $2\times 2$ minors of $P_1(z)$, both as subsets of $\R[z_1,\dots,z_d]$. Conversely, by (\ref{eqr1a5.7}) the rows of $P_1(z)$ can be represented as linear combinations of the rows of $P_2(z)$. Therefore, exactly the same argument shows the opposite inclusion in (\ref{eqr1a5}).
\end{proof}

To simplify notation, given $a, b\in \R^d$, we define $a\odot b\in\R[z_1,\dots,z_d]$ by
\begin{equation}
\label{tseq22}
a\odot b:=\lt(a\cdot z\rt)\lt(b\cdot z\rt).
\end{equation}
If $a\in\R^d$ and $W$ is a subspace of $\R^d$ we define
\begin{equation}
\label{tseq18}
a\odot W:=\lt\{a\odot w:w\in W\rt\}.
\end{equation}
Further given two subspaces $W, R\subset \R^d$ we define
\begin{equation}
\label{tseq23}
W\odot R:=\lt\{w\odot r: w\in W, r\in R\rt\}.
\end{equation}

\begin{a1}
	\label{l17}
	Let $d\in\{1,2,3\}$ and $K\subset M^{m\times n}$ be a $d$-dimensional subspace without Rank-$1$ connections. Assume that $K=P_K(\R^d)$ with $P_K$ represented by (\ref{eqr1a1}). If
	\begin{equation}
	\label{tseq7}
	\mathrm{dim}\lt(\mathrm{Span}\lt\{a_{i_0l}:l=1,2,\dots, n\rt\}\rt)= 1\text{ for some }i_0\in \lt\{1,2,\dots, m\rt\}
	\end{equation}
	or
	\begin{equation}
	\label{tseq8}
	\mathrm{dim}\lt(\mathrm{Span}\lt\{a_{lj_0}:l=1,2,\dots, m\rt\}\rt)=1\text{ for some }j_0\in \lt\{1,2,\dots, n\rt\},
	\end{equation}
	then there exists some $\beta\in\pg{\R^{q_0}\setminus\{0\}}$ such that (\ref{sseqc1}) is satisfied.
\end{a1}

\begin{proof}
	It is enough to establish (\ref{sseqc1}) assuming (\ref{tseq7}). Under the assumption (\ref{tseq8}), the conclusion follows in exactly the same way. Recall the definition (\ref{eqp103}). By performing row and column operations to $P_K(z)$ we have
	\begin{equation}
	\label{tseq9}
	P_K(z)\sim\wt{P}(z):=\lt(\begin{array}{cccc} \ti a_{11}\cdot z & 0 & \dots  & 0 \\
	\ti a_{21}\cdot z & \ti a_{22}\cdot z & \dots  & \ti a_{2n}\cdot z  \\
	\dots & \dots \\
	\ti a_{m1}\cdot z & \ti a_{m2}\cdot z & \dots  & \ti a_{mn}\cdot z
	\end{array}\rt).
	\end{equation}
	It is clear that $\mathrm{Span}\{M_k(\wt P(z)):k=1,\dots,q_0\}$ contains (recalling the notation (\ref{tseq22}))
	\begin{equation*}
	\lt\{\ti a_{11}\odot \ti a_{ij}: i\in \lt\{2,3,\dots, m\rt\},  j\in \lt\{2,3\dots, n\rt\} \rt\}.
	\end{equation*}
	Let 
	\begin{equation*}
	U_0:=\mathrm{Span}\lt\{\ti a_{ij}: i\in \lt\{2,3,\dots, m\rt\},  j\in \lt\{2,3,\dots, n\rt\}\rt\}.
	\end{equation*}
	It follows that (recalling the notation (\ref{tseq18}))
	\begin{equation}
	\label{eqp104}
	\ti a_{11}\odot U_0 \subset \mathrm{Span}\{M_k(\wt P(z)):k=1,\dots,q_0\}.
	\end{equation}
	
	Suppose $\mathrm{dim}(U_0)\leq d-1$ (when $d=1$, $U_0$ simply contains all scalar zeros), then there exists some $z_0\in \R^{d}$ such that $\ti a_{ij}\cdot z_0=0$ for all $i\in \lt\{2,3, \dots, m\rt\},  j\in \lt\{2,3, \dots, n\rt\}$ and thus from (\ref{tseq9}) we have $\mathrm{Rank}(\wt P(z_0))=1$. By Lemma \ref{r1lemma2}, since $K$ does not contain Rank-$1$ connections, we know that the subspace parametrized by $\wt P$ also does not contain Rank-$1$ connections. This contradicts the fact that $\mathrm{Rank}(\wt P(z_0))=1$. Hence we have $\mathrm{dim}(U_0)=d$, and thus $\ti a_{11}\in U_0$.
	It is then clear from (\ref{eqp104}) that $\ti a_{11}\odot\ti a_{11} \in \mathrm{Span}\{M_k(\wt P(z)):k=1,\dots,q_0\}$. As $\ti a_{11}$ is non-trivial, $\ti a_{11}\odot\ti a_{11}=\lt(\ti a_{11}\cdot z\rt)^2$ provides an element in $\mathrm{Span}\{M_k(\wt P(z)):k=1,\dots,q_0\}$ that is non-negative and non-trivial. By (\ref{eqr1a5}), we know that $\ti a_{11}\odot\ti a_{11} \in \mathrm{Span}\{M_k(P_K(z)):k=1,\dots,q_0\}$ and this establishes (\ref{sseqc1}).
\end{proof}

\begin{a1}
	\label{L18}
	Assume that $\min\{m,n\} =2$. Let $d\in\{2,3\}$ and $K\subset M^{m\times n}$ be a $d$-dimensional subspace without Rank-$1$ connections, then there exists some $\beta\in\pg{\R^{q_0}\setminus\{0\}}$ such that (\ref{sseqc1}) is satisfied.
\end{a1}

\begin{proof}
	Without loss of generality, we assume that $n=2$, and the case $m=2$ can be dealt with in an identical manner. Let $K=P_K(\R^d)$ with $P_K$ represented by (\ref{eqr1a1}). We claim that 
	\begin{equation*}
	\dim \lt(\mathrm{Span}\lt\{a_{ij}: i=1,\dots,m\rt\}\rt) =\red{d}\text{ for } j = 1,2.
	\end{equation*}
	If not, suppose without loss of generality that $\dim\lt(\mathrm{Span}\{a_{i1}: i=1,\dots,m\}\rt) \leq \red{d-1}$. Then there exists some $z_0\ne 0$ such that $a_{i1}\cdot z_0 = 0$ for all $i=1,\dots,m$, and hence $P_K(z_0)$ forms a Rank-$1$ direction in $K$, which is a contradiction. 
	
	By row operations \org{on} $P_K(z)$, we may without loss of generality assume that 
	\begin{equation*}
	\dim \lt(\mathrm{Span}\lt\{a_{i1}: \red{i=1,\dots, d}\rt\}\rt)=\red{d}. 
	\end{equation*}
	By further row operations \red{we} eliminate the remaining terms in the first column in $P_K(z)$ (viewed as a matrix with entries in $\R[z_1,\dots,z_d]$) \pg{and get}
	\begin{equation}
	P_K(z)\sim\wt P(z):=\lt(\begin{array}{cc} \ti a_{11}\cdot z & \ti a_{12}\cdot z  \nn\\
	\dots  & \dots  \nn\\
	\red{\ti a_{d 1}} \cdot z & \red{\ti a_{d 2}}\cdot z  \nn\\
	0 & \ti a_{\red{(d+1)} 2}\cdot z   \nn\\
	\dots &  \dots  \nn\\
	\end{array}\rt), 
	\end{equation}
	that is, $\red{\ti{a}_{i1}=0}$ for all $i\geq d+1$. We may assume that 
	\begin{equation}
	\label{eqp105}
	\ti a_{i2}=0 \text{ for all } i\geq d+1,
	\end{equation} 
	as otherwise there would be a row, say the $i_0$-th row with $i_0\geq d+1$, in $\wt P(z)$ such that (\ref{tseq7}) is satisfied for this row. Then we are done in this case by Lemmas \ref{l17} and \ref{r1lemma2}. So assuming (\ref{eqp105}), $\wt P(z)$ is isomorphic to a $d$-dimensional subspace in $M^{d\times 2}$. When $d=3$, by Proposition 4.4 in \cite{bhat}, all three-dimensional subspaces in $M^{3\times 2}$ must contain Rank-$1$ connections and thus this is a contradiction. When $d=2$, 
	$\red{\wt P(z)=\lt(\begin{array} {cc} \ti a_{11}\cdot z & \ti a_{12}\cdot z\\ \ti a_{21}\cdot z & \ti a_{22}\cdot z \end{array}\rt)}$ (\pg{up to an isomorphism}). By Lemma \ref{r1lemma2}, $\wt P(z)$ does not contain Rank-$1$ connections and hence $\det \wt P(z)\ne 0$ for all $z\ne 0$. Then (possibly after multiplying by $-1$) $\det \wt P(z)$ is a non-negative and non-trivial minor in $\mathrm{Span}\{M_k(\wt P(z)):k=1,\dots,q_0\}$, and we are done by (\ref{eqr1a5}).
\end{proof}	

\begin{a1}
	\label{l18}
	Let $d\in\{2,3\}$ and $K\subset M^{m\times n}$ be a $d$-dimensional subspace without Rank-$1$ connections. Assume that $K=P_K(\R^d)$ with $P_K$ represented by (\ref{eqr1a1}). If
	\begin{equation}
	\label{eqp100}
	\mathrm{dim}\lt(\mathrm{Span}\lt\{a_{i_0l}:l=1,2,\dots, n\rt\}\rt)= 2\text{ for some }i_0\in \lt\{1,2,\dots, m\rt\}
	\end{equation}
	or
	\begin{equation*}
	\mathrm{dim}\lt(\mathrm{Span}\lt\{a_{lj_0}:l=1,2,\dots, m\rt\}\rt)=2\text{ for some }j_0\in \lt\{1,2,\dots, n\rt\},
	\end{equation*}
	then there exists some $\beta\in\pg{\R^{q_0}\setminus\{0\}}$ such that (\ref{sseqc1}) is satisfied.
\end{a1}

\begin{proof}
	By Lemma \ref{L18}, we may assume that $\min\{m,n\}\geq 3$. Without loss of generality, we assume (\ref{eqp100}) with $i_0=1$ and $\dim\lt(\mathrm{Span}\lt\{a_{11},a_{12}\rt\}\rt)=2$. We perform column operations to $P_K(z)$ (as a matrix with entries in $\R[z_1,\dots,z_d]$) to eliminate the remaining terms in the first row to get
	\begin{equation}
	\label{tseq45}
	\red{P_K(z)}\sim \wt P(z):=\lt(\begin{array}{ccccc} \ti a_{11}\cdot z & \ti a_{12}\cdot z & 0 & \dots  & 0 \\
	\ti a_{21}\cdot z & \ti a_{22}\cdot z & \ti a_{23}\cdot z & \dots & \ti a_{2n}\cdot z  \\
	\dots & \dots \\
	\ti a_{m1}\cdot z & \ti a_{m2}\cdot z & \ti a_{m3}\cdot z & \dots & \ti a_{mn}\cdot z
	\end{array}\rt),
	\end{equation}
	where 
	\begin{equation*}
	\dim\lt(\mathrm{Span}\lt\{\ti a_{11},\ti a_{12}\rt\}\rt)=2. 
	\end{equation*}
	If $\ti a_{ij}=0$ for all $i\geq 2$ and $j\geq 3$, then $\wt P(z)$ is isomorphic to a $d$-dimensional subspace in $M^{m\times 2}$, and thus we are done by Lemmas \ref{L18} and \ref{r1lemma2}. So in the following we assume that there exists some $j_0\in\{3,\dots,n\}$ such that the $j_0$-th column in $\red{\wt P(z)}$ is non-trivial. By Lemma \ref{l17} we may assume that $\dim\lt(\mathrm{Span}\lt\{\ti a_{2j_0},\dots,\ti a_{mj_0}\rt\}\rt)\geq2$ \red{since \pg{otherwise} there would be a column of $P_k(z)$ for which 
		(\ref{tseq8}) \pg{holds true}}. We define
	\begin{equation*}
	W_1:=\mathrm{Span}\{\ti a_{11}, \ti a_{12}\} \text{ and } U_1:= \mathrm{Span}\lt\{\ti a_{2j_0},\dots,\ti a_{mj_0}\rt\}.
	\end{equation*}
	Note that $\dim(W_1)+\dim(U_1) = 4>d$ and hence there exists some non-trivial $\psi\in \R^d$ such that $\psi\in W_1\cap U_1$. In particular, recalling the notation (\ref{tseq22}) and (\ref{tseq23}), we have $\psi\odot\psi\in W_1\odot U_1$. It is clear from (\ref{tseq45}) that $W_1 \odot U_1 \subset \mathrm{Span}\{M_k(\wt P(z)):k=1,\dots,q_0\}$ and hence $\psi\odot\psi\in \mathrm{Span}\{M_k(\wt P(z)):k=1,\dots,q_0\}$ is non-negative and non-trivial. Finally the conclusion of the lemma follows from (\ref{eqr1a5}).
\end{proof}

\begin{a1}
	\label{redlemma}
	Let $K\subset M^{m\times n}$ be a three-dimensional subspace without Rank-$1$ connections. Assume that $K=P_K(\R^d)$ with $P_K$ represented by (\ref{eqr1a1}). If
	\begin{equation}
	\label{tseq40}
	\mathrm{dim}\lt(\mathrm{Span}\lt\{a_{i_0 l}:l=1,2,\dots, n\rt\}\rt)= 3\text{ for some }i_0\in \lt\{1,2,\dots, m\rt\}
	\end{equation}
	or
	\begin{equation*}
	\mathrm{dim}\lt(\mathrm{Span}\lt\{a_{l j_0}:l=1,2,\dots, m\rt\}\rt)= 3\text{ for some }j_0\in \lt\{1,2,\dots, n\rt\},
	\end{equation*}
	then there exists some $\beta\in\pg{\R^{q_0}\setminus\{0\}}$ such that (\ref{sseqc1}) is satisfied.
\end{a1}

\begin{proof} 
	As at the beginning of the proof of Lemma \ref{l18}, we assume without loss of generality (\ref{tseq40}) with $i_0=1$ and find $P_K(z)\sim \wt P(z)$   \org{where $\lt[\wt P(z)\rt]_{ij}= \ti a_{ij}\cdot z$} and $\{\ti a_{ij}\}\subset \R^3$ satisfies
	\begin{equation}
	\label{tseq70}
	\dim\lt(\mathrm{Span}\lt\{\ti a_{11},\ti a_{12}, \ti a_{13}\rt\}\rt)=3\text{ and }\ti a_{1j}=0\text{ for all }j=\blue{4},\dots, n,
	\end{equation}
	provided $n\geq 4$. If there exists $i_0\in \lt\{2,3,\dots, m\rt\}, j_0\in \lt\{4,5\dots, n\rt\}$ with $\ti a_{i_0 j_0}\not=0$, \red{then since $\ti{a}_{i_0 j_0}\odot \mathrm{Span}\lt\{\ti a_{11},\ti a_{12}, \ti a_{13}\rt\}\subset \mathrm{Span}\{M_k(\wt P(z)):k=1,\dots,q_0\}$},  it follows from (\ref{tseq70}) that $\ti a_{i_0j_0}\odot \ti a_{i_0j_0}\in \mathrm{Span}\{M_k(\wt P(z)):k=1,\dots,q_0\}$ and we are done. So in the following we assume that
	\begin{equation*}
	\ti a_{ij}=0\text{ for all }i\in \lt\{2,3,\dots, m\rt\}, j\in \lt\{4,5\dots, n\rt\},
	\end{equation*}
	and hence $\wt P(z)$ is isomorphic to a three-dimensional subspace in $M^{m\times 3}$. Note that this also includes the case $n=3$.
	
	By Lemmas \ref{l17} and \ref{l18}, we only have to deal with the case when 
	\begin{equation*}
	\mathrm{dim}\lt(\mathrm{Span}\{\ti a_{i1},\ti a_{i2},\ti a_{i3}\}\rt) = 3 \text{ for all } i
	\end{equation*}
	and
	\begin{equation*}
	\mathrm{dim}\lt(\mathrm{Span}\{\ti a_{ij}: i=1,2,\dots,m\}\rt) = 3 \text{ for all } j.
	\end{equation*}
	By Lemma \ref{L18}, we may assume that there are at least three non-zero rows in $\wt P(z)$. If $\wt P(z)$ has more than three non-zero rows, we can perform row operations to eliminate all but three entries in the first column in $\wt P(z)$ and obtain $\wt P(z)\sim \hat P(z):=\lt(\hat a_{ij}\cdot z\rt)$, where $\hat P(z)$ (up to an isomorphism) is an $m\times 3$ matrix with entries in $\R[z_1,\dots,z_3]$ and $\hat a_{i1}=0$ for all $i\geq 4$. If there exists $i_0\geq 4$ such that the $i_0$-th row of $\hat P(z)$ is non-trivial, then (noting that $\hat P(z)$ has only three non-zero columns) \pb{Lemma} \ref{l18} \red{can be applied} to give a non-negative and non-trivial element in $\mathrm{Span}\{M_k(\hat P(z)):k=1,\dots,q_0\}$, which \org{by Lemma \pgr{\ref{r1lemma2}}}  also belongs to $\mathrm{Span}\{M_k(P_K(z)):k=1,\dots,q_0\}$. Hence the only remaining case is where $\wt P(z)$ is isomorphic to a three-dimensional subspace in $M^{3\times 3}$. We prove this case in Lemma \ref{redlem2}. This concludes the proof of Lemma \ref{redlemma}.
\end{proof}

\begin{a1}
	\label{redlem2}
	Let $K\subset M^{3\times 3}$ be a \blue{three-dimensional} subspace without Rank-$1$ connections, then there exists some $\beta\in \pg{\R^{q_0}\setminus\{0\}}$ such 
	that (\ref{sseqc1}) holds true. 
\end{a1}

\begin{proof} 
	We assume that $K=P_K(\R^3)$ for $P_K(z)$ given in (\ref{eqr1a1}). By Lemmas \ref{l17} and \ref{l18}, we may assume that
	\begin{equation*}
	\mathrm{dim}\lt(\mathrm{Span}\{a_{1j},a_{2j},a_{3j}\}\rt) = 3 \text{ for all } j.
	\end{equation*}
	Thus we can perform row operations to \red{clean up} the first column of $P_K(z)$ to get
	\begin{equation*}
	P_K(z)\sim\wt P(z):=\lt(\begin{array}{ccc} e_1\cdot z & \ti{a}_{12}\cdot z & \ti{a}_{13}\cdot z\\
	e_2\cdot z & \ti{a}_{22}\cdot z& \ti{a}_{23}\cdot z \\
	e_3\cdot z & \ti{a}_{32}\cdot z & \ti{a}_{33}\cdot z \\
	\end{array}\rt).
	\end{equation*}
	Again by Lemmas \ref{l17} and \ref{l18} we may assume that $\lt\{e_1, \ti{a}_{12}, \ti{a}_{13}\rt\}$ is linearly independent. So we can perform column operations to $\wt P(z)$, using the first column to eliminate the $e_1$ component in $\ti{a}_{12}$ and $\ti{a}_{13}$, and then performing column operations to the second and third columns to find
	\begin{equation}
	\label{tseq231}
	\wt P(z)\sim\hat P(z):=\lt(\begin{array}{ccc} e_1\cdot z & e_2\cdot z & e_3\cdot z \\
	e_2\cdot z & \hat{a}_{22}\cdot z& \hat{a}_{23}\cdot z \\
	e_3\cdot z & \hat{a}_{32}\cdot z &\hat{a}_{33}\cdot z \\
	\end{array}\rt).
	\end{equation}
	
	Next we will show that if $\hat{a}_{32}\ne\hat{a}_{23}$ then (\ref{sseqc1}) holds true. To this end, we examine the minors of $\hat P(z)$. Note that (recalling (\ref{eqr1a70}))
	\begin{equation*}
		M^{1,3}_{1,2}(\hat P(z))=z_1 (\hat{a}_{23}\cdot z)-z_3 z_2,\qd M^{1,2}_{1,3}(\hat P(z))=z_1 (\hat{a}_{32}\cdot z)-z_2 z_3, 
	\end{equation*}
	\begin{equation*}
		M^{2,3}_{1,2}(\hat P(z))=z_2 (\hat{a}_{23}\cdot z)-z_3  (\hat{a}_{22}\cdot z),\qd M^{2,3}_{1,3}(\hat P(z))=z_2 (\hat{a}_{33}\cdot z)-z_3  (\hat{a}_{32}\cdot z),
	\end{equation*}
	\begin{equation*}
		M^{1,3}_{2,3}(\hat P(z))=z_2 (\hat{a}_{33}\cdot z)-z_3  (\hat{a}_{23}\cdot z),\qd M^{1,2}_{2,3}(\hat P(z))=z_2 (\hat{a}_{32}\cdot z)-z_3  (\hat{a}_{22}\cdot z).
	\end{equation*}
	So letting $b=(b_1,b_2,b_3)=\hat{a}_{23}-\hat{a}_{32}$ we have
	\begin{equation*}
		M^{1,3}_{1,2}(\hat P(z))-M^{1,2}_{1,3}(\hat P(z))=z_1\lt((\hat{a}_{23}-\hat{a}_{32} )\cdot z\rt)=z_1 \lt(b\cdot z\rt),
	\end{equation*}
	\begin{equation*}
		M^{2,3}_{1,2}(\hat P(z))-M^{1,2}_{2,3}(\hat P(z))=z_2\lt((\hat{a}_{23}-\hat{a}_{32} )\cdot z\rt)=z_2 \lt(b\cdot z\rt)
	\end{equation*}
	and 
	\begin{equation*}
		M^{2,3}_{1,3}(\hat P(z))-M^{1,3}_{2,3}(\hat P(z))=z_3\lt((\hat{a}_{23}-\hat{a}_{32} )\cdot z\rt)=z_3 \lt(b\cdot z\rt).
	\end{equation*}
	Thus 
	\begin{eqnarray*}
		&~&b_1 \lt(M^{1,3}_{1,2}(\hat P(z))-M^{1,2}_{1,3}(\hat P(z))\rt)+b_2\lt( M^{2,3}_{1,2}(\hat P(z))-M^{1,2}_{2,3}(\hat P(z))\rt)\nn\\
		&~&\qd\qd\qd+b_3\lt(M^{2,3}_{1,3}(\hat P(z))-M^{1,3}_{2,3}(\hat P(z)) \rt)=(b\cdot z)^2\nn
	\end{eqnarray*}
	and we have a non-negative and non-trivial element in $\mathrm{Span}\lt\{M_k(\hat P(z)):k=1,\dots,q_0\rt\}$ as $b\ne 0$. So we are done by (\ref{eqr1a5}).
	
	Finally, if $\hat{a}_{32}=\hat{a}_{23}$, then from (\ref{tseq231}), $\hat P(z)$ would define a three-dimensional subspace in $M^{3\times 3}_{sym}$ without Rank-$1$ connections. By Lemma \ref{svlem} (from Appendix \ref{AP1}), any three-dimensional subspace in $M^{3\times 3}_{sym}$ must contain a Rank-$1$ connection, which is a contradiction \red{by Lemma \ref{r1lemma2}}. This completes the proof of Lemma \ref{redlem2}.
\end{proof}

\begin{proof}[Proof of Theorem \ref{l16}]
	If $K$ is a three-dimensional subspace, then one of the assumptions in Lemmas \ref{l17}, \ref{l18} and \ref{redlemma} is trivially satisfied and hence the conclusion follows from these lemmas. Similarly, for two-dimensional subspaces, the conclusion follows from Lemmas \ref{l17} and \ref{l18} and for one-dimensional subspaces the conclusion follows from Lemma \ref{l17}.
\end{proof}

We conclude this section with the proof of Theorem \ref{CC1}. The necessity part is trivial. The sufficiency part makes use of Theorem \ref{l16} and the ideas are very similar to those in the proof of the sufficiency part in Theorem \ref{C1}.
\begin{proof}[Proof of Theorem \ref{CC1}]
	Suppose that $K$ has Rank-$1$ connections, then there exists a non-trivial $A\in K$ such that $\mathrm{Rank}(A)=1$. Let $\mu:=\frac{1}{2} \delta_{A}+\frac{1}{2}\delta_{-A}$. Note that for any minor $M_{\pg{k}}$ of $M^{m\times n}$, $M_{\pg{k}}(A)=M_{\pg{k}}(-A)=M_{\pg{k}}(\frac{A-A}{2})=0$. It follows that 
	$$
	\int_{K} M_{\pg{k}}(X) d\mu(X)=\frac{1}{2}M_{\pg{k}}(A)+\frac{1}{2}M_{\pg{k}}(-A)=0=M_{\pg{k}}\lt(\frac{A-A}{2}\rt)=M_{\pg{k}}\lt(\int_K X d\mu(X)\rt).
	$$
	Hence $\mu\in \mathcal{M}^{pc}(K)$ and $\mathcal{M}^{pc}(K)$ contains non-trivial measures.
	
	Next suppose that $K$ has no Rank-$1$ connections, and we show that $\mathcal{M}^{pc}(K)$ consists of Dirac measures. We assume that $K$ is a three-dimensional subspace and provide the detailed proof, part of which can be used to prove the cases for lower dimensional subspaces. 
	
	We first show that $\MI^{pc}_{K}(0)$ is trivial, where
	\begin{equation*}
	\MI^{pc}_{K}(0):=\lt\{\mu\in \MI^{pc}(K): \overline{\mu}=0\rt\}.
	\end{equation*}
	To this end, we apply Theorem \ref{l16} to the subspace $K$ to find $\beta\in\pg{\R^{q_0}\setminus\{0\}}$ such that
	\begin{equation}
	\label{eqp102}
	\sum_{k=1}^{q_0} \beta_k \pg{M}_k(X)\geq 0\text{ for all }X\in K\text{ and }\sum_{k=1}^{q_0} \beta_k \pg{M}_k\not\equiv 0\text{ on }K.
	\end{equation}
	Let $\mu\in \MI^{pc}_{K}(0)$. It follows from the definition of $\MI^{pc}_{K}(0)$ that
	\red{$\int_{K} M_k(X) d\mu(X)=M_k\lt(\overline{\mu}\rt)=0$ and hence 
		$\int_{K}\sum_{k=1}^{q_0} \beta_k \pg{M}_k(X)d\mu(X) = 0$}. By (\ref{eqp102}), it is clear that $\mathrm{Spt}(\mu)\subset K_1$ where
	\begin{equation*}
	K_1:=\lt\{X\in K: \sum_{k=1}^{q_0} \beta_k \pg{M}_k(X)=0\rt\}.
	\end{equation*}
	\pg{We claim that $K_1$ is a subspace of $K$. Since $\mathrm{dim}(K)=3$, there exists a linear isomorphism $\sigma:\R^3\rightarrow K$. Define $f(z):=\sum_{k=1}^{q_0} \beta_k \pg{M}_k(\sigma(z))$ which is a homogeneous quadratic function as all $M_k$'s are $2\times 2$ minors. Because of (\ref{eqp102}), $f$ is convex on $\R^3$. It follows that its zero set $\sigma^{-1}(K_1)$ is a convex cone, which is a subspace in $\R^3$, and hence $K_1$ is also a subspace of $K$ as $\sigma$ is a linear isomorphism.} Further since $\sum_{k=1}^{q_0} \beta_k \pg{M}_k(X)$ is non-trivial, we know that $K_1$ is a proper subspace of $K$. Now $\mathrm{Spt}(\mu)\subset K_1$ and $\mathrm{dim}\lt(K_1\rt)\leq 2$. Since $K$ has no Rank-$1$ connections, the same holds for $K_1$. Repeating the above arguments using Theorem \ref{l16} at most three times, we conclude that $\mu=\delta_{0}$ and hence $\MI^{pc}_{K}(0)$ is trivial.

	Now let $\mu\in \MI^{pc}(K)$ and $\overline{X}:=\int_K X\,d\mu(X)$. Define the translation $\PI^{\ol{X}}:M^{m\times n}\rightarrow M^{m\times n}$ by $\PI^{\ol{X}}(X):=X-\ol{X}$. Letting $\nu:=\lt(\PI^{\ol{X}}\rt)_\sharp\mu$, i.e., the push forward of the measure $\mu$ under the mapping $\PI^{\ol{X}}$, we claim that $\nu\in\MI^{pc}_K(0)$. First we have
	\begin{equation}
	\label{eqp106}
	\int_K X\,d\nu(X) = \int_{K}\lt(X-\ol{X}\rt)\,d\mu(X) = 0.
	\end{equation}
	Next recall that, given $2\times2$ matrices $A$ and $B$, \cblue{we have} that
	\begin{equation}\label{eqp201}
	\det(A-B) = \det(A) - A:\mathrm{Cof}(B) + \det(B).
	\end{equation}
	It follows that 
	\begin{equation*}
	\begin{split}
	&\int_{K} M_k(X)\,d\nu(X) = \int_{K} M_k\lt(X-\ol{X}\rt)\,d\mu(X)\\
	&\qd\qd\qd \overset{(\ref{eqp201})}{=}\int_K M_k(X)\,d\mu(X) - \int_{K} X:\mathrm{Cof}(\ol{X})\,d\mu(X) + M_k(\ol{X})\\
	&\qd\qd\qd \overset{\mu\in \MI^{pc}(K)}{=} M_k(\ol{X}) - \ol{X}:\mathrm{Cof}(\ol{X}) + M_k(\ol{X})\\
	&\qd\qd\qd\overset{(\ref{eqp201})}{=} M_k(\ol{X}-\ol{X}) = 0 \overset{(\ref{eqp106})}{=} M_k\lt(\int_K X\,d\nu(X)\rt)
	\end{split}
	\end{equation*}
	for all $k=1,\dots,q_0$. Hence we have established that $\nu\in\MI^{pc}_K(0)$, and thus $\nu=\delta_{0}$. It follows immediately that $\mu=\delta_{\ol{X}}$ and therefore $\MI^{pc}(K)$ consists of Dirac measures.
\end{proof}

%
%

\section{Proof of \cblue{Theorem} \ref{C17}} \label{t4}

\cblue{We first recall the notion of Grassmannian which is needed in the discussions of this section. Let $p, k$ be fixed integers with $p\geq 0$ and $0\leq k\leq p$. We denote by $G(k,p)$ the set of all $k$-dimensional subspaces of $\R^p$, and it is called the Grassmannian of $k$-dimensional subspaces of $\R^p$. We have the following property regarding $G(k,p)$, whose proof can be found, for example, in \cite{mas}:}

\begin{a1}
	\label{graslem1}
	\cblue{The Grassmannian $G(k,p)$ is a real analytic, compact and connected manifold of dimension $k(p-k)$. }
\end{a1}

\cblue{One can view $G(k,p)$ as a differentiable manifold in the following way. We fix a pair of transversal subspaces $(W_0,W_1)$ of $\R^p$, i.e., $W_0\cap W_1 = \{0\}$, where $\dim(W_0)=k$ and $\dim(W_1)=p-k$. Then one can view elements in $G^0(k,p,W_1)$ as the graphs of linear maps from $W_0$ to $W_1$, where
	\begin{equation*}
		G^0(k,p,W_1):=\lt\{V\in G(k,p):V\cap W_1 =\lt\{0\rt\}\rt\}.
	\end{equation*}
	Specifically,} we pick a basis $\lt\{a_1,a_2,\dots, a_k\rt\}$ for $W_0$ and a basis $\lt\{b_1,b_2,\dots, b_{p-k}\rt\}$ for $W_1$. We identify $\R^{k(p-k)}$ with $M^{\cblue{(p-k)\times k}}$ in the obvious way, \cblue{i.e., identify $x\in\R^{k(p-k)}$ with $A_{x}\in M^{(p-k)\times k}$ where $[A_x]_{ij}=\xgreen{[x]_{(i-1)k+j}}$}.  For each $\cblue{A}\in   \R^{k(p-k)}\simeq M^{\cblue{(p-k)\times k}}$, let $T_A:W_0\rightarrow W_1$ be the linear map defined by $A$ and the 
choices of bases $\lt\{a_1,a_2,\dots, a_{k}\rt\}$, $\lt\{b_1,b_2,\dots, b_{p-k}\rt\}$. We 
define 
\begin{equation}
\label{graseq1}
\phi_{W_0,W_1}(\cblue{A}):=\lt\{v+T_A(v):v\in W_0\rt\}.
\end{equation}
Note that the mapping $\phi_{W_{\cblue{0}},W_{\cblue{1}}}$ is one to one from $\R^{k(p-k)}$ onto $G^0(k,p,W_1)$. \cblue{Hence it defines a chart on $G(k,p)$ that covers $G^0(k,p,W_1)$}. As noted in Remark \cblue{2.2.4} in \cite{mas}, the charts defined by (\ref{graseq1}) actually form a real analytic \cblue{atlas for $G(k,p)$. Further, it is shown in Corollary 2.4.3 in \cite{mas} that $G(k,p)$ is compact and connected.}

%
%
%

\begin{deff}
	\label{grasdef1}
	We say that a property holds \bf  generically \rm \em  for \cblue{$k$-dimensional subspaces of $\R^p$} if there exist
	finitely many smooth manifolds $\Gamma_1, \Gamma_2, \dots, \Gamma_{\cblue{r}}$ \xgreen{in $\R^{k(p-k)}$} of dimension less \cblue{than} $k(p-k)$ and Lipschitz mappings $P_j:\Gamma_j\rightarrow G(k,p)$ for $j=1,2,\dots, \cblue{r}$ such that the property holds true 
	for every $V\in G(k,p)\backslash \lt(\bigcup_{j=1}^r P_j(\Gamma_j) \rt)$.
\end{deff}

%
%
%

\begin{remark}
	We let $G(k,M^{m\times n})$ denote the space of $k$-dimensional subspaces in $M^{m\times n}$. In an obvious way we can uniquely identify any $V\in G(k,M^{m\times n})$ with some $\cblue{W}\in G(k,mn)$. For this reason we will not distinguish between $G(k,M^{m\times n})$ and $G(k,mn)$. 
\end{remark}

%
%
%

\cblue{Let $k, m, n$ be positive integers with $k\leq mn$, and $(W_0,W_1)$ be a pair of transversal subspaces of $\R^{mn}$ with $\dim(W_0)=k$ and $\dim(W_1)=mn-k$. Further let $T: W_0\rightarrow W_1$ be a linear mapping. We fix a basis $\mathcal{B}_{0}=\lt\{a_1, a_2, \dots, a_k\rt\}$ for $W_0$. Recall that $M_1,M_2,\dots, M_{q_0}$ denote all $2\times 2$ minors in $M^{m\times n}$. 
	With the linear mapping $T$ and the basis $\mathcal{B}_{0}$ we can define the set of quadratics $Q_1, Q_2,\dots, Q_{q_0}$ on $\R^k$ by 
	\begin{equation}
	\label{graseq23}
	Q_j(y):=M_j\lt(\sum_{l=1}^{k} y_l \lt(a_l+T(a_l)\rt)\rt)\text{ for }j=1,2,\dots, q_0. 
	\end{equation}
	Then for each $Q_j$, there exists a unique $X_j\in M^{k\times k}_{sym}$ that represents $Q_j$. We need the following auxiliary lemma. }

\begin{a1} 
	\label{graslem1.5}
	
	\cblue{Let $k, m, n$ be positive integers with $k\leq mn$.} Suppose that for $\varpi=1,2$, $(W^{\varpi}_0, W^{\varpi}_1)$ is a pair of transversal subspaces of $\R^{mn}$ with $\mathrm{dim}(W^{\varpi}_0)=k$, $\mathrm{dim}(W^{\varpi}_1)=nm-k$, and 
	$T^{\varpi}:W^{\varpi}_0\rightarrow W^{\varpi}_1$ is a linear mapping. Suppose also that for some $V\in G(k,mn)$ we have that 
	\begin{equation}
	\label{graseq22}
	\lt\{v+T^1(v):v\in W^1_0\rt\}=V=\lt\{v+T^2(v):v\in W^2_0\rt\}.
	\end{equation}
	Let $\mathcal{B}_0^{\varpi}=\lt\{a^{\varpi}_1, a^{\varpi}_2, \dots, a^{\varpi}_k\rt\}$ be a basis of $W^{\varpi}_0$ for $\varpi=1,2$. \cblue{We denote by $X_j^{\varpi}\in  M^{k\times k}_{sym}$ the symmetric matrix that represents the quadratic $Q_j^{\varpi}$ given by (\ref{graseq23}) with respect to the linear mapping $T^{\varpi}$ and the basis $\mathcal{B}_0^{\varpi}$.}  Then we have
	\begin{equation}
	\label{graseq230}
	\mathrm{Span}\lt\{X^1_1,X^1_2,\dots, X^1_{q_0}\rt\}=M^{k\times k}_{sym}\iff \mathrm{Span}\lt\{X^2_1,X^2_2,\dots, X^2_{q_0}\rt\}=M^{k\times k}_{sym}.
	\end{equation}
	
\end{a1}

\begin{proof} 
	We begin by showing that for $\varpi=1,2$, 
	\begin{equation}
	\label{graseq41}
	\mathcal{B}^{\varpi}:=\lt\{a^{\varpi}_l+T^{\varpi}(a^{\varpi}_l):l=1,2,\dots, k\rt\}\text{ forms a basis of }V.
	\end{equation}
	From \cblue{(\ref{graseq22}) and the fact that $\mathcal{B}_0^{\varpi}$ is a basis of $W_0^{\varpi}$,} it is immediate that $\mathcal{B}^{\varpi}$ spans 
	$V$. Assume $\sum_{l=1}^k \lm_l\lt(a^{\varpi}_l+T^{\varpi}(a_l^{\varpi})\rt)=0$. As 
	$\sum_{l=1}^k \lm_l a^{\varpi}_l\in W^{\cblue{\varpi}}_0$, $\sum_{l=1}^k \lm_l T^{\varpi}(a^{\varpi}_l)\in W^{\cblue{\varpi}}_1$, \cblue{and $W^{\cblue{\varpi}}_0$ and $W^{\cblue{\varpi}}_{1}$ are transversal, it follows that
		$\sum_{l=1}^k \lm_l a^{\varpi}_l=0$ and $\sum_{l=1}^k \lm_l T^{\varpi}(a^{\varpi}_l)=0$.} Hence 
	$\lm_l=0$ for all $l=1,2,\dots,k$ and thus (\ref{graseq41}) is established.

	Let $A\in M^{k\times k}$ denote the change of basis matrix between the two bases $\mathcal{B}^1$ and $\mathcal{B}^2$ of $V$, i.e., letting $\alpha_{ij}\in \R$ denote the 
	$(i,j)$ entry of $A$, we have that 
	\begin{equation*}
		a^1_i+T^1(a^1_i)=\sum_{l=1}^k \alpha_{il}(a^2_l+T^2(a^2_l))\text{ for }i=1,2,\dots, k.
	\end{equation*}
	So for any $j\in \lt\{1,2,\dots, q_0\rt\}$, we have
	\begin{equation*}
		\begin{split}
			y^T X^1_j y&=M_j \lt( \sum_{i=1}^k y_i\lt(a^1_i+T^1(a^1_i)\rt)   \rt)\\
			&=M_{j}\lt(\sum_{i=1}^k y_i \lt(\sum_{l=1}^k \alpha_{il} \lt(a^2_l+T^2(a^2_l)\rt)  \rt)    \rt)\\
			&=M_{j}\lt(\sum_{l=1}^k \lt(\sum_{i=1}^k y_i \alpha_{i l} \rt)\lt( a^2_l+T^2(a^2_l) \rt) \rt)\\
			&=M_{j}\lt(\sum_{l=1}^k \lt[\cblue{A^T}y\rt]_l\lt(a^2_l+T^2(a^2_l) \rt) \rt)\\
			&= \lt[\cblue{A^T} y\rt]^T X^2_j \lt[\cblue{A^T} y\rt]= y^T \cblue{A} X^2_{j} \cblue{A^T} y,
		\end{split}
	\end{equation*}
	and thus
	\begin{equation}
	\label{graseq24.77}
	X^1_j=\cblue{A} X^2_{j} \cblue{A^T} \text{ for }j=1,2,\dots, q_0. 
	\end{equation}

	Suppose $\mathrm{Span}\lt\{X^1_1,X^1_2,\dots, X^1_{q_0}\rt\}=M^{k\times k}_{sym}$. For any $S\in M^{k\times k}_{sym}$, note that $\cblue{A S A^T}\in M^{k\times k}_{sym}$. Therefore we have $\sum_{j=1}^k \lm_j X^1_j=\cblue{A S A^T}$ for some 
	$\lm_1,\lm_2, \dots, \lm_k\in \R$. It follows from (\ref{graseq24.77}) that 
	\cblue{$\sum_{j=1}^k \lm_j X^2_j=\sum_{j=1}^k \lm_j A^{-1} X_j^1 (A^T)^{-1}=S$}, and thus $M^{k\times k}_{sym}=\mathrm{Span}\lt\{X^2_1,X^2_2, \dots, X^2_{\cblue{q_0}}\rt\}$. Exactly the same argument shows the other implication in (\ref{graseq230}). 
\end{proof}

%
%
%

\begin{a1} 
	\label{graslem2}
	
	\cblue{Let $k, m, n$ be positive integers with $m,n\geq 2$ and} $k\leq \frac{1}{2}\min\lt\{m, n\rt\}$, and $M_1,M_2,\dots, M_{q_0}$ denote all the $2\times 2$ minors of $M^{m\times n}$. Generically for $V\in G(k,M^{m\times n})$ (in the sense of Definition \ref{grasdef1}) there exists $\beta\in S^{q_0-1}$ such that 
	\begin{equation}
	\label{graseq3}
	\sum_{j=1}^{q_0} \beta_j M_j(X)>0\text{ for all }X\in V\backslash \lt\{0\rt\}. 
	\end{equation}
\end{a1}

\begin{proof}
	\cblue{By Lemma \ref{graslem1}}, $G(k,M^{m\times n})$ is a $k(mn-k)$-\cblue{dimensional} compact manifold. Therefore we can find finitely many charts \cblue{of the form (\ref{graseq1})} whose images 
	cover $G(k,M^{m\times n})$. Formally we can find finitely many pairs of \cblue{transversal} subspaces 
	\begin{equation*}
		\lt\{(W^1_0,W^1_1), (W^2_0,W^2_1), \dots, (W^{\cblue{p}_0}_0,W^{\cblue{p}_0}_1)\rt\}
	\end{equation*}
	\cblue{with $\dim(W_0^i)=k$ and $\dim(W^i_1)=mn-k$} such that 
	\begin{equation}
	\label{graseq39}
	G(k,M^{m\times n})\subset \bigcup_{i=1}^{p_0} \phi_{W^i_0,W^i_1}(\R^{k(mn-k)}). 
	\end{equation}
	
	For each $i\in \{1,2,\dots,p_0\}$, we fix a basis $\mathcal{B}_0^i=\{a_1^i,\dots,a_k^i\}$ for $W_0^i$ and a basis $\mathcal{B}_1^i=\{b_1^i,\dots,b_{mn-k}^i\}$ for $W_1^i$. Given $x\in \R^{k(mn-k)}$, define $A_x\in M^{\cblue{(mn-k)\times k}}$ to be the matrix \cblue{with $[A_x]_{ij}=\xgreen{[x]_{(i-1)k+j}}$}, and let $T_x^{\cblue{i}}:W^{i}_0\rightarrow W^{i}_1$ be the linear mapping defined 
	from $A_x$ given the bases $\mathcal{B}_0^i$ and $\mathcal{B}_1^i$. For each $x\in \R^{k(mn-k)}$, exactly the same arguments used to establish (\ref{graseq41}) show that the 
	set 
	\begin{equation}
	\label{graseq7}
	\cblue{\mathcal{B}_x^i:=\lt\{a_l^i+T_x^i(a_l^i):l=1,2,\dots, k\rt\}}
	\end{equation}
	forms a basis of the subspace $\phi_{W_0^{i}, W_1^{i}}(x)=\lt\{(v,T_x^i(v)):v\in W^{i}_0\rt\}$. For any $j\in \lt\{1,2,\dots, q_0\rt\}$, \cblue{as in (\ref{graseq23})},  the mapping 
	\begin{equation}
	\label{graseq9}
	\cblue{Q_{j,x}^i}(y):=M_j\lt(\sum_{l=1}^k y_l \lt(a_l^i+T_x^i(a_l^i)\rt) \rt)
	\end{equation}
	defines a quadratic mapping on $\R^k$ and hence can be represented by a matrix $X^i_{\cblue{j,x}}\in M^{k\times k}_{sym}$. \cblue{Given $i\in \{1,2,\dots,p_0\}$, recall that we have fixed the bases $\mathcal{B}_0^i$ and $\mathcal{B}_1^i$. Thus, each entry of the $m\times n$ matrix $a_l^i+T_x^i(a_l^i)$ is either linear in $x$ or constant. As $M_j$ is a $2\times 2$ minor, it follows that the coefficients in the quadratic $Q_{j,x}^i(y)$ are polynomials (of degree less than or equal to two) of $x$, and so are all the entries of the matrix $X^i_{j,x}$.} \nl
	
	\em Step 1. \rm For each $i\in \{1,2,\dots,p_0\}$, we show that there exists a \cblue{polynomial} function 
	$\Lambda^{i}:\R^{k(mn-k)}\rightarrow \R$ such that
	\begin{equation}
	\label{graseq34}
	\mathrm{Span}\lt\{X^i_{1,x}, X^i_{2,x}, \dots, X^i_{q_0,x}\rt\}=M^{k\times k}_{sym}\text{ for any }x\in \R^{k(mn-k)}\backslash \lt\{x:\Lambda^{i}(x)=0 \rt\}.
	\end{equation}
	
	\em Proof of Step 1. \rm \cblue{First note that, since $X^i_{j,x}$ is a symmetric $k\times k$ matrix, it can be uniquely identified with a vector $v^i_{j,x}\in \R^{\frac{k(k+1)}{2}}$. Then it is clear that
		\begin{equation}\label{graseq101}
		\mathrm{Span}\lt\{X^i_{1,x}, X^i_{2,x}, \dots, X^i_{q_0,x}\rt\}=M^{k\times k}_{sym} \Longleftrightarrow \mathrm{Span}\lt\{v^i_{1,x}, v^i_{2,x}, \dots, v^i_{q_0,x}\rt\}=\R^{\frac{k(k+1)}{2}}.
		\end{equation}
		Recall that $q_0$ is the number of $2\times 2$ minors in $M^{m\times n}$, and therefore $q_0 = \frac{m(m-1)}{2}\, \frac{n(n-1)}{2}$. Without loss of generality, we may assume $m\leq n$, and by assumption, we have $k\leq \frac{m}{2}$. In particular, $m\geq 2k \geq k+1$. \ared{In order 
			to apply Lemma \ref{grasaux1} later in the proof we observe the following inequality} 
		\begin{equation}
		\label{graseq100}
		q_0 = \frac{m(m-1)}{2}\, \frac{n(n-1)}{2} \geq \frac{m^2(m-1)^2}{4}\geq \frac{(2k)(k+1)}{4} (m-1)^2 \geq \frac{k(k+1)}{2}.
		\end{equation}}
	
	Now, viewing $v^i_{j,x}$ as column vectors for all $j$, we define 
	\begin{equation*}
		\Pi^i(x):=\lt(v^i_{1,x}, v^i_{2,x}, \dots, v^i_{q_0,x}\rt) \in M^{\frac{k(k+1)}{2}\times q_0},
	\end{equation*}
	and
	\begin{equation*}
		\Lambda^{i}(x):=\ared{\det\lt(\Pi^i(x)\lt( \Pi^i(x)\rt)^T\rt)}\text{ for any }x\in \R^{k(mn-k)}.
	\end{equation*}
	\cblue{Note that each entry of $\Pi^i(x)$ is an entry of $X_{j,x}^i$ for some $j\in\{1,\dots,q_0\}$. By previous discussions, we know that each entry of $\Pi^i(x)$ is a polynomial of $x$, and hence $\Lambda^{i}(x)$ is also a polynomial of $x$. Further, by Lemma \ref{grasaux1} (note the relation (\ref{graseq100})), we have that
		\begin{equation}
		\label{graseq12}
		\mathrm{dim}\lt(\mathrm{Span}\lt\{v^i_{1,x}, v^i_{2,x}, \dots, v^i_{q_0,x}\rt\}\rt)=\frac{k(k+1)}{2} \iff \Lambda^{i}(x)\not =0. 
		\end{equation}
		This together with (\ref{graseq101}) gives (\ref{graseq34}).}\nl
	
	\em Step 2. \rm \cblue{There exists $i_0\in\{1,2,\dots,p_0\}$ such that $\Lambda^{i_0}$ is non-trivial.}
	
	\em Proof of Step 2. \em We define a subspace 
	$V_0\in G(k,M^{m\times n})$ to be
	\begin{equation}
	\label{graseq5}
	V_0:=\lt\{\lt(\begin{array}{ccccccccccc} \xgreen{y_1} & 0 & 0 & 0 & \dots & 0 & 0 & 0 & \dots &  0 \\
	0 & \xgreen{y_1} & 0 & 0 & \dots & 0 & 0 & 0  &\dots & 0\\ 
	0 & 0 & \xgreen{y_2} & 0 & \dots & 0 & 0 & 0  & \dots &0 \\
	0 & 0 & 0 & \xgreen{y_2} & \dots & 0 & 0 & 0 & \dots &0\\
	\dots \\
	0 & 0 & 0 & 0 & \dots & \xgreen{y_k} & 0 & 0 & \dots &0\\
	0 & 0 & 0 & 0 & \dots & 0 & \xgreen{y_k} & 0 & \dots &0\\
	0 & 0 & 0 & 0 & \dots & 0 & 0 & 0 & \dots &0\\
	\dots \\
	0 & 0 & 0 & 0 & \dots & 0 & 0 & 0 & \dots &0\\
	\end{array}\rt): (\xgreen{y_1},\dots,\xgreen{y_k})\in \R^k\rt\}.
	\end{equation}
	\cblue{Note that the assumption $k\leq \frac{1}{2}\min\{m,n\}$ allows to construct the subspace $V_0$ in $M^{m\times n}$.} Now for some $i_0\in \lt\{1,2,\dots, p_0\rt\}$ and $x_0\in \R^{k(mn-k)}$ we have $\phi_{W_0^{i_0},W_1^{i_0}}(x_0)=V_0$. Thus by (\ref{graseq1}) we have 
	$V_0=\lt\{(v,T_{x_0}^{i_0}(v)):v\in W^{i_0}_0\rt\}$. Since $\mathcal{B}_{x_0}^{i_0}$ (recall (\ref{graseq7})) is a basis \cblue{of $V_0$}, 
	the mapping $H:\R^k\rightarrow V_0$ defined 
	by 
	\begin{equation}
	\label{graseq13}
	H(y):=\sum_{l=1}^k y_l \lt(a_l^{i_0}+T_{x_0}^{i_0}(a_l^{i_0})\rt)
	\end{equation}
	is a linear isomorphism onto $V_0$. Thus there exist $h_{s,t}\in \R^k$ for $s=1,2,\dots, m$, $t=1,2,\dots, n$ such that 
	\begin{equation*}
		H(y)=\lt(\begin{array}{cccc} h_{1,1}\cdot y & h_{1,2}\cdot y & \dots & h_{1,n}\cdot y\\
			h_{2,1}\cdot y & h_{2,2}\cdot y & \dots & h_{2,n}\cdot y\\
			\dots \\
			h_{m,1}\cdot y & h_{m,2}\cdot y & \dots & h_{m,n}\cdot y \end{array}\rt).
	\end{equation*}
	By definition of $V_0$ (recall (\ref{graseq5})) we have that 
	\begin{equation}
	\label{graseq15}
	h_{s,t}=0\text{ for }s\not=t, h_{2s-1,2s-1}=h_{2s,2s}\text{ for }s=1,2,\dots, k,\text{ and }h_{s,s}=0\text{ for }\cblue{s>2k}.
	\end{equation}
	
	Now we claim that $\lt\{h_{1,1}, h_{3,3}, \dots, h_{2k-1,2k-1}\rt\}$ are linearly independent. Suppose this is false, \xgreen{then 
		pick $y_1\in \bigcap_{s=1}^k \lt(  h_{2s-1,2s-1}\rt)^{\perp}\backslash \lt\{0\rt\}$. Thus $H(y_1)=0\in M^{m\times n}$, contradicting the fact that 
		$H$ is an isomorphism.} Thus the claim is established. 
	
	Note that $Q_{j,x_0}^{i_0}(y)\overset{(\ref{graseq9}), (\ref{graseq13})}{=}M_j(H(y))$ for $j=1,2,\dots, q_0$. \cblue{Using (\ref{graseq15}), the set $\OI_M$ (\pb{as subset of the polynomial ring $\R[y_1,\dots,y_k]$}) of all the $2\times 2$ minors on $V_0$ is simply}
	\begin{equation*}
		\begin{split}
			\OI_M&:=\lt\{Q_{j,x_0}^{i_0}(y):j=1,2,\dots, q_0\rt\}\\
			&= \lt\{ (h_{\ared{2s_1-1,2s_1-1}}\cdot y)  (h_{\ared{2s_2-1,2s_2-1}}\cdot y)  :\xgreen{s_1<s_2\in \lt\{1,2,\dots, k\rt\}}\rt\}\\
			&~\qd\qd \bigcup \lt\{ (h_{\ared{2t-1,2t-1}}\cdot y)^2 :t\in \lt\{1,2,\dots, k\rt\}\rt\}. 
		\end{split}
	\end{equation*}
	Now we claim that $\OI_M$ (\pb{as subset of the polynomial ring $\R[y_1,\dots,y_k]$}) is linearly independent. So see this, let $\lm_{s_1,s_2}\in \R$ and $\lm_{t}\in \R$ be such that 
	\begin{equation}
	\label{graseq17}
	0=\sum_{\xgreen{s_1<s_2\in \lt\{1,2,\dots, k\rt\}}}  
	\lm_{s_1,s_2} (h_{\ared{2s_1-1,2s_1-1}}\cdot y)  (h_{\ared{2s_2-1,2s_2-1}}\cdot y)+\sum_{t\in \lt\{1,2,\dots, k\rt\}} \lm_t  (h_{\ared{2t-1,2t-1}}\cdot y)^2.
	\end{equation}
	For every $t\in \lt\{1,2,\dots, k\rt\}$, pick $y\in \bigcap_{s\in \lt\{1,2,\dots, k\rt\}\backslash \lt\{t\rt\}} (h_{\ared{2s-1,2s-1}})^{\perp}$ such that $h_{\ared{2t-1,2t-1}}\cdot y\not=0$. \cblue{Note that such $y$ exists because $\lt\{h_{1,1}, h_{3,3}, \dots, h_{2k-1,2k-1}\rt\}$ are linearly independent}. Putting this into (\ref{graseq17}) we get that 
	$\lm_{t} (h_{\ared{2t-1,2t-1}}\cdot y)^2=0$ and so $\lm_t=0$. Thus $\lm_t=0$ for all $t\in\{1,2,\dots,k\}$. Next, let $\xgreen{s_1<s_2\in \lt\{1,2,\dots, k\rt\}}$. Pick 
	$$
	y\in \bigcap_{s\in \lt\{1,2,\dots, k\rt\}\backslash \lt\{s_1,s_2\rt\}} (h_{\ared{2s-1,2s-1}})^{\perp}
	$$ 
	such that 
	$y\cdot h_{\ared{2s_1-1,2s_1-1}}\not=0$ and $y\cdot h_{\xgreen{2s_2-1,2s_2-1}}\not=0$. \cblue{Such $y$ exists for the same reason as above.} Putting this into (\ref{graseq17}) we have that 
	$\lm_{s_1,s_2} (h_{\ared{2s_1-1,2s_1-1} }\cdot y)  (h_{\ared{2s_2-1,2s_2-1}}\cdot y)=0$ and so $\lm_{s_1,s_2} =0$. Thus \xgreen{linear} independence \cblue{of $\OI_M$} is 
	established. 
	
	It is easy to see that
	\begin{equation*}
		\ca{\OI_M}=\frac{k!}{2(k-2)!}+k=\frac{k(k+1)}{2}.
	\end{equation*}
	Thus 
	$$\mathrm{dim}\lt(\cblue{\mathrm{Span}}\lt\{X_{j,x_0}^{i_0}:j=1,2,\dots, q_0\rt\}\rt)=\mathrm{dim}\lt(\cblue{\mathrm{Span}}\lt\{Q_{j,x_0}^{i_0}:j=1,2,\dots, q_0\rt\}\rt)=\frac{k(k+1)}{2},$$ 
	and so by \cblue{(\ref{graseq101})} and (\ref{graseq12}) we have that $\Lambda^{i_0}(x_0)\not=0$. This completes the proof Step 2. \nl
	
	\em Step 3. \rm \cblue{We show that $\Lambda^i(x)$ is non-trivial on $R^{k(mn-k)}$ for all $i\in\{1,2,\dots,p_0\}$.}
	
	\em Proof of Step 3. \rm We first denote the zero set of $\Lambda^i(x)$ by
	\begin{equation*}
		Z^i:=\lt\{x\in \R^{k(mn-k)}:\Lambda^i(x)=0\rt\}.
	\end{equation*}
	\cblue{Note that, as $\Lambda^i(x)$ is a polynomial function, the set $Z^i$ is a real algebraic variety. A classical result of Whitney \cite{wh} states that $Z^i$ can be decomposed as a disjoint union of finitely many connected analytic submanifolds of dimension less than $k(mn-k)$, provided that $\Lambda^i$ is non-trivial on $\ared{\R}^{k(mn-k)}$.} 
	
	As the collection of charts $\lt\{\phi_{W^1_0,W^1_1}, \phi_{W^2_0,W^2_1}, \dots,  \phi_{W^{p_0}_0,W^{p_0}_1}\rt\}$ satisfy 
	(\ref{graseq39}) \xgreen{and $G(k,M^{m\times n})$ is connected}, it is clear that we can find $i_1\in \lt\{1,2,\dots, p_0\rt\}$ such that 
	\begin{equation}
	\label{graseq40}
	\UI_{i_0,i_1}:=\phi_{W^{i_0}_0, W^{i_0}_1}\lt(\R^{k(mn-k)}\rt)\bigcap \phi_{W^{i_1}_0, W^{i_1}_1}\lt(\R^{k(mn-k)}\rt)\not=\emptyset.
	\end{equation}
	We begin by showing that $\Lambda^{i_1}$ is non-trivial. \cblue{Since $\UI_{i_0, i_1}$ is a nonempty open subset of $G(k,M^{m\times n})$ and $\phi_{W^{i_0}_0, W^{i_0}_1}$ is a Lipschitz mapping, we know that $(\phi_{W^{i_0}_0, W^{i_0}_1})^{-1}(\UI_{i_0, i_1})$ is a nonempty open subset of $\ared{\R}^{k(mn-k)}$. As $\Lambda^{i_0}$ is non-trivial, we know from Whitney's result \cite{wh} that} $Z^{i_0}$ is the \cblue{disjoint} union of finitely many submanifolds of dimension less than $k(mn-k)$. It follows that
	\begin{equation*}
		\cblue{(\phi_{W^{i_0}_0, W^{i_0}_1})^{-1}(\UI_{i_0, i_1})\not \subset Z^{i_0}.}
	\end{equation*}
	Thus we must be able to find $\xgreen{\ti x_0}\in \R^{k(mn-k)}$ such that $\phi_{W^{i_0}_0, W^{i_0}_1}(\xgreen{\ti x_0})\in U_{i_0,i_1}$ and $\Lambda^{i_0}(\xgreen{\ti x_0})\not=0$. From (\ref{graseq40}), as $U_{i_0,i_1}\subset \phi_{W^{i_1}_0, W^{i_1}_1}\lt(\R^{k(mn-k)}\rt)$, we can find some
	$x_1\in \R^{k(mn-k)}$ such that $\phi_{W^{i_0}_0, W^{i_0}_1}(\xgreen{\ti x_0})= \phi_{W^{i_1}_0, W^{i_1}_1}(x_1)=:V_1$. Thus 
	\begin{equation*}
		\lt\{v+T_{\xgreen{\ti x_0}}^{i_0}(v):v\in W^{i_0}_0\rt\}=V_1=\lt\{v+T_{x_1}^{i_1}(v):v\in W^{i_1}_0\rt\}.
	\end{equation*}
	Since $\Lambda^{i_0}(\xgreen{\ti x_0})\not=0$, we know from Step 1 that $\mathrm{Span}\lt\{X^{i_0}_{1,\xgreen{\ti x_0}}, X^{i_0}_{2,\xgreen{\ti x_0}}, \dots, X^{i_0}_{q_0,\xgreen{\ti x_0}}\rt\}=M^{k\times k}_{sym}$. By Lemma \ref{graslem1.5}, we have that 
	$\mathrm{Span}\lt\{X^{i_1}_{1,x_1}, X^{i_1}_{2,x_1}, \dots, X^{i_1}_{q_0,x_1}\rt\}=M^{k\times k}_{sym}$. From Step 1 again, we have $\Lambda^{i_1}(x_1)\ne 0$ and hence $\Lambda^{i_1}$ is non-trivial. \cblue{From Lemma \ref{graslem1}, $G(k,M^{m\times n})$ is connected. Therefore, \xgreen{because of (\ref{graseq39})} the above arguments can be repeated to all charts to conclude that $\Lambda^i(x)$ is non-trivial for all $i\in\{1,2,\dots,p_0\}$.} \nl

	\em Proof of Lemma \ref{graslem2} completed. \rm Let $M^{k\times k}_{sym, +}$ denote the cone of positive definite matrices in $M^{k\times k}_{sym}$, i.e.,  
	$A\in M^{k\times k}_{sym, +}$ if and only if $y^T A y>0$ for all $y\in \R^k\cblue{\setminus\{0\}}$. If $V\in \phi_{W^{i}_0,W^{i}_1}(\R^{k(mn-k)}\backslash Z^{i})$, then there exists $x\in \R^{k(mn-k)}$ such that 
	$\phi_{W^{i}_0,W^{i}_1}(x)=V$ and $\Lambda^i(x)\ne 0$.  Now by (\ref{graseq34})  we must be able to 
	find some $\beta\in S^{q_0-1}$ such that $\sum_{j=1}^{q_0} \beta_j X_{j,x}^i\in  M^{k\times k}_{sym, +}$. By (\ref{graseq9}) this 
	implies that $\sum_{j=1}^{q_0} \beta_j M_j(X)>0$ for all $X\in V\backslash \lt\{0\rt\}$. Since $\Lambda^{i}$ is non-trivial \cblue{by Step 3}, we know from \cblue{Whitney's result \cite{wh}} that $Z^{i}$ is the \cblue{disjoint} union of finitely many submanifolds of dimension less than $k(mn-k)$. Thus (\ref{graseq3}) holds \em generically \rm (recall Definition \ref{grasdef1}) for 
	\cblue{$V\in  \phi_{W^{i}_0,W^{i}_1}(\R^{k(mn-k)})$ for all $i\in\{1,2,\dots,p_0\}$. We conclude the proof of the lemma by noting (\ref{graseq39}).} 
\end{proof}

\begin{proof}[Proof of Theorem \ref{C17}]
	\zzred{Recall that $M_1, M_2, \dots, M_{q_0}$ denote all $2\times 2$ minors in $M^{m\times n}$. 
		By Lemma \ref{graslem2}, we have that for generic subspace $V\in G(k,M^{m\times n})$ there exists 
		$\beta\in S^{q_0-1}$ such that (\ref{graseq3}) holds true. Now for any $\mu\in \MI^{pc}(V)$, using the 
		expansion (\ref{eqp201}) we know 
		\begin{equation}
		\label{graseq60}
		\int \sum_{j=1}^{q_0} \beta_j M_j(X-\ol{\cblue{X}}) d\mu=\sum_{j=1}^{q_0} \beta_j M_j(0)=0. 
		\end{equation}
		However by (\ref{graseq3}), unless $\mu=\delta_{\ol{X}}$  the left hand side of (\ref{graseq60}) is strictly positive, which is a contradiction. }
\end{proof}

\section{Preliminaries for Theorems \ref{t102}}
\label{prelims}

In this section, we gather some preliminary lemmas that will be useful in dealing with $\mathcal{M}^{pc}(\mathcal{K}_1)$ in the following section. First, we introduce some notation that will be used repeatedly. \pblue{Given a matrix $A\in M^{m\times n}$, recall the notation (\ref{eqp2}) and (\ref{eqp3})}. 
Let $A\in M^{m\times 2}$ with $m\geq 2$ and $1\leq i<j\leq m$, we define 
\begin{equation}
\label{noteq2}
X_{ij}(A):=\lt(\begin{array}{cc} \lt[A\rt]_{i1} &   \lt[A\rt]_{i2} \\ \lt[A\rt]_{j1} &   \lt[A\rt]_{j2}   \end{array}\rt)
\end{equation}
and
\begin{equation}
\label{noteq3}
M_{ij}(A):=\xblue{R_i(A)\wedge R_j(A)=}\det\lt(X_{ij}(A)\rt).
\end{equation}
\indent Recall the definitions of $\K_1$ and $P_1$ in (\ref{eq205}) and (\ref{auxxeq7.4}), respectively. Further, given $\alpha\in\R^2$, define
\begin{equation}
\label{bvv21}
P_1^{\alpha}(u,v):=\lt(\begin{array}{cc}   u-\alpha_1 & v-\alpha_2 \\
a(v)-a(\alpha_2) & u-\alpha_1\\
\lt(u-\alpha_1\rt)\lt(a(v)-a(\alpha_2)\rt) & \frac{\lt(u-\alpha_1\rt)^2}{2}+F(v)-F(\alpha_2)-a(\alpha_2)(v-\alpha_2) \end{array}\rt)
\end{equation}
and
\begin{equation}
\label{bvv22}
\mathcal{K}^{\alpha}_1:=\lt\{P_1^{\alpha}(u,v):u,v\in \R \rt\}.
\end{equation}
Finally, given a measure $\mu$ and a function $f$ which is integrable with respect to the measure $\mu$, define 
\begin{equation*}
	\overline{f}:=\int f(z) d\mu(z). 
\end{equation*}
We prove a couple of lemmas that will be essential in the following section.
%
%

\begin{a1} 
	\label{LBBB3}
	\bblue{For all} $\alpha\in\R^2$ the push forward mapping $(P_1^{\alpha})_{\sharp}:\mathcal{P}(\R^2)\rightarrow \mathcal{P}(\mathcal{K}_1^{\alpha})$ defined by $\nu\mapsto \mu:=(P_1^{\alpha})_{\sharp}\nu$ forms a bijection. Moreover, $\mu\in \mathcal{M}^{pc}(\mathcal{K}_1^{\alpha})$ if and only if 
	\begin{equation}
	\label{auxxeq3}
	\int_{\R^2} M_{ij}(P_1^{\alpha}(u,v)) d\nu =M_{ij}\lt( \int_{\R^2} P_1^{\alpha}(u,v) d\nu  \rt)\text{ for }i<j\in \orgg{\lt\{1,2,3\rt\}}.
	\end{equation}
	Further, given $\delta>0$, if $\spt \mu\subset \mathcal{K}_1^{\alpha}\cap B_{\delta}(0)$ then $\spt \nu\subset B_{\delta}(\alpha)$, and conversely, if $\spt \nu\subset B_{\delta}(\alpha)$ then $\spt \mu\subset \mathcal{K}_1^{\alpha}\cap B_{C\delta}(0)$ for some constant $C$ depending on the function $a$, $\alpha_2$ and $\delta$.
\end{a1}

\begin{proof} 
	\blue{First note that, since the first row of $P_1^{\alpha}(u,v)$ is $(u-\alpha_1,v-\alpha_2)$, it is clear that $P_1^{\alpha}:\R^2\rightarrow \mathcal{K}_1^{\alpha}$ is a bijection. Therefore it is straightforward to check that $\lt((P_1^{\alpha})^{-1}\rt)_{\sharp}$ is the inverse mapping of $(P_1^{\alpha})_{\sharp}$ and hence $(P_1^{\alpha})_{\sharp}$ is a bijection.}
	
	Let $\nu\in\mathcal{P}(\R^2)$ and $\mu\in\mathcal{P}(\mathcal{K}_1^{\alpha})$ be related by $\mu=(P_1^{\alpha})_{\sharp}\nu$. By change of variable formula for push forward measures, we have
	\begin{equation*}
		\int_{\R^2} M_{ij}(P_1^{\alpha}(u,v)) d\nu= \int_{\mathcal{K}_1^{\alpha}} M_{ij}(\zeta) d\mu(\zeta)
	\end{equation*}
	and
	\begin{equation*}
		M_{ij}\lt(\int_{\R^2} P_1^{\alpha}(u,v)) d\nu\rt)= M_{ij}\lt(\int_{\mathcal{K}_1^{\alpha}} \zeta d\mu(\zeta)\rt).
	\end{equation*}
	It follows that $\mu\in \mathcal{M}^{pc}(\mathcal{K}_1^{\alpha})$ if and only if (\ref{auxxeq3}) holds.
	
	\blue{Next, assume $\spt \mu\subset \mathcal{K}_1^{\alpha}\cap B_{\delta}(0)$. Since the first row of $P_1^{\alpha}(u,v)$ is $(u-\alpha_1,v-\alpha_2)$, it is clear that $\lVert(u,v)-\alpha \rVert \leq \lVert P_1^{\alpha}(u,v) \rVert$. Therefore $(P_1^{\alpha})^{-1}(\mathcal{K}_1^{\alpha}\cap B_{\delta}(0)) \subset B_{\delta}(\alpha)$. \ggreen{As $\spt\nu=\lt(P^{\alpha}_1\rt)^{-1} \spt\mu$}, \xblue{it} follows that $\spt \nu\subset B_{\delta}(\alpha)$. Conversely, assume $\spt \nu\subset B_{\delta}(\alpha)$. From the expression of $P_1^{\alpha}(u,v)$ in (\ref{bvv21})  it is clear that the absolute value of each component of $P_1^{\alpha}(u,v)$ is bounded above by  $C\lVert(u,v)-\alpha \rVert$ for some constant $C$ depending on the function $a$, $\alpha_2$ and $\delta$, provided $\delta$ is sufficiently small. Therefore $\lVert P_1^{\alpha}(u,v) \rVert \leq \ti C\lVert(u,v)-\alpha \rVert$ and hence $P_1^{\alpha}(B_{\delta}(\alpha)) \subset \mathcal{K}_1^{\alpha}\cap B_{\ti C\delta}(0)$. It follows that $\spt \mu\subset \mathcal{K}_1^{\alpha}\cap B_{\ti C\delta}(0)$. This completes the proof of the lemma.}
\end{proof}

%
%
\xred{The following lemma is implicitly stated in \cite{kms}. We thank S. M\"{u}ller  \cite{steph1} for providing us with the elegant proof presented in this section. 
	\begin{a1}[Kirchheim-M\"{u}ller-\v{S}ver\'{a}k \cite{kms}]
		\label{LAUX33}
		Given $\nu\in\mathcal{P}(\R^2)$, for all $\alpha \in \R^2$ we have
		\begin{equation*}
		(P_1)_{\sharp}\nu \in \mathcal{M}^{pc}(\mathcal{K}_1)\Longleftrightarrow  (P_1^{\alpha})_{\sharp}\nu\in   \mathcal{M}^{pc}(\mathcal{K}^{\alpha}_1).
		\end{equation*}
\end{a1}}

We break the proof into several steps. The first lemma is standard. 

\begin{a1}
	\label{LAUX34}
	Given $\nu\in\mathcal{P}(\R^2)$, \bblue{for all} $\alpha \in \R^2$, we have
	\begin{equation}
	\label{eqyytt4}
	(P_1^{\alpha})_{\sharp}\nu \in \mathcal{M}^{pc}(\mathcal{K}_1^{\alpha})\Longleftrightarrow (\ti{P_1^{\alpha}})_{\sharp}\nu\in     \mathcal{M}^{pc}(\tilde{\mathcal{K}}^{\alpha}_1),
	\end{equation}
	where 
	\begin{equation}
	\label{eqyytt301}
	\ti{P}_1^{\alpha}(u,v):=\lt(\begin{array}{cc}   u & v \\
	a(v) & u\\
	u a(v)-\alpha_1 a(v)-u a(\alpha_2) & \frac{u^2}{2}-u \alpha_1+F(v)-v a(\alpha_2)\end{array}\rt)
	\end{equation}
	and
	\begin{equation*}
		\tilde{\mathcal{K}}^{\alpha}_1:=\lt\{ \ti{P}_1^{\alpha}(u,v)   :u,v\in \R\rt\}.
	\end{equation*}
\end{a1}	

\begin{proof}
	Recall the definition of $P_1^{\alpha}$ given in (\ref{bvv21}). Direct calculations show that
	\begin{eqnarray}
	\label{eqyytt3004}
	P_1^{\alpha}(u,v)=\ti{P}_1^{\alpha}(u,v)-\ti{E}^{\alpha}
	\end{eqnarray}	
	for
	\begin{equation*}
		\ti{E}^{\alpha}:=\lt(\begin{array}{cc}  \alpha_1 & \alpha_2 \\
			a(\alpha_2) & \alpha_1 \\
			-\alpha_1 a(\alpha_2) & \bblue{-\frac{\ap_1^2}{2}}+F(\alpha_2)-\alpha_2 a(\alpha_2)\end{array}\rt).
	\end{equation*}
	Note that $\ti{E}^{\alpha}$ is the constant part of $P_1^{\alpha}(u,v)$.

	Given $\nu\in \mathcal{P}(\R^2)$, by arguing exactly as in Lemma \ref{LBBB3}, we have that 
	\begin{equation}
	\label{eqyytt1.8}
	(\ti{P}_1^{\alpha})_{\sharp}\nu \in \mathcal{M}^{pc}(\tilde{\mathcal{K}}^{\alpha}_1)\Longleftrightarrow \int_{\R^2} M_{ij}(\ti{P}_1^{\alpha}(u,v)) d\nu =M_{ij}\lt( \int_{\R^2} \ti{P}_1^{\alpha}(u,v) d\nu  \rt)\text{ for } \orgg{i<j\in \lt\{1,2,3\rt\}}.
	\end{equation}
	Using \bblue{the formula (\ref{eqp201})}, the fact that $\ti{E}^{\alpha}$ is a constant matrix and the notation $M_{ij}(\cdot)=\det\lt(X_{ij}(\cdot)\rt)$ \rred{(recalling (\ref{noteq2}))}, we have
	\begin{equation}
	\label{eqyytt303}
	\begin{split}
	&\int M_{ij}\lt(P_1^{\alpha}(u,v)\rt) d\nu \overset{(\ref{eqyytt3004})}{=}\int M_{ij}\lt(\ti{P}_1^{\alpha}(u,v)-\ti{E}^{\alpha} \rt) d\nu \\
	&  = \int \lt[   M_{ij}\lt( \ti{P}_1^{\alpha}(u,v) \rt)-X_{ij}\lt( \ti{P}_1^{\alpha}(u,v) \rt):\mathrm{Cof}\lt( X_{ij}\lt(\ti{E}^{\alpha} \rt) \rt)
	+M_{ij}\lt(\ti{E}^{\alpha} \rt) \rt] d\nu\\
	&  = \int M_{ij}\lt( \ti{P}_1^{\alpha}(u,v) \rt)d\nu- \int X_{ij}\lt( \ti{P}_1^{\alpha}(u,v) \rt)d\nu:\mathrm{Cof}\lt( X_{ij}\lt(\ti{E}^{\alpha} \rt) \rt)
	+M_{ij}\lt(\ti{E}^{\alpha} \rt).
	\end{split}
	\end{equation}
	In a similar way \bblue{using (\ref{eqp201})} we have that 
	\begin{equation}
	\label{eqyytt304}
	\begin{split}
	&\det\lt(\int X_{ij}\lt(P_1^{\alpha}(u,v) \rt) d\nu\rt) \\
	&=\ggreen{\det\lt(\int X_{ij}\lt(\ti{P}_1^{\alpha}(u,v)-\ti{E}^{\alpha}\rt) d\nu\rt)}\\
	&= \det\lt(\int X_{ij}\lt(\ti{P}_1^{\alpha}(u,v)\rt)d\nu\rt)-\int X_{ij}\lt(\ti{P}_1^{\alpha}(u,v)\rt)d\nu: \mathrm{Cof}\lt(X_{ij}\lt(\ti{E}^{\alpha}\rt)\rt)+M_{ij}\lt(\ti{E}^{\alpha} \rt).
	\end{split}
	\end{equation}
	Putting (\ref{eqyytt303}) and (\ref{eqyytt304}) together \ggreen{and using Lemma \ref{LBBB3}} we have that 
	\begin{eqnarray}
	\label{eqyytt306}
	&~&(P_1^{\alpha})_{\sharp}\nu \in \mathcal{M}^{pc}(\mathcal{K}_1^{\alpha})\nn\\
	&~&\qd\qd \overset{\ggreen{(\ref{auxxeq3})}}{\Longleftrightarrow} \int M_{ij}\lt(P_1^{\alpha}(u,v)\rt) d\nu=\det\lt(\int X_{ij}\lt(P_1^{\alpha}(u,v) \rt) d\nu\rt) \text{ for all } \orgg{i<j\in \lt\{1,2,3\rt\}}\nn\\
	&~&\qd\qd\overset{(\ref{eqyytt303}),(\ref{eqyytt304})}{\Longleftrightarrow } \int M_{ij}\lt( \ti{P}_1^{\alpha}(u,v) \rt)d\nu = \det\lt(\int X_{ij}\lt(\ti{P}_1^{\alpha}(u,v)\rt)d\nu\rt) \text{ for all } \orgg{i<j\in \lt\{1,2,3\rt\}}\nn\\
	&~&\qd\qd\overset{\ggreen{(\ref{eqyytt1.8})}}{\Longleftrightarrow} (\ti{P}_1^{\alpha})_{\sharp}\nu\in     \mathcal{M}^{pc}(\tilde{\mathcal{K}}^{\alpha}_1)\nn.
	\end{eqnarray}
	This establishes (\ref{eqyytt4}).
\end{proof}

\begin{a1}[\xred{M\"{u}ller} \xblue{\cite{steph1}}]
	\label{le:minors}
	Every row $R_i\lt(\ti{P}_1^{\alpha}(u,v)\rt)$ of the matrix $\ti{P}_1^{\alpha}(u,v)$ can be expressed as a linear combination of the 
	rows of $P_1(u,v)$, and conversely every row $R_i\lt(P_1(u,v)\rt)$ of the matrix $P_1(u,v)$ can be expressed as a linear combination of the rows of $\ti{P}_1^{\alpha}(u,v)$, and the coefficients depend only on $\alpha$, but not on $(u,v)$, i.e.,
	\begin{equation}\label{eqpp103}
	R_i\lt(\ti{P}_1^{\alpha}(u,v)\rt) = \sum_{i'=1}^{\orgg{3}}  c_{ii'}(\alpha) \, R_{i'}\lt(P_1(u,v)\rt)   \quad \text{ for all } (u,v)\in \R^2
	\end{equation}
	and 
	\begin{equation}\label{eqpp104}
	R_i\lt(P_1(u,v)\rt) = \sum_{i'=1}^{\orgg{3}}  \ti{c}_{ii'}(\alpha) \, R_{i'}\lt(\ti{P}_1^{\alpha}(u,v)\rt)   \quad \text{ for all } (u,v)\in \R^2.
	\end{equation}
\end{a1}

\begin{proof} From the definitions of $P_1$ and $\ti{P}_1^{\alpha}$ in (\ref{auxxeq7.4}) and (\ref{eqyytt301}), we see that
	\begin{equation}\label{eqpp107}
	R_1\lt(\ti{P}_1^{\alpha}(u,v)\rt)=R_1\lt(P_1(u,v)\rt), \quad R_2\lt(\ti{P}_1^{\alpha}(u,v)\rt)=R_2\lt(P_1(u,v)\rt).
	\end{equation}
	Now we calculate
	\begin{equation}\label{eqpp105} 
	\begin{split}
	R_3\lt(\ti{P}_1^{\alpha}(u,v)\rt)
	&=  \lt(u a(v)-\alpha_1 a(v)-u a(\alpha_2), \frac{u^2}{2}-u \alpha_1+F(v)-v a(\alpha_2)\rt) \\
	&= R_3\lt(P_1(u,v)\rt) - \alpha_1  R_2\lt(P_1(u,v)\rt) - a(\ap_2) R_1\lt(P_1(u,v)\rt).
	\end{split}
	\end{equation}
	This proves (\ref{eqpp103}). Conversely, we have
	\begin{equation*}
		R_3\lt(P_1(u,v)\rt) \overset{(\ref{eqpp105}),(\ref{eqpp107})}{=} R_3\lt(\ti{P}_1^{\alpha}(u,v)\rt) + \alpha_1  R_2\lt(\ti{P}_1^{\alpha}(u,v)\rt) + a(\ap_2) R_1\lt(\ti{P}_1^{\alpha}(u,v)\rt)
	\end{equation*}
	and therefore we have (\ref{eqpp104}).
\end{proof}

\begin{proof}[Proof of Lemma \ref{LAUX33}]
	By Lemma \ref{LAUX34}, it suffices to show that
	\begin{equation}\label{eqpp102}
	(P_1)_{\sharp}\nu \in \mathcal{M}^{pc}(\mathcal{K}_1)\Longleftrightarrow      (\ti{P}_1^{\alpha})_{\sharp}\nu\in     \mathcal{M}^{pc}(\tilde{\mathcal{K}}^{\alpha}_1).
	\end{equation}
	We first show the implication ``$\Longrightarrow$''. 
	We denote $\mu:=(P_1)_{\sharp}\nu$ and $\mu^{\ap}:=(\ti{P}_1^{\alpha})_{\sharp}\nu$, and assume $\mu\in \mathcal{M}^{pc}(\mathcal{K}_1)$. Recall the notation $M_{ij}$ given by (\ref{noteq3}), and denote $\ol \zeta := \int_{\K_1} \zeta \, d \mu(\zeta)$.
	Then by the change of variable formula for push forward measures, Lemma \ref{le:minors} and bilinearity of the minor we have
	\begin{equation}\label{eqpp101}
	\begin{split}
	\int_{\tilde{\mathcal{K}}^{\alpha}_1} M_{ij}(\zeta) \, d\mu^\alpha &= \int_{\R^2} M_{ij}\lt(\ti{P}_1^{\alpha}(u,v)\rt) \, d\nu \\
	&=  \int_{\R^2}  \sum_{i', j' =1}^{\orgg{3}}  c_{ii'} c_{jj'}   \,  M_{i'j'}\lt(P_1(u,v)\rt) \, d\nu \\
	&=  \sum_{i', j' =1}^{\orgg{3}}   c_{ii'} c_{jj'}   \int_{\K_1}    M_{i'j'}\lt(\zeta\rt) \, d\mu \\
	&\overset{\mu \in \mathcal{M}^{pc}(\K_1)}{=}  \sum_{i', j' =1}^{\orgg{3}}   c_{ii'} c_{jj'}     M_{i'j'}(\ol \zeta).
	\end{split}
	\end{equation}
	On the other hand bilinearity of the minor implies that
	\begin{equation*}
		\begin{split}
			&M_{ij}\lt(\int_{\tilde{\mathcal{K}}^{\alpha}_1}\zeta\, d\mu^\alpha \rt) =  M_{ij} \left(   \int_{\R^2}  \ti{P}_1^{\alpha}(u,v)   \, d\nu \right) \\
			&\qd\qd\qd\qd\overset{(\ref{noteq3})}{=}   \int_{\R^2}R_i\lt(\ti{P}_1^{\alpha}(u,v)\rt)\,d\nu \wedge \int_{\R^2}R_j\lt(\ti{P}_1^{\alpha}(u,v)\rt)\,d\nu\\
			&\qd\qd\qd\qd= \sum_{i',j'=1}^{\orgg{3}} c_{ii'} c_{jj'}   \int_{\R^2}R_{i'}\lt(P_1(u,v)\rt)\,d\nu \wedge \int_{\R^2}R_{j'}\lt(P_1(u,v)\rt)\,d\nu\\
			&\qd\qd\qd\qd=  \sum_{i',j'=1}^{\orgg{3}} c_{ii'} c_{jj'}  M_{i'j'} \left( \int_{\R^2} P_1(u,v) \, d\nu \right) 
			=  \sum_{i',j'=1}^{\orgg{3}} c_{ii'} c_{jj'}  M_{i'j'} \left( \int_{\K_1} \zeta  \, d\mu \right) \\
			&\qd\qd\qd\qd=  \sum_{i',j'=1}^{\orgg{3}} c_{ii'} c_{jj'}  M_{i'j'}(\overline \zeta)\overset{(\ref{eqpp101})}{=}  	\int_{\tilde{\mathcal{K}}^{\alpha}_1} M_{ij}(\zeta) \, d\mu^\alpha
		\end{split}
	\end{equation*}
	as desired. The proof of the converse implication is analogous using (\ref{eqpp104}). This completes the proof of (\ref{eqpp102}), and hence Lemma \ref{LAUX33}.
\end{proof}

\section{Existence of non-trivial measure in $\mathcal{M}^{pc}(\mathcal{K}_1)$}\label{s7}
In this section, we first construct non-trivial measures in $\mathcal{M}^{pc}(\mathcal{K}^{\bblue{\ti\alpha}}_1)$ in the case $a'(\ti\alpha_2)>0$. Then it follows from Lemma \ref{LAUX33} that we also have non-trivial elements in $\mathcal{M}^{pc}(\mathcal{K}_1)$. More precisely, we construct non-trivial measures supported at five points that belong to the space $\mathcal{M}^{pc}(\mathcal{K}^{\ti\alpha}_1)$. To begin with, given $s_0,t_0>0$, \rred{recalling (\ref{bvv21})}, we set
\begin{equation*}
	\zeta_0:\blue{=P_1^{\ti\alpha}(\ti\alpha_1,\ti\alpha_2)}=\lt(\begin{matrix} 0 & 0 \\ 0 & 0 \\ 0 & 0 \end{matrix}\rt),\; \zeta_1:\blue{=P_1^{\ti\alpha}(\ti\alpha_1+s_0,\ti\alpha_2)}=\lt(\begin{matrix} s_0 & 0 \\ 0 & s_0 \\ 0 & \frac{1}{2}s_0^2 \end{matrix}\rt),
\end{equation*}
\begin{equation*} 
	\qd \zeta_2:=P_1^{\ti\alpha}(\ti\alpha_1-s_0,\ti\alpha_2)= \lt(\begin{matrix} -s_0 & 0 \\ 0 & -s_0 \\ 0 & \frac{1}{2}s_0^2 \end{matrix}\rt),
\end{equation*}
\begin{equation*}
	\zeta_3:=P_1^{\ti\alpha}(\ti\alpha_1,\ti\alpha_2+t_0)= \lt(\begin{matrix} 0 & t_0 \\ a(\ti\alpha_2+t_0)-a(\ti\alpha_2) & 0 \\ 0 & F(\ti\alpha_2+t_0)-F(\ti\alpha_2)-a(\ti\alpha_2)t_0 \end{matrix}\rt),
\end{equation*}
and
\begin{equation*}
	\zeta_4:=P_1^{\ti\alpha}(\ti\alpha_1,\ti\alpha_2-t_0)= \lt(\begin{matrix} 0 & -t_0 \\ a(\ti\alpha_2-t_0)-a(\ti\alpha_2) & 0 \\ 0 &  F(\ti\alpha_2-t_0)-F(\ti\alpha_2)+a(\ti\alpha_2)t_0  \end{matrix}\rt).
\end{equation*}
\blue{We first prove
	\begin{a2}\label{t103}
		\xgreen{Suppose $a\in C^2(\R)$}. Let $\bblue{\ti\alpha}\in\R^2$ be such that $a'(\ti\alpha_2)>0$. Given $s_0,t_0>0$ sufficiently small depending on the function $a$ and $\ti\alpha_2$, there exists $0<\ep_0<1$ depending on the function $a$, $\ti\alpha_2$, $s_0$ and $t_0$ such that, for all $\ep\leq \ep_0$, there exists a collection of weights $\{\xgreen{[\gamma^{\ep}]_j}\}_{j=0}^4 \subset \R_{+}$ with $\sum_{j=1}^{4}\xgreen{[\gamma^{\ep}]_j}=\ep$ and $\xgreen{[\gamma^{\ep}]_0}=1-\ep$ such that
		\begin{equation*}
			\mu^{\ep}:=\sum_{j=0}^{4} \xgreen{[\gamma^{\ep}]_j}\delta_{\zeta_j} \in \mathcal{M}^{pc}(\mathcal{K}^{\ti\alpha}_1).
		\end{equation*}
\end{a2}}

The proof of Theorem \ref{t103} will rely on a couple of crucial lemmas. Let us first introduce some notations. We denote by $D_1, D_2, D_3$ the $(1,2), (2,3), (1,3)$ minors of a $3\times 2$ matrix, respectively. We set the matrix
\begin{equation}\label{eq2102}
A:= \lt(\begin{matrix} D_1(\zeta_1) & D_1(\zeta_2) & D_1(\zeta_3) & D_1(\zeta_4) \\ D_2(\zeta_1) & D_2(\zeta_2) & D_2(\zeta_3) & D_2(\zeta_4) \\ D_3(\zeta_1) & D_3(\zeta_2) & D_3(\zeta_3) & D_3(\zeta_4) \\ 1 & 1 & 1 & 1 \end{matrix}\rt).
\end{equation}
For any $\ep>0$ and $\gamma\in \R^4$, define
\begin{equation}
\label{eqdfg22}
L^{\ep}(\gamma) := A\gamma - \lt(\begin{matrix} 0 \\ 0 \\ 0 \\ \ep \end{matrix}\rt), \qd Q(\gamma) := \lt(\begin{matrix} D_1\lt(\sum_{j=1}^{4}\xgreen{[\pb{\gamma}]_j}\zeta_j\rt) \\ D_2\lt(\sum_{j=1}^{4}\xgreen{[\gamma]_j}\zeta_j\rt) \\ D_3\lt(\sum_{j=1}^{4}\xgreen{[\gamma]_j}\zeta_j\rt) \\ 0 \end{matrix}\rt),
\end{equation}
and
\begin{equation}
\label{eqdfg21}
G^{\ep}(\gamma) := L^{\ep}(\gamma) - Q(\gamma).
\end{equation}

\begin{a1}\label{l105}
	\xgreen{Suppose $a\in C^2(\R)$}. \bblue{Let $\bblue{\ti\alpha}\in\R^2$ be such that $a'(\ti\alpha_2)>0$. Given $s_0,t_0>0$ sufficiently small depending on the function $a$ and $\ti\alpha_2$}, the matrix $A$ defined in (\ref{eq2102}) is invertible. Moreover, for any $0<\ep<1$, the unique solution $\gamma_0^{\ep}$ of the system 
	\begin{equation}\label{eq2104}
	L^{\ep}(\gamma) = 0
	\end{equation}
	is non-negative componentwise. Further, there exist constants $0<\lambda<\Lambda<\infty$ depending on the function $a$, $\ti\alpha_2$, $s_0$ and $t_0$ such that 
	\begin{equation}\label{eq2103}
	\lambda\ep \leq\bblue{[\gamma_0^{\ep}]_i} \leq \Lambda\ep \ \text{ for }\  i=1,2,3,4.
	\end{equation}
\end{a1}

\begin{proof}
	To simplify notation define $a_{\ti\alpha_2}(t):=a(\ti\alpha_2+t)-a(\ti\alpha_2)$ \blue{ and $F_{\ti\alpha_2}(t):=F(\ti\alpha_2+t)-F(\ti\alpha_2)-a(\ti\alpha_2)t$}. First, explicit calculations using the formulas for $\zeta_j,\; j=1,2,3,4$, give
	\begin{equation}
	\label{eq2101.5}
	A=\lt(\begin{matrix} s_0^2 & s_0^2 & -t_0 a_{\ti\alpha_2}(t_0)  & t_0  a_{\ti\alpha_2}(-t_0) \\ 0 & 0 &   
	a_{\ti\alpha_2}(t_0)F_{\ti\alpha_2}(t_0)  &   a_{\ti\alpha_2}(-t_0)F_{\ti\alpha_2}(-t_0)     \\ \frac{1}{2}s_0^3 & -\frac{1}{2}s_0^3 & 0 & 0 \\ 1 & 1 & 1 & 1 \end{matrix}\rt).
	\end{equation}
	We claim that for any $(y_1,y_2,y_3)\ne 0\in\R^3$, we have
	\begin{equation}\label{eq2101}
	\min_{j}\lt\{\sum_{i=1}^{3} y_iD_i(\zeta_j)\rt\} < 0
	\end{equation}
	and 
	\begin{equation}\label{eq2106}
	\max_{j}\lt\{\sum_{i=1}^{3} y_iD_i(\zeta_j)\rt\} > 0.
	\end{equation}
	We check (\ref{eq2101}) by an enumerative argument. Note that \red{since $a'(\ti\alpha_2)>0$, assuming $t_0>0$ is small enough, we have that  
		\begin{equation}
		\label{eq2101.3}
		a_{\ti\alpha_2}(t_0)>0\text{ and }a_{\ti\alpha_2}(-t_0)<0. 
		\end{equation}
		\xgreen{Recall that $F'=a$.} We also know that $F$ is \xgreen{strictly} convex in small neighborhood of $\ti\alpha_2$ and so 
		\begin{equation}
		\label{eq2101.6}
		\rred{F_{\ti{\alpha}_2}(t_0)=}F(\ti\alpha_2+t_0)-F(\ti\alpha_2)-a(\ti\alpha_2)t_0>0\text{ and }\rred{F_{\ti{\alpha}_2}(-t_0)=}F(\ti\alpha_2-t_0)-F(\ti\alpha_2)+a(\ti\alpha_2)t_0>0.
		\end{equation}
		By carefully checking out the columns of $A$ and using (\ref{eq2101.3}), (\ref{eq2101.6}) we see that }
	\blue{\begin{enumerate}
			\item If $y_1>0, y_2\geq0, y_3\in\R$, we have $\sum_{i=1}^3 y_iD_i(\zeta_4)<0$.
			\item If $y_1>0, y_2\leq0, y_3\in\R$, we have $\sum_{i=1}^3 y_iD_i(\zeta_3)<0$.
			\item If $y_1<0, y_3\geq0, y_2\in\R$, we have $\sum_{i=1}^3 y_iD_i(\zeta_2)<0$.
			\item If $y_1<0, y_3\leq0, y_2\in\R$, we have $\sum_{i=1}^3 y_iD_i(\zeta_1)<0$.
			\item If $y_1=0$, and
			\begin{enumerate}
				\item $y_2>0, y_3\in\R$, we have $\sum_{i=1}^3 y_iD_i(\zeta_4)=y_2D_2(\zeta_4)<0$;
				\item $y_2<0, y_3\in\R$, we have $\sum_{i=1}^3 y_iD_i(\zeta_3)=y_2D_2(\zeta_3)<0$;
				\item $y_2=0$, $y_3>0$, we have $\sum_{i=1}^3 y_iD_i(\zeta_2)=y_3D_3(\zeta_2)<0$;
				\item $y_2=0$, $y_3<0$, we have $\sum_{i=1}^3 y_iD_i(\zeta_1)=y_3D_3(\zeta_1)<0$.
			\end{enumerate}
		\end{enumerate}
		The above (1), (2) cover all cases when $y_1>0$, and (3), (4) cover all cases when $y_1<0$. For the case $y_1=0$, we have either $y_2\ne 0$ or $y_2=0$. The former is covered by (5a), (5b), and the latter is covered by (5c), (5d). Therefore the above enumerative argument shows (\ref{eq2101})}, and (\ref{eq2106}) is equivalent to (\ref{eq2101}).  \pgr{Let $a_1,\dots,a_4$ denote the columns of the matrix $A$}. \pgr{Hence} given any $y\in \R^4\pgr{\setminus\{0\}}$ such that $y\cdot a_i\geq 0$ for \pgr{all} $i=1,2,3,4$, we must have $y_4>0$ \pgr{and thus $y\cdot\lt(0,0,0,\ep\rt)^{T}>0$}. \xblue{We deduce from the Farkas-Minkowski Lemma (\pgr{Lemma \ref{fm}}) that (\ref{eq2104})} has a  non-negative solution.
	
	To see the matrix $A$ is invertible, consider the following system
	\begin{equation}\label{eq2105}
	A^{T}y=0.
	\end{equation}
	We claim that the system (\ref{eq2105}) has only the trivial solution. Indeed, let $y=(y_1,y_2,y_3,y_4)\in \R^4$ be a solution of (\ref{eq2105}), i.e.,
	\begin{equation*}
		\sum_{i=1}^3 y_iD_i(\zeta_j) + y_4 = 0 \text{ for } j=1,2,3,4.
	\end{equation*}
	It follows from (\ref{eq2101}) and (\ref{eq2106}) that $y_4=0$, and therefore \red{by (\ref{eq2101}) and (\ref{eq2106}) again}, we have $y_i=0$ for $i=1,2,3$. It follows that $A^{T}$, and hence $A$, are invertible.
	
	Finally, we show (\ref{eq2103}). Let $\gamma_0^{\ep}$ be the unique solution of (\ref{eq2104}). We already know that $\gamma_0^{\ep}$ is non-negative componentwise. \rred{We will show} that all components of $\gamma_0^{\ep}$ are strictly positive. We argue by contradiction. Suppose $[\gamma_0^{\ep}]_1=0$ or $[\gamma_0^{\ep}]_2=0$. Then using the third row of (\ref{eq2104}), we have $\frac{s_0^3}{2}([\gamma_0^{\ep}]_1-[\gamma_0^{\ep}]_2)\overset{(\ref{eq2101.5})}{=}0$ and therefore we have $[\gamma_0^{\ep}]_1=[\gamma_0^{\ep}]_2=0$. Now using the first two rows of (\ref{eq2104}) 
	\rred{(see (\ref{eq2101.5}))}, we have
	\begin{equation*}
		\lt(\begin{matrix} -t_0 a_{\ti\alpha_2}(t_0)  & t_0  a_{\ti\alpha_2}(-t_0)\\ a_{\ti\alpha_2}(t_0)F_{\ti\alpha_2}(t_0)  &   a_{\ti\alpha_2}(-t_0)F_{\ti\alpha_2}(-t_0)\end{matrix}\rt) \lt(\begin{matrix} [\gamma_0^{\ep}]_3\\ [\gamma_0^{\ep}]_4 \end{matrix}\rt) = \lt(\begin{matrix} 0\\0 \end{matrix}\rt).
	\end{equation*}
	It is clear from (\ref{eq2101.3}), (\ref{eq2101.6}) that the matrix $\lt(\begin{matrix} -t_0 a_{\ti\alpha_2}(t_0)  & t_0  a_{\ti\alpha_2}(-t_0)\\ a_{\ti\alpha_2}(t_0)F_{\ti\alpha_2}(t_0)  &   a_{\ti\alpha_2}(-t_0)F_{\ti\alpha_2}(-t_0)\end{matrix}\rt)$ is invertible, and hence $[\gamma_0^{\ep}]_3=[\gamma_0^{\ep}]_4=0$. But now we have $\gamma_0^{\ep}=0$, which contradicts the fourth row of (\ref{eq2104}). This contradiction implies $[\gamma_0^{\ep}]_1>0$ and $[\gamma_0^{\ep}]_2>0$. A similar argument yields $[\gamma_0^{\ep}]_3>0$ and $[\gamma_0^{\ep}]_4>0$. Since $[\gamma_0^{\ep}]_i = \ep[A^{-1}]_{i4}$, $i=1,2,3,4$, it follows that $[A^{-1}]_{i4}>0$ for all $i$. Now we define 
	\begin{equation*}
		\lambda:=\min_{i}\{[A^{-1}]_{i4}\} \qd\text{ and }\qd \Lambda:=\max_{i}\{[A^{-1}]_{i4}\}.
	\end{equation*}
	It is clear that $0<\lambda\leq\Lambda<\infty$ and (\ref{eq2103}) is satisfied. Note that $A^{-1}$ is a fixed matrix independent of $\ep$, and so are $\lambda$ and $\Lambda$ independent of $\ep$.
\end{proof} 

\begin{a1}\label{l106}
	\xgreen{Suppose $a\in C^2(\R)$}. \bblue{Let $\bblue{\ti\alpha}\in\R^2$ be such that $a'(\ti\alpha_2)>0$. Given $s_0,t_0>0$ sufficiently small depending on the function $a$ and $\ti\alpha_2$}, there exists $0<\ep_0<1$ sufficiently small such that for all $0<\ep\leq \ep_0$, the system
	\begin{equation}\label{eq2107}
	G^{\ep}(\gamma)=0
	\end{equation}
	has a non-negative solution.
\end{a1}

\begin{proof}
	Given a sufficiently small $0<\ep<1$ whose size will be specified later, by Lemma \ref{l105}, the linear system $L^{\ep}(\gamma)=0$ has a unique non-negative solution $\gamma^{\ep}_0$ that satisfies the estimate (\ref{eq2103}). We will find a solution to (\ref{eq2107}) by iteration. For all $k\in\Nb^{+}$, define
	\begin{equation}\label{eq2116}
	\Delta^{\ep}_{k} := A^{-1}\lt(-G^{\ep}(\gamma^{\ep}_{k-1})\rt) \text{ and } \gamma^{\ep}_k := \gamma^{\ep}_{k-1}+\Delta^{\ep}_{k}.
	\end{equation}
	Then we have
	\begin{equation}\label{eq2115}
	\begin{split}
	G^{\ep}(\gamma^{\ep}_k) &\overset{\red{(\ref{eqdfg21})}}{=} L^{\ep}(\gamma^{\ep}_k) - Q(\gamma^{\ep}_k)\\
	&\overset{\ggreen{(\ref{eqdfg22})},\red{(\ref{eq2116})}}{=}A\lt(\gamma^{\ep}_{k-1}+\Delta^{\ep}_{k}\rt) - \lt(0,0,0,\ep\rt)^{T} - Q(\gamma^{\ep}_k)\\
	&\overset{(\ref{eqdfg22})}{=} L^{\ep}(\gamma^{\ep}_{k-1})+A(\Delta^{\ep}_{k})-Q(\gamma^{\ep}_k)\\
	&\overset{\red{(\ref{eqdfg21})}}{=} G^{\ep}\lt(\gamma^{\ep}_{k-1}\rt) + A\lt(\Delta^{\ep}_{k}\rt) + Q(\gamma^{\ep}_{k-1}) - Q(\gamma^{\ep}_k)\\
	&\overset{\red{(\ref{eq2116})}}{=}G^{\ep}\lt(\gamma^{\ep}_{k-1}\rt) -G^{\ep}\lt(\gamma^{\ep}_{k-1}\rt) + Q(\gamma^{\ep}_{k-1}) - Q(\gamma^{\ep}_k)\\
	& = Q(\gamma^{\ep}_{k-1}) - Q(\gamma^{\ep}_k).
	\end{split}
	\end{equation}
	
	Now let us estimate the sizes of $\Delta^{\ep}_{k}$ and $\gamma^{\ep}_k$. First note that $D_i\lt(\sum_{j=1}^{4}\xgreen{[\gamma]_j}\zeta_j\rt)$, $i=1,2,3$, is a fixed quadratic function of $\gamma$ whose coefficients depend only on the function $a$, $\ti\alpha_2$, $s_0$ and $t_0$. Therefore, \bblue{for all $r>0$ and} $\gamma,\ti\gamma\in B_{r}(0)\subset\R^4$, we have
	\begin{equation}\label{eq2109}
	\|Q(\gamma)\| \leq C_1 \|\gamma\|^2
	\end{equation}
	and 
	\begin{equation}\label{eq2110}
	\|Q(\gamma)-Q(\ti\gamma)\| \leq \sup_{z\in B_{r}(0)} \|DQ(z)\|\cdot \|\gamma-\ti\gamma\| \leq C_1 r \|\gamma-\ti\gamma\|,
	\end{equation}
	where the above constant $C_1$ depends only on the coefficients of $Q$ and therefore does not depend on $\ep$ or $\gamma, \ti\gamma, r$. Let $\theta>0$ be sufficiently small such that 
	\begin{equation}\label{eq2117}
	\sum_{p=1}^{\infty} 2^{p-1}\theta^p \leq \frac{\lambda}{4\Lambda}<1.
	\end{equation}
	Clearly such $\theta$ exists. We denote 
	\begin{equation}
	\label{eq2117.3}
	\rred{C_2:=\lt\lVert A^{-1}\rt\rVert  C_1.} 
	\end{equation}
	Let $\ep_0:=\frac{\theta}{2\red{C_2}\Lambda}>0$. Now for all $0<\ep\leq \ep_0$, it follows from (\ref{eq2103}) that 
	\begin{equation}\label{eq2112}
	\red{C_2}\|\gamma^{\ep}_0\|\leq \red{C_2}\cdot 2\Lambda \ep_0=\theta.
	\end{equation}
	We claim that
	\begin{equation}\label{eq2111}
	\lt\lVert \Delta^{\ep}_{k}\rt\rVert \leq 2^{k-1}\theta^k\lt\lVert \gamma^{\ep}_{0}\rt\rVert
	\end{equation}
	and 
	\begin{equation}\label{eq2111.2}
	\lt\lVert \gamma^{\ep}_{k}\rt\rVert\leq \lt(1+\sum_{p=1}^{k} 2^{p-1}\theta^p\rt)\lt\lVert \gamma^{\ep}_{0}\rt\rVert < 2\lt\lVert \gamma^{\ep}_{0}\rt\rVert.
	\end{equation}
	We show this by induction. \red{Recall that $L^{\ep}(\gamma^{\ep}_{0})=0$}. We deduce from (\ref{eq2109}) that
	\begin{equation}
	\label{eq2117.4}
	\lt\lVert -G^{\ep}(\gamma^{\ep}_0)\rt\rVert \overset{\red{(\ref{eqdfg21})}}{=} \lt\lVert Q(\gamma^{\ep}_0)\rt\rVert \overset{\red{(\ref{eq2109})}}{\leq} C_1\lt\lVert \gamma^{\ep}_{0}\rt\rVert^2.
	\end{equation}
	It follows from this, (\ref{eq2116}) and  \ared{(\ref{eq2117.3})},  (\ref{eq2112}) that
	\begin{equation}
	\label{eq2117.7}
	\lt\lVert \Delta^{\ep}_{1}\rt\rVert \overset{\red{(\ref{eq2116}),(\ref{eq2117.4})}}{\leq} \lt\lVert A^{-1}\rt\rVert \cdot \red{C_1}\lt\lVert \gamma^{\ep}_{0}\rt\rVert^2 \overset{\ared{(\ref{eq2117.3})}}{=} \red{C_2} \lt\lVert \gamma^{\ep}_{0}\rt\rVert^2 \overset{\red{(\ref{eq2112})}}{\leq} \theta \lt\lVert \gamma^{\ep}_{0}\rt\rVert
	\end{equation}
	and therefore
	\begin{equation}
	\label{eq2117.6}
	\lt\lVert \gamma^{\ep}_{1}\rt\rVert \overset{(\ref{eq2116})}{\leq} (1+\theta)\lt\lVert \gamma^{\ep}_{0}\rt\rVert.
	\end{equation}
	So by (\ref{eq2117.7}), (\ref{eq2117.6}) we have that (\ref{eq2111}), (\ref{eq2111.2}) hold for $k=1$. Now suppose (\ref{eq2111}), (\ref{eq2111.2}) hold for $k\geq 1$. Using (\ref{eq2115}), (\ref{eq2110}) and the induction assumption, we have
	\begin{equation}
	\label{eq2114}
	\begin{split}
	&\lt\lVert G^{\ep}(\gamma_k^{\ep})\rt\rVert \overset{(\ref{eq2115})}{=} 
	\|Q(\gamma^{\ep}_{k-1})- Q(\gamma^{\ep}_{k})  \|\\
	&\qd\qd\qd\ared{\overset{(\ref{eq2110}),(\ref{eq2111.2})}{\leq} C_1\cdot 2\| \gamma^{\ep}_{0}\|
		\| \gamma^{\ep}_{k-1}- \gamma^{\ep}_{k}\|}\\
	&\qd\qd\qd\overset{ (\ref{eq2116})}{\leq} C_1\cdot 2\| \gamma^{\ep}_{0}\|\cdot \|\Delta^{\ep}_k \|
	\overset{(\ref{eq2111})}{\leq} C_1 2^k \theta^k \|\gamma_0^{\ep}\|^2 .
	\end{split}
	\end{equation}
	It follows from (\ref{eq2116}) and (\ref{eq2112}) that
	\begin{equation}
	\label{eq2114.4}
	\lt\lVert \Delta^{\ep}_{k+1}\rt\rVert \overset{\red{(\ref{eq2116})}}{\leq} \lt\lVert A^{-1}\rt\rVert \cdot \|G^{\ep}(\gamma_k^{\ep})\| \overset{\rred{(\ref{eq2114}), (\ref{eq2117.3})}}{\leq} C_2 2^{k}\theta^k\lt\lVert \gamma^{\ep}_{0}\rt\rVert^2  \overset{\red{(\ref{eq2112})}}{\leq} 2^{k}\theta^{k+1}\lt\lVert \gamma^{\ep}_{0}\rt\rVert
	\end{equation}
	and \red{
		\begin{equation*}
			\lt\lVert \gamma^{\ep}_{k+1}\rt\rVert \overset{\rred{(\ref{eq2116})}}{\leq}  \lt\lVert \gamma^{\ep}_{0} \rt\rVert+\sum_{p=1}^{k+1}  \lt\lVert \Delta^{\ep}_p   \rt\rVert \overset{\xgreen{(\ref{eq2111})}, (\ref{eq2114.4})}{\leq}\lt(1+\sum_{p=1}^{k+1} 2^{p-1}\theta^p\rt)\lt\lVert \gamma^{\ep}_{0}\rt\rVert \overset{(\ref{eq2117})}{\leq} 2 \lt\lVert \gamma^{\ep}_{0}\rt\rVert.
		\end{equation*}
		Thus we have established (\ref{eq2111}), (\ref{eq2111.2}) for general k}.
	
	Since $\{\gamma^{\ep}_k\}_{k}$ forms a bounded sequence, it has a convergent subsequence such that (without relabeling) 
	\begin{equation*}
		\lim_{k\rightarrow \infty} \gamma^{\ep}_k = \bar{\gamma}^{\ep}
	\end{equation*}
	for some $\bar{\gamma}^{\ep}$. We claim that $\bar{\gamma}^{\ep}$ is a non-negative solution to (\ref{eq2107}). From the estimates (\ref{eq2114}) and (\ref{eq2117}), we have
	\begin{equation*}
		\lt\lVert G^{\ep}(\gamma_k^{\ep})\rt\rVert \overset{(\ref{eq2114})}{\leq} C_1 2^{k}\theta^k\lt\lVert \gamma^{\ep}_{0}\rt\rVert^2 \rightarrow 0 \qd\text{ as }\qd k\rightarrow \infty.
	\end{equation*}
	Since $G^{\ep}$ is continuous, we have
	\begin{equation*}
		\|G^{\ep}(\bar{\gamma}^{\ep})\| = \lim_{k\rightarrow \infty} \|G^{\ep}(\gamma_k^{\ep})\| = 0.
	\end{equation*}
	It only remains to show that $\bar{\gamma}^{\ep}$ is non-negative componentwise. We deduce from (\ref{eq2111}), (\ref{eq2103}) and (\ref{eq2117}) that
	\begin{equation}\label{eq2118}
	\lt\lVert \gamma^{\ep}_k-\gamma^{\ep}_0\rt\rVert \overset{\red{(\ref{eq2116})}}{\leq} \sum_{p=1}^{k}\lt\lVert \Delta^{\ep}_p\rt\rVert \overset{\red{(\ref{eq2111})}}{\leq} \sum_{p=1}^{k}2^{p-1}\theta^p\lt\lVert \gamma^{\ep}_{0}\rt\rVert \overset{\rred{(\ref{eq2117}), (\ref{eq2103})}}{\leq} \frac{\lambda}{4\Lambda}\cdot 2\Lambda\ep = \frac{\lambda}{2}\ep.
	\end{equation}
	We know from (\ref{eq2103}) that each component of $\gamma^{\ep}_{0}$ is bounded below by $\lambda\ep$. This together with (\ref{eq2118}) shows that all components of $\gamma_k^{\ep}$ are bounded below by $\frac{\lambda}{2}\ep$ for all $k$. Therefore, the same holds for $\bar{\gamma}^\ep$. In particular, $\bar{\gamma}^\ep$ is non-negative.
\end{proof}

\begin{proof}[Proof of Theorem \ref{t103}]
	Given $0<\ep\leq \ep_0<1$, let $\bar{\gamma}^{\ep}=(\xgreen{[\bar{\gamma}^\ep]_1},\xgreen{[\bar{\gamma}^\ep]_2},\xgreen{[\bar{\gamma}^\ep]_3},\xgreen{[\bar{\gamma}^\ep]_4})$ be the non-negative solution \bblue{of (\ref{eq2107})} found in Lemma \ref{l106}. Then we have $\sum_{j=1}^{4}\xgreen{[\bar{\gamma}^\ep]_j}=\ep$. Define $\xgreen{[\bar{\gamma}^\ep]_0}:=1-\ep$. Then we have $\xgreen{[\bar{\gamma}^\ep]_j}\geq 0$ for all $j=0,1,2,3,4$ and $\sum_{j=0}^{4}\xgreen{[\bar{\gamma}^\ep]_j}=1$. Now we define
	\begin{equation*}
		\mu^{\ep}:=\sum_{j=0}^{4} \xgreen{[\bar{\gamma}^\ep]_j}\delta_{\zeta_j}.
	\end{equation*}
	It is clear that $\mu^{\ep}$ is a probability measure. Since $0<\ep<1$, $\mu^{\ep}$ is non-trivial. Since $\zeta_0$ is the trivial matrix and $\bar{\gamma}^{\ep}$ solves the system (\ref{eq2107}), we have
	\begin{equation*}
		\sum_{j=0}^{4}\xgreen{[\bar{\gamma}^\ep]_j}D_i(\zeta_j) = \sum_{j=1}^{4}\xgreen{[\bar{\gamma}^\ep]_j}D_i(\zeta_j) \overset{\red{(\ref{eq2107}), (\ref{eqdfg21})}}{=} D_i\lt(\sum_{j=1}^4 \xgreen{[\bar{\gamma}^\ep]_j}\zeta_j\rt) = D_i\lt(\sum_{j=0}^4 \xgreen{[\bar{\gamma}^\ep]_j}\zeta_j\rt)
	\end{equation*}
	for all $i=1,2,3$. This shows that $\mu^{\ep}\in \mathcal{M}^{pc}(\mathcal{K}_1^{\ti\alpha})$.
\end{proof}

\begin{proof}[Proof of Theorem \ref{t102} completed]
	We first consider the case $a'(\ti\alpha_2)>0$. \blue{Given $0<\ep\leq \ep_0<1$, let $\mu^{\ep}\in\mathcal{M}^{pc}(\mathcal{K}_1^{\ti\alpha})$ be the measure constructed in Theorem \ref{t103}. Let $\nu^{\ep}:=\lt((P_1^{\ti\alpha})^{-1}\rt)_{\sharp} \mu^{\ep}$. Note that since $P_1^{\ti\alpha}$ is a bijection, we have $\mu^{\ep}=(P_1^{\ti\alpha})_{\sharp} \nu^{\ep}$. Define $\ti{\mu}^{\ep}:=(P_1)_{\sharp} \nu^{\ep}$. Since $(P_1^{\ti\alpha})_{\sharp} \nu^{\ep} = \mu^{\ep}\in\mathcal{M}^{pc}(\mathcal{K}_1^{\ti\alpha})$, it follows from Lemma \ref{LAUX33} that $\ti{\mu}^{\ep}=(P_1)_{\sharp} \nu^{\ep}\in\mathcal{M}^{pc}(\mathcal{K}_1)$. Since $P_1$ and $P_1^{\ti\alpha}$ are both bijections, it is clear that $\ti{\mu}^{\ep}$ is also supported at five points, and hence is non-trivial. Further, by choosing $s_0, t_0$ sufficiently small in Theorem \ref{t103}, one can make the support of $\mu^{\ep}$ sufficiently small. It follows from Lemma \ref{LBBB3} that the support of $\ti{\mu}^{\ep}$ can be made sufficiently small.} This establishes the case where 
	$a'(\ti\alpha_2)>0$.
	
	\red{Now suppose $a'(\ti\alpha_2)<0$, then for some $\delta>0$ sufficiently small we have that
		\begin{equation}\label{eqddcc1}
		(v_2-v_1)(a(v_2)-a(v_1))<0 \text{ for any } v_1,v_2\in \lt(a(\ti\alpha_2)-\delta,a(\ti\alpha_2)+\delta\rt).
		\end{equation}
		Let \bblue{$\mathcal{K}_0:=\lt\{\lt(\begin{smallmatrix}  u & v \\ a(v) & u\end{smallmatrix}\rt):u,v\in \R \rt\}$}. Note that if $\det\lt(\begin{smallmatrix}  u_2-u_1 & v_2-v_1 \\ a(v_2)-a(v_1) & u_2-u_1\end{smallmatrix}  \rt)=0$ \bblue{for some $(u_1,v_1)$ and $(u_2,v_2)$ in $B_{\delta}(\xgreen{\ti{\alpha}})$}, then 
		\begin{equation*}
			(u_2-u_1)^2-(v_2-v_1)(a(v_2)-a(v_1))=0,
		\end{equation*}
		which by (\ref{eqddcc1}) implies $u_1=u_2$ and $v_1=v_2$. Thus, for sufficiently small neighborhood \bblue{$\ti U$} of $\xgreen{\lt(\begin{smallmatrix}  \ti \alpha_1 & \ti \alpha_2 \\ a(\ti \alpha_2) & \ti \alpha_1\end{smallmatrix}\rt)}$, $\mathcal{K}_0\cap \ti U$ does not contain Rank-$1$ connections \bblue{and \xgreen{therefore} $\det(X-Y)$ does not change sign on $(\mathcal{K}_0\cap \ti U)\times(\mathcal{K}_0\cap \ti U)$ \xgreen{by Lemma 1 in \cite{sv2}}}. 
		By \cite{sv2} \bblue{Lemma 3} we have that $\mathcal{M}^{pc}(\mathcal{K}_0\cap \ti U)$ consists of Dirac measures only. As $\mathcal{M}^{pc}(\mathcal{K}_1\cap U)$ can be embedded in  
		$\mathcal{M}^{pc}(\mathcal{K}_0\cap \ti U)$, this completes the proof of the case $a'(\ti\alpha_2)<0$, and hence the proof of Theorem \ref{t102}. } 
\end{proof}

\section{Appendix}
\label{AP}

In this appendix, we put together various auxiliary results used in the main body of the paper. 

\subsection{Auxiliary lemmas for Theorems \ref{T1} and \ref{C1}}
\label{APP}
\begin{a1}
\label{L6}
\pblue{Let $\FI=\lt\{f_1, f_2, \dots, f_{M_1}\rt\}$ be a collection of polynomials satisfying property $R$ (see Definition \ref{def2})}. For any $\varpi\in \R^n$, \pblue{let the translation $P^{\varpi}$ be defined by (\ref{eqbb35}).} Then we have
\begin{equation*}
\mu\in \Mb_{\FI}^{pc}(\varpi)\Longleftrightarrow \lt(P^{\varpi}\rt)_{\sharp} \mu\in \Mb_{\FI}^{pc}(0). 
\end{equation*}
\end{a1}

\begin{proof}
Note that since \pblue{$\FI$ satisfies property $R$}, for any 
$k\in \lt\{1,2,\dots, M_1\rt\}$ and any $z_0\in \R^n$ there exist $\alpha_0^{k,z_0}, \alpha_1^{k,z_0}, \dots, \alpha_{M_1}^{k,z_0}$ 
such that 
\begin{equation}
\label{eqp108}
f_k(z-z_0)=\sum_{i=1}^{M_1}  \alpha_i^{k,z_0} f_i(z)+\alpha_0^{k,z_0}.
\end{equation} 
Let $\mu\in \Mb_{\FI}^{pc}(\varpi)$. To simplify notation let $\ti{\mu}:=\lt(P^{\varpi}\rt)_{\sharp} \mu$. We have $\int_{\R^n} z d \ti{\mu}(z)=\int_{\R^n} (z-\varpi)\; d\mu (z)=0$. Further for any $k=1,2,\dots, M_1$ we have
\begin{equation}
\label{eqbb45}
\begin{split}
\int f_k(z) d\ti{\mu}(z)&=\int  f_k(z-\varpi)  d\mu(z)\\
&\overset{(\ref{eqp108})}{=}\int  \lt(\sum_{i=1}^{M_1}  \alpha_i^{k,\varpi} f_i(z)+\alpha_0^{k,\varpi}\rt)     d\mu(z)\\
&=\sum_{i=1}^{M_1}  \alpha_i^{k,\varpi} \int f_i(z)  d\mu(z)  +\alpha_0^{k,\varpi} \\
&\overset{\mu\in \Mb_{\FI}^{pc}(\varpi)}{=}\sum_{i=1}^{M_1}  \alpha_i^{k,\varpi}  f_i\lt(\varpi\rt)    +\alpha_0^{k,\varpi} \overset{(\ref{eqp108})}{=}f_k(0).
\end{split}
\end{equation}
Thus $\ti{\mu}\in \Mb_{\FI}^{pc}(0)$. If $\lt(P^{\varpi}\rt)_{\sharp} \mu\in \Mb_{\FI}^{pc}(0)$, in exactly the same way, an argument like (\ref{eqbb45}) gives that $\mu\in\Mb_{\FI}^{pc}(\varpi)$. 
\end{proof}

%
%

\em Notation. \rm For $i\in \lt\{1,2,\dots, m\rt\}$, $j\in \lt\{1,2,\dots, n\rt\}$ \pg{and $A\in M^{m\times n}$} let 
\begin{equation*}
\lt[A\rt]^{sb}_{i,j}\in M^{m-1,n-1} 
\end{equation*}
denote the matrix obtained from deleting the $i$-th row and the $j$-th column \pblue{of $A$}. Further let 
\begin{equation*}
\lt[A\rt]^{sb}_{0,j}\in M^{m,n-1}\text{ and }\lt[A\rt]^{sb}_{i,0}\in M^{m-1,n}
\end{equation*}
respectively denote the matrices obtained by deleting the $j$-th column and deleting the $i$-th row of $A$.

\begin{a1} 
\label{LA2}
Suppose $A,X\in M^{n\times n}$ and $r_0\in \lt\{1,2,\dots, n-1\rt\}$. \pblue{Let $I:=\lt\{i_1,i_2,\dots, i_{r_0}\rt\}$ and $J:=\lt\{j_1,j_2,\dots, j_{n-r_0}\rt\}$ be such that $I\cup J=\{1,2,\dots,n\}$}. Let $M^{I,J}_{A,X}\in M^{n\times n}$ denote the matrix whose first $r_0$ rows are given by 
\pblue{$R_{i_1}(A), \dots, R_{i_{r_0}}(A)$} and the remaining $n-r_0$ rows 
given by \pblue{$R_{j_1}(X), \dots, R_{j_{n-r_0}}(X)$}. Further let 
$\Xb^J\in M^{n-r_0,n}$ denote the matrix whose rows are given by 
$R_{j_1}(X), \dots, R_{j_{n-r_0}}(X)$. Then there exists $\Ib=\Ib(J)\subset \lt\{1,2,\dots, q(n-r_0,n)\rt\}$ such that 
\begin{equation}
\label{eqap2}
\det\lt(M^{I,J}_{A,X}\rt)=\sum_{l\in \Ib} \PI_l(A) M_l^{n-r_0, n}(\Xb^J)
\end{equation}
where $\lt\{\PI_l(A)\rt\}$ are polynomial functions \pblue{of} the entries of the matrix $A$.
\end{a1}

\begin{proof}
We prove this by induction on $n$. The lemma is immediate for 
$n=2$. Assume it is true for $n-1$. Let $A, X\in M^{n\times n}$, $r_0\in \lt\{1,2,\dots, n-1\rt\}$ and 
\pblue{$I:=\lt\{i_1,i_2,\dots, i_{r_0}\rt\}$, $J:=\lt\{j_1,j_2,\dots, j_{n-r_0}\rt\}$ be such that $I\cup J=\{1,2,\dots,n\}$}. 

To simplify notation let $B=M^{I,J}_{A,X}$. We expand $\det(B)$ along its first row. If $r_0=1$ this 
immediately gives the results because 
\begin{equation}
\label{eqap3}
\det\lt(B\rt)=\sum_{k=1}^{n} (-1)^{k+1} \pblue{[A]_{i_1k}} \det\lt(\lt[B \rt]^{sb}_{1,k}\rt)
\end{equation}
which is of the form (\ref{eqap2}). 

Now assume $r_0>1$. We apply the inductive hypothesis to the matrix $\lt[B \rt]^{sb}_{1,k}\in 
M^{n-1,n-1}$ for each $k=1,2,\dots, n$. So there exist $\Ib_k\subset \lt\{1,2, \dots, q(n-r_0,n-1)\rt\}$ and polynomials 
$\PI^k_1, \PI^k_2, \dots $ such that 
\begin{equation}
\label{eqap4}
\det\lt(\lt[B \rt]^{sb}_{1,k}\rt)=\sum_{l\in \Ib_k} \PI^k_l(A) M^{n-r_0, n-1}_l\lt( \lt[\Xb^J \rt]^{sb}_{0,k}  \rt).
\end{equation}
But there exists injective function 
\begin{equation*}
\varpi_k:\lt\{1,2,\dots, q(n-r_0,n-1)\rt\}\rightarrow \lt\{1,2,\dots, q(n-r_0,n)\rt\}
\end{equation*}
such that 
\begin{equation*}
M^{n-r_0,n-1}_l(\lt[A\rt]^{sb}_{0,k})=M^{n-r_0,n}_{\varpi_k(l)}(A)\text{ for all }A\in M^{n-r_0,n}.
\end{equation*}
So we can rewrite (\ref{eqap4}) as 
\begin{equation*}
\det\lt(\lt[B \rt]^{sb}_{1,k}\rt)=\sum_{l\in \Ib_k} \PI^k_l(A) M^{n-r_0, n}_{\pg{\varpi_k}(l)}\lt(\Xb^{J}\rt).
\end{equation*}
Putting this into \pblue{(\ref{eqap3})} we have that 
\begin{eqnarray}
\det(B)&=&\sum_{k=1}^n (-1)^{k+1} \pblue{[A]_{i_1k}}\lt( \sum_{l\in \Ib_k} \PI^k_l(A) M^{n-r_0, n}_{\pg{\varpi_k}(l)}\lt(\Xb^{J}\rt)   \rt)\nn\\
&=&\sum_{k=1}^n \sum_{l\in \Ib_k}(-1)^{k+1} \pblue{[A]_{i_1k}} \PI^k_l(A) M^{n-r_0, n}_{\pg{\varpi_k}(l)}\lt(\Xb^{J}\rt) \nn
\end{eqnarray}
which is of the form (\ref{eqap2}). 
\end{proof}

%
%

\begin{a1}
\label{LA3}
Given $A,X\in M^{n\times n}$, we have
\begin{equation}
\label{eqap12}
\det(A+X)=\det(A)+\det(X)+\sum_{k=1}^{q(n,n)} \PI_k(A) M^{n,n}_k(X)
\end{equation}
where $\PI_1, \PI_2, \dots, \PI_{q(n,n)}$ are polynomials functions of the entries of $A$. 
\end{a1}
\begin{proof} 
Note that \orgg{(recalling the notation (\ref{eqp2}))}
\begin{equation}
\label{eqap13}
\begin{split}
\det(A+X)&=\pblue{R_1(A+X)}\wedge R_2(A+X) \wedge \dots \wedge R_n(A+X)\\
&=\lt(R_1(A)+R_1(X)\rt)\wedge \lt(R_2(A)+R_2(X)\rt)\wedge \dots \wedge\lt(R_n(A)+R_n(X)\rt).
\end{split}
\end{equation}
Expanding this sum produces a number of terms, all but two of which are of the form 
\begin{equation}
\label{eqap14}
c R_{i_1}(A)\wedge \dots  \wedge R_{i_k}(A)\wedge  R_{j_1}(X)\wedge\dots\wedge  R_{j_{n-k}}(X)
\end{equation}
for some constant $c$, where $k\in \lt\{1,2,\dots, n-1\rt\}$, \pblue{$I:=\lt\{i_1,i_2,\dots, i_{k}\rt\}$ and $J:=\lt\{j_1,j_2,\dots, j_{n-k}\rt\}$ are such that $I\cup J=\{1,2,\dots,n\}$.} The only 
two terms in the expansion (\ref{eqap13}) that are not of the form (\ref{eqap14}) are $\det(X)$ and $\det(A)$. Now 
\begin{equation*}
R_{i_1}(A)\wedge \dots  \wedge R_{i_k}(A)\wedge  R_{j_1}(X)\wedge\dots\wedge  R_{j_{n-k}}(X)
=\det(M^{I,J}_{A,X})
\end{equation*}
and so by applying Lemma \ref{LA2} establishes (\ref{eqap12}). 
\end{proof}

\subsection{An improved \v{S}ver\'{a}k estimate for subspaces in $M^{3\times 3}_{sym}$}
\label{AP1}

\pblue{Following equation (4.9) in \cite{bhat}, we denote by $l(m,n)$ the maximum possible dimension of a linear subspace in $M^{m\times n}$ which contains no Rank-$1$ elements. As could be expected, the estimates on $l(m,n)$ can be improved in the case where the subspace is in $M^{3\times 3}_{sym}$ (symmetric $3\times 3$ matrices).} Note that the authors of \cite{bhat} state that Proposition 4.4 in the paper is \blue{due to} \v{S}ver\'{a}k. So following essentially exactly the arguments of Proposition 4.4 in
\cite{bhat} it is straightforward to obtain the the following estimate. 

\begin{a1}[\v{S}ver\'{a}k]
	\label{svlem}
	Let $K\subset M^{3\times 3}_{sym}$ be a \blue{three-dimensional} subspace, then $K$ must contain a Rank-$1$ element.
\end{a1}
\begin{proof}
	We argue by contradiction. Suppose $K\subset M^{3\times 3}_{sym}$ is a \blue{three-dimensional} subspace without Rank-$1$ connections. Note that $\mathrm{dim}\lt( M^{3\times 3}_{sym}\rt)=6$, and thus $\mathrm{dim}\lt( K^{\perp}\cap M^{3\times 3}_{sym}\rt)=3$. Let 
	$\{E_1, E_2, E_3\}$ be a basis of $K^{\perp}\cap M^{3\times 3}_{sym}$. For $a,b\in \R^3$ define 
	\begin{equation*}
	\lt[\Phi(a,b)\rt]_i:=a\cdot \lt(E_i b\rt)=E_{\pblue{i}} : (a\otimes b). 
	\end{equation*}
	So if $\Phi(a,b)=\lt(\begin{array}{c} 0 \\ 0 \\ 0 \end{array}\rt)$ then $a\otimes b\in K$. As we are assuming that $K$ has no Rank-$1$ connections, \pblue{this implies that $\Phi$ forms a non-singular bilinear mapping in the sense that if $\Phi(a,b)=0$ then 
		either $a=0$ or $b=0$}. \pblue{As noted in \cite{bhat} such mappings have been studied in the topological literature. The estimate we prove may indeed be 
		known in some form in those literature, however for the convenience of the reader we give a proof. }
	
	Note that for each $a\in \R^3\backslash \lt\{0\rt\}$ the mapping $x\mapsto \Phi(x,a)$ is linear and as such can be 
	represented by a matrix $M_a\in M^{3\times 3}$. Further as $\Phi$ is non-singular we have that 
	\begin{equation}
	\label{tseq362}
	\det(M_a)\not=0\text{ for all }a\in \R^3\backslash \lt\{0\rt\}. 
	\end{equation}
	Further note that by bilinearity we have that 
	\begin{equation*}
	x\mapsto \Phi(x,\lm_1 a_1+\lm_2 a_2)=M_{\lm_1 a_1+\lm_2 a_2} x=\lt(\lm_1 M_{a_1}+\lm_2 M_{a_2} \rt)x\text{ for all }x\in \R^3. 
	\end{equation*}
	Thus the mapping $\PI:\R^3\rightarrow M^{3\times 3}$ defined by $\PI(a):=M_a$ is a linear mapping and so is of the form 
	\begin{equation*}
	\PI(a)=\lt(\begin{array}{ccc} \varpi_{11}\cdot a &  \varpi_{12}\cdot a &  \varpi_{13}\cdot a \\
	\varpi_{21}\cdot a &  \varpi_{22}\cdot a &  \varpi_{23}\cdot a\\
	\varpi_{31}\cdot a &  \varpi_{32}\cdot a &  \varpi_{33}\cdot a
	\end{array}\rt). 
	\end{equation*}
	Thus $\det(\PI(a))$ is a $3$-homogeneous polynomial on $\R^3$. Let $a=(\alpha,\beta,\beta)$, then 
	\begin{equation*}
	\det(\PI(a))=c_0 \beta^3+c_1 \beta^2 \alpha+ c_2 \beta \alpha^2+c_3 \alpha^3\text{ for some }c_0, c_1, c_2, c_3\in \pblue{\R}. 
	\end{equation*}
	If $c_0=0$ then $\det\lt(\PI(0,1,1)\rt)=0$ which contradicts (\ref{tseq362}). Thus $c_0\not=0$. \blue{Similarly $c_3\ne 0$ as otherwise $\det\lt(\PI(1,0,0)\rt)=0$}. Then 
	\begin{equation*}
	\beta\mapsto \det\lt(\PI(1,\beta,\beta)\rt)=c_0\beta^3+c_1\beta^2+c_2 \beta+c_3
	\end{equation*}
	is of degree $3$ and so has a non-zero \blue{real} root $\beta_0$. Thus $\det\lt(\PI(1,\beta_0,\beta_0)\rt)=0$ which contradicts (\ref{tseq362}). This completes the proof of Lemma \ref{svlem}.
\end{proof}

\subsection{A counter example from \cite{bhat}}
\label{S11.2}
It was known already in \cite{bhat} that having no Rank-$1$ connections is in general not a sufficient condition for triviality of Null Lagrangian measures in subspaces in $M^{m\times n}$. In Proposition 4.2 in \cite{bhat}, a counter example of a four-dimensional subspace in $M^{3\times 3}$ is given. The subspace does not contain Rank-$1$ connections, yet it supports non-trivial Null Lagrangian measures. Fairly minor adaptions of \pgr{their example allow to construct counter examples for $\pgr{d}$-dimensional subspaces for all $\pgr{d}\geq 4$}.  Here for the convenience of the readers, we provide \pgr{the more general} counter example \pgr{focusing on} the adaptions needed.  First \pgr{we introduce some notation}. Given $A\in M^{n\times n}$ for $n\geq 3$, let $\Sbb(A)\in M^{3\times 3}$ denote the matrix defined by 
\begin{equation*}
\lt[\Sbb(A)\rt]_{ij}=\lt[A\rt]_{ij}\text{ for }i,j\in \lt\{1,2,3\rt\}.
\end{equation*}
\pgr{We denote by
\begin{equation*}
B(\alpha,\beta,\gamma,\delta):=\lt(\begin{array}{ccc} \beta+\delta & \alpha-\gamma & \gamma  \\
\alpha+\gamma & 0 & \delta  \\
\alpha & \beta & 0 \\
\end{array}\rt).
\end{equation*}
Given any non-negative integer $r$, let $P^r(\alpha,\beta,\gamma,\delta,\sigma_1,\dots,\sigma_r):\R^{4+r}\rightarrow M^{(3+2r)\times (3+2r)}$ be defined by
\begin{equation*}
\Sbb(P^r)=B(\alpha,\beta,\gamma,\delta),\, \lt[P^r\rt]_{2+2k, 2+2k}=\lt[P^r\rt]_{3+2k, 3+2k}=\sigma_{k} \text{ for } k=1,\dots,r
\end{equation*}
and all other entries of $P^r$ vanish. Further define
\begin{equation}
\label{cx1}
K^r:=\lt\{ P^r(\alpha,\beta,\gamma,\delta,\sigma_1,\dots,\sigma_r): \alpha,\beta,\gamma,\delta,\sigma_k\in\R\rt\}.
\end{equation}
Note that $K^0=\lt\{B(\alpha,\beta,\gamma,\delta):\alpha,\beta,\gamma,\delta\in\R\rt\}$ is exactly the subspace given  in Proposition 4.2 in \cite{bhat}. Then} we have
\begin{a5}[Bhattacharya-Firoozye-James-Kohn \cite{bhat}]\label{cx} 
\pgr{Given any non-negative interger $r$}, the \pgr{subspace} $K^r\subset M^{(3+2r)\times (3+2r)}$ defined in (\ref{cx1}) is a $(4+r)$-dimensional subspace, does not contain Rank-$1$ connections and $\mathcal{M}^{pc}(K^{\pgr{r}})$ is non-trivial.
\end{a5}
\begin{proof}
We first show that $\mathcal{M}^{pc}(K^{\pgr{r}})$ is non-trivial. \pgr{As in \cite{bhat}}, let $(\alpha_i,\beta_i,\gamma_i,\delta_i)\in\R^4$ for $i=1,2,3,4$ to be chosen later, and define 
$H_i\in M^{(3+2r)\times (3+2r)}$ to be \pgr{such that}
\begin{equation*}
\Sbb(H_i):= B(\alpha_i,\beta_i,\gamma_i,\delta_i)
\end{equation*}
\pgr{and all other entries of $H_i$ vanish}. Next we define
\begin{equation}
\label{cx2}
F_i=\begin{cases}
H_{\frac{i+1}{2}} & \text{if } i \text{ is odd},\\
-H_{\frac{i}{2}} & \text{if } i \text{ is even}.
\end{cases}
\end{equation}
We define the probability measure $\mu$ to be
\begin{equation}
\label{eqp300}
\mu := \sum_{i=1}^8 \frac{1}{8} \delta_{F_i}. 
\end{equation}
By (\ref{cx2}), it is clear that $\overline{\mu} = 0$. \pgr{So we need to show that $\int_{K^r}M_k(\zeta)d\mu = 0$ for all minors $M_k$ to conclude that $\mu\in\MI^{pc}(K^r)$. From (\ref{eqp300}) we have $\int_{K^r}M_k(\zeta)d\mu = \sum_{i=1}^8\frac{1}{8}M_k(F_i)$.} \org{Note that all the $F_i\in \lt\{P^r(\alpha,\beta,\gamma,\delta,0,\dots ,0):\alpha,\beta,\gamma,\delta\in \R\rt\}$. If $M_k$ is a minor for which 
$M_k(\pb{\zeta})$ involves elements \pb{$\lt[\pb{\zeta}\rt]_{lj}$ of the matrix $\zeta$} for \pb{some $l\geq 4$ or $j\geq 4$,} then $M_k(F_i)$ is the determinant 
of a submatrix of $F_i$ that contains at least one zero row or column,} \pgr{and thus $M_k(F_i)=0$ for all $i$. Hence we only have to check all minors in $K^0$. This is done in the proof of Proposition 4.2 in \cite{bhat} by choosing $(\alpha_i,\beta_i,\gamma_i,\delta_i)$ appropriately. In particular, thinking of $\alpha=(\alpha_1,\alpha_2,\alpha_3,\alpha_4)$ et al. as vectors in $\R^4$, one can choose $\alpha, \beta, \gamma, \delta$ to be unit vectors that are mutually perpendicular. Then it is straightforward to check that the measure $\mu$ defined in (\ref{eqp300}) commutes with all minors in $K^0$ and thus $\mu\in\MI^{pc}(K^r)$.}
	
	Finally we show that $K^{\pgr{r}}$ has no Rank-$1$ connections. \pgr{Suppose not, then $K^{r}$ would have a Rank-$1$ matrix $P(\alpha^0,\beta^0,\gamma^0,\delta^0,\sigma_1^0,\dots,\sigma_r^0)$. We must have $\sigma_k^0=0$ for all $k$, as otherwise the $(2+2k)$-th and $(3+2k)$-th rows (and columns) would be linearly independent. Thus the Rank-$1$ matrix is isomorphic to $B(\alpha^0,\beta^0,\gamma^0,\delta^0)$. However, in the proof of Proposition 4.2 in \cite{bhat}, it is shown that $K^0$ has no Rank-$1$ connections}, which is a contradiction. Hence $K^{\pgr{r}}$ has no Rank-$1$ connections.
\end{proof}

\subsection{An auxiliary lemma on \cblue{linear algebra}}
\label{grasman}

\begin{a1} 
	\label{grasaux1}
	Suppose $A\in M^{m\times n}$ and $m\leq n$. Then $\mathrm{Rank}(A)=m$ if and only if $\det\lt(A A^T \rt)\not =0$.
\end{a1}

\begin{proof} 
	By singular value decomposition $A=PBQ$ where $P\in O(m)$, $Q\in O(n)$ and $B$ is a diagonal matrix in $M^{m\times n}$. Let $\xgreen{\mathrm{Rank}(B)=p}$. 
	
	First assume $\mathrm{Rank}(A)=m$. \cblue{As $Q$ is 
		invertible, we have $\mathrm{Rank}(A Q^{-1})=m$} and it follows that $p=\mathrm{Rank}(PB)=m$. Now 
	\begin{equation}
	\label{graseq20}
	\det(A A^T)=\det\lt(P B Q Q^T B^T P^T\rt)=\det\lt(P B B^T P^T   \rt)=\det\lt(B B^T\rt)\not=0.
	\end{equation}
	
	Conversely if $\det(A A^T)\not=0$, \cblue{by calculations in (\ref{graseq20}) we have that $\det(B B^T)\ne 0$ and thus $p=m$. Therefore we have
		\begin{equation*}
			m= \mathrm{Rank}(PB)=\mathrm{Rank}(A Q^{-1})=\mathrm{Rank}(A).
		\end{equation*}
		This completes the proof of the lemma.}
\end{proof}

\end{document}